\documentclass[10pt,reqno]{amsart}
\usepackage{booktabs}
\usepackage{amssymb}
\usepackage{graphicx}
\usepackage{epstopdf} 
\usepackage[headings]{fullpage}
\usepackage{url} 
\usepackage{color}
\usepackage[usenames,dvipsnames,svgnames,table]{xcolor}
\usepackage{equipert}
\usepackage[square,numbers,sort&compress]{natbib}

\usepackage{tikz}
\usetikzlibrary{calc,positioning}

\def\visible<#1>{}  

\usepackage[utf8]{inputenc}
\usepackage{xr-hyper}

\usepackage[english]{babel}
\usepackage{amsfonts}
\usepackage{amsmath}
\usepackage{latexsym}
\usepackage{color}
\usepackage{subfigure}
\usepackage{enumerate}

\usepackage{hyperref}  

\usepackage{ifpdf}

\graphicspath{{figures/}}

\DeclareMathOperator    \aff                    {aff}

\DeclareMathOperator    \cl                     {cl}
\DeclareMathOperator    \conv           {conv}

\DeclareMathOperator    \intr                   {int}

\DeclareMathOperator    \proj           {proj}

\DeclareMathOperator    \relint         {rel\,int}

\DeclareMathOperator    \verts          {vert}

\DeclareMathOperator    \Aff {Aff}  

\newcommand{\old}[1]{{}}

\newcommand{\bb}{\mathbb}

\newcommand{\R}{\bb R}
\newcommand{\Q}{\bb Q}
\newcommand{\Z}{\bb Z}

\newcommand{\C}{\bb C}

\newcommand\st{\mid}
\newcommand\bigst{\mathrel{\big|}}

\def\ve#1{\mathchoice{\mbox{\boldmath$\displaystyle\bf#1$}}
{\mbox{\boldmath$\textstyle\bf#1$}}
{\mbox{\boldmath$\scriptstyle\bf#1$}}
{\mbox{\boldmath$\scriptscriptstyle\bf#1$}}}

\newcommand{\bpi}{\bar \pi}
\newcommand{\varphiD}{\psi_{q,\diag}}
\newcommand{\psiPoint}{\psi^m_{q,\point}}
\newcommand{\I}{\mathcal{P}}  
\newcommand{\Itri}[1][q]{\I_{#1,\tri}}
\newcommand{\Idiag}[1][q]{\I_{#1,\diag}}
\newcommand{\Ivert}[1][q]{\I_{#1,\ver}}
\newcommand{\Iverthor}[1][q]{\I_{#1,\ver\,\hor}}
\newcommand{\Ihor}[1][q]{\I_{#1,\hor}}
\newcommand{\Iedge}[1][q]{\I_{#1,\edge}}
\newcommand{\Ipoint}[1][q]{\I_{#1,\point}}
\newcommand{\Ipointedge}[1][q]{\I_{#1,\point\,\edge}}
\newcommand{\Ipointdiag}[1][q]{\I_{#1,\point\,\diag}}
\newcommand{\Ipointdiagtri}[1][q]{\I_{#1,\point\,\diag\,\tri}}
\renewcommand{\P}{\mathcal{P}}

\renewcommand{\S}{\mathcal{S}}
\newcommand{\Stri}{\S_{q,\tri}}
\newcommand{\barStri}{\bar\S_{q,\tri}}

\newcommand{\point}{{\EquiPoint}}
\newcommand{\ver}{{\EquiVertical}}
\newcommand{\hor}{{\EquiHorizontal}}
\newcommand{\diag}{{\EquiDiagonal}}
\newcommand{\edge}{{\EquiEdge}}

\newcommand{\tri}{{\EquiTriangle}} 
\newcommand{\triup}{{\EquiTriangleLower}}
\newcommand{\tridown}{{\EquiTriangleUpper}}

\newcommand\FundaTriangleLower{{}^{}_{\ve0}\EquiTriangleLower}  
\newcommand\FundaTriangleUpper{{}^{}_{\ve0}\EquiTriangleUpper}  

\newcommand{\E}{\mathcal{E}}
\newcommand{\G}{\mathcal{G}}

\newcommand\EqClass[1]{[#1]} 

\newcommand{\rx}{{\ve r}}
\newcommand{\x}{{\ve x}}
\newcommand{\y}{{\ve y}}
\newcommand{\z}{{\ve z}}
\renewcommand{\v}{{\ve v}}
\newcommand{\g}{{\ve g}}

\renewcommand{\u}{{\ve u}}
\renewcommand{\a}{{\ve a}}
\newcommand{\f}{{\ve f}}
\newcommand{\0}{{\ve 0}}

\newcommand{\p}{{\ve p}}
\renewcommand{\t}{{\ve t}}
\newcommand{\w}{{\ve w}}
\renewcommand{\b}{{\ve b}}

\newcommand{\cve}{{\ve c}}

\newcommand{\rr}{{\ve r}}

\def\st{\mid}

\newenvironment{psmallmatrixbig}{\bigl(\smallmatrix}{\endsmallmatrix\bigr)}
\newcommand\InlineFrac[2]{#1/#2}  
\newcommand\ColVec[3][\relax]
{
  \ifx#1\relax
  \bgroup\let\frac=\InlineFrac\begin{psmallmatrixbig}#2\vphantom{/}\\#3\vphantom{/}\end{psmallmatrixbig}\egroup
  \else
  \bgroup\let\frac=\InlineFrac\begin{psmallmatrixbig}\ifx#200\else#2/#1\fi\\\ifx#300\else#3/#1\fi\end{psmallmatrixbig}\egroup
  \fi
}

\newcommand\CVbaseline{11}
\newcommand\CVbaselineNew{25}

\newcommand\CVcoordleft[4]{#1}
\newcommand\CVcoordright[4]{#2}
\newcommand\CVcoordbottom[4]{#3}
\newcommand\CVcoordtop[4]{#4}

\newcommand{\diagGrid}[1][0101]{
\node[anchor=east,draw=none,font=\tiny,inner sep=2pt,yshift=1pt] at (0,0) {$\CVcoordbottom#1$};
\node[anchor=east,draw=none,font=\tiny,inner sep=2pt,yshift=-1pt] at (0,5) {$\CVcoordtop#1$};
\node[anchor=north,draw=none,font=\tiny,inner sep=2pt,xshift=1pt] at (0,0) {$\CVcoordleft#1$};
\node[anchor=north,draw=none,font=\tiny,inner sep=2pt,xshift=-1pt] at (5,0) {$\CVcoordright#1$};
\draw[step=1.0,black,thin] (0,0) grid (5,5);
\draw[black,thin]  (0,1) -- (1,0);
\draw[black,thin] (0,2) -- (2,0);
\draw[black,thin] (0,3)--(3,0);
\draw[black,thin] (0,4)--(4,0);
\draw[black,thin] (0,5)--(5,0);
\draw[black,thin]  (5,1) -- (1,5);
\draw[black,thin] (5,2) -- (2,5);
\draw[black,thin] (5,3)--(3,5);
\draw[black,thin] (5,4)--(4,5);}

\newcommand{\myCVone}[3][0101]{%
\begin{tikzpicture}[scale=0.2,baseline=\CVbaseline,sharp corners]
\diagGrid[#1]
\fill [red]  (#2,#3)  circle[radius=7pt];
\end{tikzpicture}}

\newcommand{\myCVtwo}[5][0101]{%
\begin{tikzpicture}[scale=0.2,baseline=\CVbaseline,sharp corners]
\diagGrid[#1]
\fill [red]  (#2,#3)  circle[radius=7pt];
\fill [red]  ((#4,#5)  circle[radius=7pt];
\draw[red,ultra thick] (#2,#3) -- (#4,#5); 
\end{tikzpicture}}

\newcommand{\myCVthree}[8][0101]{%
\begin{tikzpicture}[scale=0.2,baseline=\CVbaseline,sharp corners]
\draw[thick, top color=#8,bottom color=#8] (#2,#3) -- (#4,#5) -- (#6,#7) -- cycle; 
\diagGrid[#1]
\end{tikzpicture}}

\newcommand{\myCVfigSevenOne}[1][0101]{
\begin{tikzpicture}[scale=0.7, baseline=\CVbaselineNew]
\diagGrid[#1]
\node at (2.68, 0.7) {{\tiny $I_2$}};
\fill [red]  (2,1)  circle[radius=5pt];
\node at (0.35, 2.3) {{\tiny $I_1$}};
\fill [red]  (1,2)  circle[radius=5pt];
\node at (1.35, 1.3) {{\tiny $I$}};
\draw[red,ultra thick] (2,1) -- (1,2); 
\end{tikzpicture}}

\newcommand{\myCVfigSevenTwo}[1][0101]{
\begin{tikzpicture}[scale=0.7, baseline=\CVbaselineNew]
\node at (2.7, 0.7) {{\tiny $J$}};
\draw[thick, top color=red,bottom color=red] (2,0) -- (3,0) -- (2,1) -- cycle; 
\diagGrid[#1]
\end{tikzpicture}
}
\newcommand{\myCVfigSevenThree}[1][0101]{
\begin{tikzpicture}[scale=0.7, baseline=\CVbaselineNew]
\draw[thick, top color=pink,bottom color=pink, opacity=0.6] (3,2) -- (3,3) -- (5,1) -- (4,1) -- cycle; 
\draw[thick, top color=red,bottom color=red] (4,2) -- (4,1) -- (3,2) -- cycle; 
\diagGrid[#1]
\node at (2.68, 2.7) {{\tiny $K_1$}};
\node at (3.35, 1.3) {{\tiny $K_3$}};
\node at (4.7, 1.7) {{\tiny $K_2$}};
\end{tikzpicture}
}

\newtheorem{theorem}{Theorem}[section]

\makeatletter
\newcommand\MkNewTheorem[2]{%
  \newtheorem{#1}{#2}
  \expandafter\def\csname c@#1\endcsname{\c@theorem}
  \expandafter\def\csname p@#1\endcsname{\p@theorem}
  \expandafter\def\csname the#1\endcsname{\thetheorem}
  \expandafter\def\csname #1name\endcsname{#2}
}

\MkNewTheorem{corollary}{Corollary}
\MkNewTheorem{lemma}{Lemma}
\MkNewTheorem{proposition}{Proposition}
\MkNewTheorem{prop}{Proposition}
\MkNewTheorem{claim}{Claim}
\MkNewTheorem{observation}{Observation}
\MkNewTheorem{obs}{Observation}
\MkNewTheorem{conjecture}{Conjecture}
\MkNewTheorem{openquestion}{Open question}

\theoremstyle{definition}
\MkNewTheorem{example}{Example}
\MkNewTheorem{exercise}{Exercise}
\MkNewTheorem{notation}{Notation}
\MkNewTheorem{assumption}{Assumption}
\MkNewTheorem{definition}{Definition}
\MkNewTheorem{remark}{Remark}
\MkNewTheorem{goal}{Goal}

\makeatletter
\let\OurMathBbAux=\mathbb
\DeclareRobustCommand\OurMathBb{\OurMathBbAux}
\let\mathbb=\OurMathBb
\let\bfseries=\undefined
\DeclareRobustCommand\bfseries
{\not@math@alphabet\bfseries\mathbf
  \boldmath\fontseries\bfdefault\selectfont\let\OurMathBbAux=\mathbf}
\def\@thm#1#2#3{%
  \ifhmode\unskip\unskip\par\fi
  \normalfont
  \trivlist
  \let\thmheadnl\relax
  \let\thm@swap\@gobble
  \thm@notefont{\fontseries\mddefault\upshape\unboldmath}
  \thm@headpunct{.}
  \thm@headsep 5\p@ plus\p@ minus\p@\relax
  \thm@space@setup
  #1
  \@topsep \thm@preskip               
  \@topsepadd \thm@postskip           
  \def\@tempa{#2}\ifx\@empty\@tempa
    \def\@tempa{\@oparg{\@begintheorem{#3}{}}[]}%
  \else
    \refstepcounter{#2}%
    \def\@tempa{\@oparg{\@begintheorem{#3}{\csname the#2\endcsname}}[]}%
  \fi
  \@tempa
}
\makeatother

\makeatletter
\renewcommand{\pod}[1]
{\allowbreak\mathchoice{\mkern18mu}{\mkern8mu}{\mkern8mu}{\mkern8mu}(#1)}
\makeatother

\title[Equivariant Perturbation III]{Equivariant Perturbation in \\Gomory and Johnson's Infinite Group
  Problem\\[1ex] III. Foundations for the $k$-Dimensional Case\\ with  Applications to $k=2$
  }

\thanks{An extended abstract with some of the results of the paper has
    appeared as: \emph{Equivariant perturbation in Gomory and Johnson's
      infinite group problem. II. The unimodular two-dimensional case} in:
    Michel Goemans and Jos\'e Correa (eds.), Integer Programming and
    Combinatorial Optimization, Lecture Notes in Computer Science, vol.~7801,
    Springer, 2013, pp.~62--73, ISBN 978-3-642-36693-2.}
  \thanks{The authors gratefully acknowledge partial support from the National Science
    Foundation through grants DMS-0636297 (R.~Hildebrand), DMS-0914873
    (R.~Hildebrand, M.~K\"oppe), and DMS-1320051 (M.~K\"oppe). 
  }

\author{Amitabh Basu} 
\address{Amitabh Basu: Dept.~of Applied Mathematics and Statistics, The Johns Hopkins University}
\email{basu.amitabh@jhu.edu}

\author{Robert Hildebrand}
\address{Robert Hildebrand: Institute for Operations Research, Dept. of Mathematics, ETH Z\"urich, Switzerland}
\email{robert.hildebrand@ifor.math.ethz.ch}

\author{Matthias K\"oppe}
\address{Matthias K\"oppe: Dept.\ of Mathematics, University of California, Davis}
\email{mkoeppe@math.ucdavis.edu}

\date{$\relax$Revision: 2090 $ - \ $Date: 2016-07-28 13:56:14 -0700 (Thu, 28 Jul 2016) $ $\!\!\!}

\usepackage[normalem]{ulem}

\begin{document}
 \newcommand{\tgreen}[1]{\textsf{\textcolor {ForestGreen} {#1}}}
 \newcommand{\tred}[1]{\texttt{\textcolor {red} {#1}}}
 \newcommand{\tblue}[1]{\textcolor {blue} {#1}}

\setlength\marginparwidth{0.7in}
\newcommand\citedinsurveyas[1]{}
\newcommand\forwardreferencetosectionthree{\marginpar{\raggedright Forward
    reference to section 3.}}

\maketitle
\begin{abstract}
  We develop foundational tools for classifying the extreme valid
  functions for the $k$-dimensional infinite group problem.  In particular,
  we
  present the general regular solution to Cauchy's additive functional
  equation on restricted lower-dimensional convex domains.  This provides a $k$-dimensional
  generalization of the so-called Interval Lemma, allowing us to deduce affine
  properties of the function from certain additivity relations.   Next, we study the discrete
  geometry of additivity domains of piecewise linear functions,
  providing a framework for finite tests of minimality and extremality.
   We then give a theory of non-extremality certificates in the form of 
  perturbation functions. 

  We apply these tools in the context of minimal valid functions for the
  two-dimensional infinite group problem that are piecewise linear on a
  standard triangulation of the plane, under a regularity
  condition called diagonal constrainedness.  
  We show that the extremality of a minimal valid function is equivalent 
  to the extremality of its restriction to a certain finite two-dimensional
  group problem.  This gives an algorithm for testing the
  extremality of a given minimal valid function.
\end{abstract}

\clearpage
\tableofcontents
\clearpage

\section{Introduction}
\label{sec:introduction}
Over 40 years ago,  Gomory and
Johnson introduced an elegant infinite-dimensional relaxation of integer linear
optimization problems called the \emph{infinite
  group problem}~\cite{infinite,infinite2}.
The motivation for studying it is the hope to find
effective multi-row cutting plane procedures with better performance
characteristics compared to the single-row cutting plane procedures in use
today. 

\subsection{The group problem}

Gomory's \emph{group problem}~\cite{gom} is a central object in the study of strong
cutting planes for integer linear optimization problems.  One considers an
abelian group $G$, written additively, and studies
the set of functions $s \colon G \to \R$ satisfying the following constraints:  
\begin{equation}
  \label{GP} 
  \begin{aligned}
    &\sum_{\rx \in G} \rx\, s(\rx) \in \ve f + S \\
    &s(\rx) \in \mathbb{Z}_+ \ \ \textrm{for all $\rx \in G$}  \\
    &s \textrm{ has finite support}, 
  \end{aligned}
\end{equation}
where $S$ is a subgroup of $G$ and $\ve f$ is a given element in $G\setminus S$; so $\ve
f + S$ is the coset containing the element $\ve f$.  
We will be concerned with the so-called {\em infinite group
  problem} \cite{infinite,infinite2}, where $G =\R^k$ is taken to be the group
of real $k$-vectors under addition, and $S= \Z^k$ is the subgroup of the
integer vectors. 
We are interested in studying the convex hull $R_{\ve f}(G,S)$ of all
functions satisfying the constraints in~\eqref{GP}. Observe that $R_{\ve f}(G,S)$ is
a convex subset of the infinite-dimensional vector space
$\mathcal{V}$ of functions $s \colon G \to \R$ with finite support.

 A main focus of the research in this area is to  give a description
 of $R_{\ve f}(\R,\Z)$ as the intersection of halfspaces of $\mathcal{V}$.  This
 makes a very useful connection between $R_{\ve f}(\R, \Z)$ and traditional integer
 programming, both from a theoretical, as
 well as, practical point of view. This
 arises from the fact that important classes
 of cutting planes for general integer
 programs can be viewed as finite-dimensional
 restrictions of the linear inequalities used
 to describe $R_{\ve f}(\R, \Z)$.

 \subsection{Valid inequalities and valid functions}
 Any linear inequality in
$\mathcal{V}$ is given by $\sum_{\rx \in G} \pi(\rx)s(\rx) \geq
\alpha$ where $\pi$ is a function
$\pi\colon G \to \R$ and $\alpha \in
\R$. The left-hand side of the inequality is a finite sum because $s$ has finite
support. Such an inequality is called a {\em valid inequality} for $R_{\ve f}(G,S)$
if $\sum_{\rx \in G} \pi(\rx)s(\rx) \geq \alpha$ for all $s \in R_{\ve f}(G,S)$.  It is
customary to concentrate on valid inequalities with $\pi \geq 0$;
then we can choose, after a scaling, $\alpha = 1$. Thus, we only focus on
valid inequalities of the form $\sum_{\rx \in G} \pi(\rx)s(\rx) \geq 1$ with $\pi
\geq 0$. Such functions $\pi$ will be termed {\em valid functions} for
$R_{\ve f}(G,S)$. 

 As pointed out in \cite{corner_survey}, the nonnegativity
             assumption in the definition of a valid function might seem
             artificial at first. Although there exist valid
             inequalities $\sum_{r \in \R} \pi(r)s(r) \geq \alpha$ for
             $R_{\ve f}(\R,\Z)$ such that $\pi(r) < 0$ for some~$r \in \R$, it can
             be shown that $\pi$ must be nonnegative over all
             \emph{rational} $r \in \Q$. Since data in integer programs is
             usually rational, it is natural to focus on nonnegative valid
             functions.  

 \subsection{Minimal functions} 
 Gomory and
 Johnson~\cite{infinite,infinite2} defined a hierarchy on the set of valid
 functions, capturing the strength of the corresponding valid inequalities, which we summarize now.
 
A valid function $\pi$ for $R_{\ve f}(G,S)$ is said to be
\emph{minimal} for $R_{\ve f}(G,S)$ if there is no valid function $\pi' \neq \pi$
such that $\pi'(\rx) \le \pi(\rx)$ for all $\rx \in G$.  For every valid
function $\pi$ for $R_{\ve f}(G,S)$, there exists a minimal valid function $\pi'$
such that $\pi' \leq \pi$ (cf.~\cite{bhkm}), and thus non-minimal valid
functions are redundant in the description of $R_{\ve f}(G,S)$.  Minimal
functions for $R_{\ve f}(G,S)$ were characterized by Gomory for the case where
$S$ has finite index in~$G$ in~\cite{gom}, and later for $R_{\ve f}(\R,\Z)$ by Gomory and
Johnson~\cite{infinite}. We state these results in a unified notation in the
following theorem. 

A function $\pi\colon G \rightarrow \mathbb{R}$ is \emph{subadditive} if
$\pi(\x + \y) \le \pi(\x) + \pi(\y)$ for all $\x,\y \in G$. We say that  $\pi$ is
\emph{symmetric} if $\pi(\x) + \pi(\f - \x) = 1$ for all $\x \in G$.
     
\begin{theorem}[Gomory and Johnson \cite{infinite}] \label{thm:minimal} Let
  $\pi \colon G \rightarrow \mathbb{R}$ be a nonnegative function. Then $\pi$
  is a minimal valid function for $R_{\ve f}(G,S)$ if and only if $\pi(\ve z) = 0$ for
  all $\ve z\in S$, $\pi$ is subadditive, and $\pi$ satisfies the symmetry
  condition. (The first two conditions imply that $\pi$ is periodic modulo
  $S$, that is, $\pi(\x) = \pi(\x + \ve z)$ for all $\ve z \in S$.) 
\end{theorem}

 \begin{remark} 
   Note that this implies that one can view a minimal
   valid function $\pi$ as a function from $G/S$ to $\R$, and thus
   studying $R_{\ve f}(G,S)$ is the same as studying $R_{\ve f}(G/S,\0)$.  However, we avoid
   this viewpoint in this paper.
 \end{remark}


\subsection{Extreme functions and their classification}

In polyhedral combinatorics, one is interested in classifying the
facet-defining inequalities of a polytope, which are the strongest
inequalities and provide a finite minimal description.
In the infinite group problem, the analogous notion is that of an {\em extreme
  function}. 

A~valid function~$\pi$ is \emph{extreme}
for $R_{\ve f}(G,S)$ if it cannot be written as a convex combination of two other
valid functions for $R_{\ve f}(G,S)$, i.e., $\pi = \tfrac12(\pi^1 + \pi^2)$ implies $\pi = \pi^1 = \pi^2$.  Extreme functions are
minimal.
 
Various sufficient conditions for extremality have been
proved in the previous literature~\cite{bhkm,3slope,dey1,dey2,dey3,tspace,Richard-Li-Miller-2009:Approximate-Liftings,Miller-Li-Richard2008}.  
In part I~\cite{basu-hildebrand-koeppe:equivariant} of the present series of papers, 
the authors initiated the study 
of perturbation functions
that are equivariant with 
respect to certain finitely generated reflection groups. This addressed an
inherent previously unknown \emph{arithmetic} (number-theoretic) aspect of  
the problem and allowed the authors to give an algorithm that tests extremality of
piecewise linear functions with rational breakpoints and relate extremality to a finite-dimensional problem.
\begin{theorem}[Theorems 1.3 and 1.5 in~\cite{basu-hildebrand-koeppe:equivariant}]
\label{thm:one-dim}
Consider the following problem.  
\begin{quote}
  Given a minimal 
  valid function $\pi$ for $R_f(\R,\Z)$ that is piecewise
  linear with a set of rational breakpoints with the least common
  denominator~$q$, decide if $\pi$ is extreme or not.
\end{quote}
\begin{enumerate}[\rm(i)]
\item There exists an algorithm for this problem that takes a number of
  elementary operations over the reals that is bounded by a polynomial in
  $q$. 
\item If the function~$\pi$ is continuous, then $\pi$ is extreme for
  $R_f(\R,\Z)$ if and only if the restriction $\pi \big|_{\frac1{4q}\Z}$ is extreme
  for the finite group problem $R_f(\frac1{4q}\Z,\Z)$. 
\end{enumerate}
\end{theorem}

\medbreak

\subsection{Contributions, techniques, and outline of this paper}

\begin{figure}[tp]
  \centering
\bgroup
\begin{tikzpicture} [ font = \small, align = flush center, >=stealth, thick,
  node distance = 4cm, 
  every node/.style={draw, rounded corners, text width=3.1cm, align=center, minimum height=1.3cm}]
  \node [] (sec1-introduction) 
  {1. Introduction};
  \node [fill=orange!40, below left of = sec1-introduction, text width=4cm] (sec2-cauchy) 
  {2. {Regular solutions to Pexider's functional equation on bounded domains of $\R^k$}};
  \node [fill=orange!40, below right of = sec1-introduction, text width=4cm] (sec3-complexes) 
  {3. Discrete geometry of piecewise linear minimal valid functions and their additivity domains};
  \node [fill=orange!20, right = 1cm of sec3-complexes.east] (appb-genuine) 
  {B. Genuinely $k$-dimensional functions};
  \node [fill=blue!20, below = 5cm of sec1-introduction] (sec4-class) 
  {4. {A class of minimal valid functions defined over $\R^2$}};
  \node [fill=blue!20, below = 1cm of sec4-class] (sec5-proofs) 
  {5. Proof of the main results for the two-dimensional case};
  \node [fill=blue!10, right = 3.5cm of sec5-proofs, text width=4cm] (appa2-perturbations) 
  {A.2 {Deriving the perturbation functions 
      using equivariance formulas}};
  \node [fill=orange!20, above = 1.5cm of appa2-perturbations, text width=4cm] (appa1-reflectiongroups) 
  {A.1 {Reflection groups and their fundamental domains}};
  \draw[->] (sec1-introduction) -- (sec2-cauchy);
  \draw[->] (sec2-cauchy) -- (sec4-class);
  \draw[->] (sec4-class) -- (sec5-proofs);
  \draw[->] (sec1-introduction) -- (sec3-complexes);
  \draw[->] (sec3-complexes) -- (sec4-class);
  \draw[->] (appa1-reflectiongroups) -- (appa2-perturbations);
  \draw[dotted] (sec5-proofs) -- (appa2-perturbations);
  \draw[dotted] (sec3-complexes) -- (appb-genuine);
  \coordinate (x) at ($(sec2-cauchy)!0.5!(sec4-class) - (4cm,0)$) ;
  \coordinate (x2) at ($(x) + (7.5cm,0)$) ;
  \coordinate (y) at ($(appa1-reflectiongroups)!0.5!(appa2-perturbations) + (3cm,0)$) ;
  \coordinate (y2) at ($(y) - (5cm,0)$) ;
  \draw[dashed] (x) --  (x2) -- (y2) -- (y);
  \node [draw=none,font=\itshape] at ($(sec4-class)!0.5!(sec5-proofs) - (3.7cm,0cm)$) {
    application to the two-dimensional case};
  \node [draw=none,font=\itshape,above=1cm of appb-genuine] {
    general $k$-dimensional theory};
\end{tikzpicture}
\egroup

  \caption{Structure of the paper. 
    Sections 2 and 3 can be read independently.  
    Readers who are already familiar with the Interval Lemma may find it
    convenient to start with section~2, whereas readers familiar with
    polyhedral complexes may want to start with section~3.} 
  \label{fig:paper-structure-poset}
\end{figure}
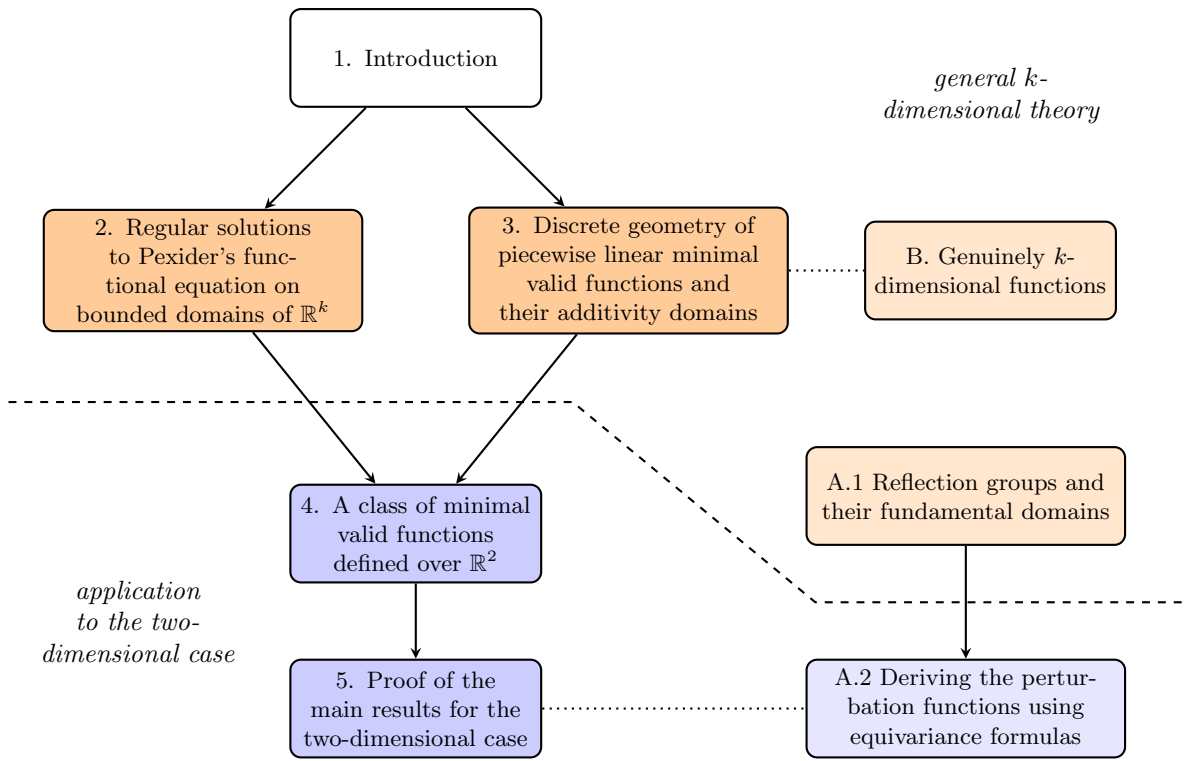

In the present paper, we continue the program
of~\cite{basu-hildebrand-koeppe:equivariant} of algorithmically studying the extemality of piecewise linear functions.  We prove several general
results that hold for arbitrary dimension~$k$ and then apply them to give an algorithm that tests 
the extremality of a large class of functions for the case~$k=2$.
The structure of the paper is outlined in \autoref{fig:paper-structure-poset}.

The main technique used to show a function $\pi$ is extreme is to assume that $\pi = \tfrac{1}{2}(\pi^1 + \pi^2)$ where $\pi^1, \pi^2$ are valid functions, and then show that $\pi = \pi^1 = \pi^2$.  
We will use three important properties of $\pi^1, \pi^2$ in our proofs, which are summarized in the following lemma.  These facts for the one-dimensional case can be found, for instance, in \cite{basu-hildebrand-koeppe:equivariant}, and are easily extended to the general $k$-dimensional case.
\begin{lemma}\label{lem:tightness}\label{lem:minimality-of-pi1-pi2}\label{lemma:tight-implies-tight}\label{Theorem:functionContinuous} \label{lem:lipschitz}
  \citedinsurveyas{Lemma 1.4}
  Let $\pi$ be minimal, $\pi = \frac12(\pi^1+\pi^2)$, and $\pi^1,\pi^2$
  valid functions.  Then the following hold:
  \begin{enumerate}[(i)]
  \item $\pi^1,\pi^2$ are minimal.
  \item All subadditivity relations $\pi(\x + \y) \le \pi(\x) + \pi(\y)$  
  that are tight for~$\pi$ are also tight for
  $\pi^1,\pi^2$.  That is, defining the \emph{additivity domain} of~$\pi$ as 
\begin{equation}
\label{eq:Epi}
  E(\pi) := \{\,(\x, \y) \st \Delta \pi(\x, \y) := \pi(\x) + \pi(\y) -\pi(\x + \y) = 0\,\},
\end{equation}
  we have $E(\pi) \subseteq E(\pi^1), E(\pi^2)$.
  \item If $\pi$ is continuous and piecewise linear, then $\pi, \pi^1, \pi^2$ are all Lipschitz continuous.
  \end{enumerate}
\end{lemma}

\subsubsection{Functional equations.}
Utilizing the set $E(\pi)$ is fundamental in the literature to classifying extreme functions.  In particular, much of the literature relies on a bounded version of a result for the classical (additive) Cauchy functional equation
\begin{equation}\label{eq:cauchy}
  \theta(u)+\theta(v) = \theta(u+v),
\end{equation}
where $u,v\in \R$ (see, e.g.,
\cite{aczel66,kuczma78,kuczma09,czerwik02,kannappan09}).  This result is known 
as the \emph{Interval Lemma} in the integer programming
community~\cite{tspace}. 
\begin{lemma}[Interval lemma~\cite{tspace,bccz08222222}] \label{lem:interval_lemma} Let $\theta \colon \R \rightarrow
\mathbb{R}$ be a function bounded on every bounded interval. Given real numbers $u_1 < u_2$ and $v_1 < v_2$, let $U = [u_1, u_2]$, $V =
[v_1, v_2]$, and $U + V = [u_1 + v_1, u_2 + v_2]$.
If $\theta(u)+\theta(v) = \theta(u+v)$ for every
$(u,v)  \in U\times V$, then $\theta$ is affine with the slope $c\in \R$ 
in each of the intervals $U$, $V$, and $U+V$.
\end{lemma}
The Interval Lemma gives a powerful dimension reduction mechanism: where it
applies, the infinite-dimensional space of functions on an interval is
replaced by a finite-dimensional space.  If this applies to all subintervals
of a piecewise linear function, testing if this function is extreme can be reduced to finite-dimensional
linear algebra.  

For the $k$-dimensional case, various authors in the integer programming community have given suitable
generalizations of this lemma~\cite{3slope,dey3,bhkm}.  There is also a parallel line of work in the functional equations literature, e.g.,~\cite{baker-rado-1987,kuczma09, kannappan09}.  In the present paper, we state and prove a certain version of these results which allows for additivity relations to hold on lower dimensional domains. To the best of our knowledge, this lower dimensional variant of such functional equations is new. We treat directly the so-called \emph{Pexider equation}, which is a simple generalization that allows for three functions instead of one that is well-studied in the functional equations community, but not as much in the integer programming community. This generalization comes at no cost in the proofs.  The utility of considering it in this generality will become apparent in a following paper~\cite{basu-hildebrand-koeppe:equivariant-general-2dim}.  

While the novelty of this paper is the lower dimensional variant of the Pexider equations proved in Theorems~\ref{lem:generalized_interval_lemma} and~\ref{lem:projection_interval_lemma}, for the expository purposes of this introduction we state two consequences whose statements are cleaner. Nonetheless, these next two results are extremely useful for understanding extremality, in our opinion.


\begin{theorem}[Higher-dimensional Interval Lemma, full-dimensional
  version]\label{theorem:generalized_interval_lemma_fulldim} 
\citedinsurveyas{Theorem 1.6}
Let $f,g,h \colon \R^k \to \R$ be bounded functions. Let $U$ and $V$ be convex
subsets of $\R^k$ such that $f(\u) + g(\v) = h(\u+\v)$ for all $(\u, \v) \in U\times V$. Assume that $\aff(U) = \aff(V) = \R^k$. Then there exists a vector $\cve\in \R^k$ such that
%
%
%
$f$, $g$ and $h$ are affine over $U$, $V$ and $W = U+V$, respectively, with the
same gradient $\cve$.
\end{theorem}

The key generalization is to consider an additivity domain specified by a general convex set $F\subseteq \R^k \times \R^k$ instead the more restrictive setting of $F = U \times V$.  

Define the projections $p_1,p_2,p_3\colon \R^k\times \R^k \to \R^k$ as
\begin{equation}
\label{eq:projections}
p_1(\x,\y) = \x, \quad p_2(\x,\y) = \y, \quad  p_3(\x,\y) = \x+\y.  
\end{equation}

\begin{theorem}[Convex additivity domain lemma, full-dimensional
  version]\label{lem:projection_interval_lemma_fulldim} 
  \citedinsurveyas{Theorem 1.7}
  Let $f,g,h \colon \R^k \to \R$ be bounded functions. 
  Let $F \subseteq \R^k \times \R^k$ be a full-dimensional convex set
  such that $f(\u) + g(\v) = h(\u+\v)$ for all $(\u, \v) \in F$. 
  Then there exists a vector $\cve\in \R^k$ such that $f, g$ and $h$ are
  affine with the same gradient $\cve$ over $\intr(p_1(F))$,
  $\intr(p_2(F))$ and $\intr(p_3(F))$, respectively. 
\end{theorem}

While Theorems~\ref{theorem:generalized_interval_lemma_fulldim} and~\ref{lem:projection_interval_lemma_fulldim} are simple corollaries of our Theorems~\ref{lem:generalized_interval_lemma} and~\ref{lem:projection_interval_lemma}, we mention here that they also follow immediately from the main result of~\cite{baker-rado-1987} (see Theorem~\ref{thm:baker-rado} for a statement of the result from~\cite{baker-rado-1987}).
It is notable that we can only deduce affine linear properties over the
\emph{interiors} of the projections.  This is best possible, as we illustrate
by examples (Remark~\ref{remark:projection_interval_lemma_interiors_only} and Remark~\ref{remark:projection_interval_lemma_interiors_only2}).


\subsubsection{Piecewise linear functions and the discrete geometry of their additivity domains.}
Piecewise linear functions form an important class of minimal valid functions.
In fact, all classes of extreme functions described in the literature are piecewise
linear, with the exception of a family of measurable functions constructed
in~\cite{bccz08222222}.  

\begin{figure}[t]
\begin{center}
\input{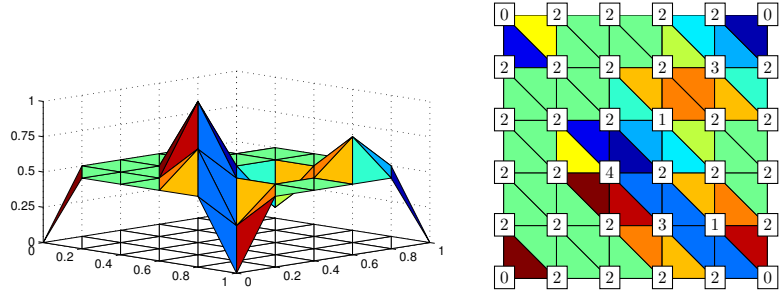}
\end{center}
\caption{A minimal valid, continuous, piecewise linear function over the
  triangulation $\P_5$, which is diagonally constrained
  .
  \emph{Left}, the three-dimensional plot of the function on 
  the unit square.
  \emph{Right},
  the triangulation $\P_5$, restricted to the unit square 
  and colored according to slopes to match the 3-dimensional plot, and labeled
  with values $v$ at each vertex of $\P_5$ where the function takes value
  $\tfrac{v}{4}$. 
}
\label{figure:diagonallyConstrained-new-figure}
\end{figure}
In the one-dimensional case ($k=1$), a continuous piecewise linear
function~$\pi$ periodic modulo~$\Z$ is given by a list of breakpoints in
$[0,1]$ and affine functions on the subintervals delimited by these
breakpoints.  If the value of $\pi$ is known on the breakpoints, then $\pi$ is
already uniquely defined everywhere by linear interpolation.

In the higher-dimensional case
($k>1$), it is not enough to give a list of breakpoints; 
rather, one needs a triangulation.  
As our prime example for $k=2$, consider the function shown
in~\autoref{figure:diagonallyConstrained-new-figure}. 
Its pieces are defined on the lower and upper triangles
$$\FundaTriangleLower = \tfrac1q \conv(\{ \ColVec{0}{0}
, \ColVec{1}{0}
, \ColVec{0}{1}
 \}) \qquad\text{and}\qquad \FundaTriangleUpper = \tfrac1q \conv(\{\ColVec{1}{0}
, \ColVec{0}{1}
,
\ColVec{1}{1}
\})$$ (with $q=5$) and their translates by elements of the lattice
$\smash[t]{\frac1q\Z^2}$.
Together these triangles form a well-known\footnote{For example,  in the context of homotopy methods
\cite{Forster:Homotopy-Methods}, this triangulation is known as 
the K1 triangulation.} triangulation of the space~$\R^2$,
which is, of course, periodic modulo~$\Z^2$.  
It has convenient geometric and arithmetic properties 
and will play an important role in the present paper; we denote it by $\P_q$. 


In general we describe piecewise linear functions $\pi \colon \R^k \to \R$ by specifying
a polyhedral complex $\P$ (a collection of polyhedra, meeting face-to-face;
see \autoref{s:prelim}) 
that covers all of $\R^k$ and affine functions on
the cells of this complex.  
The use of polyhedral complexes generalizes that of triangulations.

%
%
\smallbreak
 
 \subsubsection{Combinatorial representation of additivity domain through $\P$}
Our second main contribution in the present paper is a detailed study of the
discrete geometry of the additivity domain $E(\pi)$, as defined in~\eqref{eq:Epi}, 
of a function~$\pi$ that is continuous piecewise linear
on a polyhedral complex~$\P$.  This is missing from the previous literature on $R_\f(\R^k, \Z^k)$ for $k\geq 2$ and extends the discussion in the one-dimensional
case in \cite{basu-hildebrand-koeppe:equivariant}. In section~\ref{section:delta-p-definition}, we show that the subadditivity slack function $\Delta \pi$ (as defined in~\eqref{eq:Epi}) is continuous piecewise linear over a polyhedral complex in $\R^k\times \R^k$ that we call $\Delta \P$.  Therefore, $E(\pi)$ is composed of faces of $\Delta \P$ on which the piecewise linear function $\Delta \pi$ is constantly zero, which can be determined completely by the values of $\Delta \pi$ at the vertices of $\Delta \P$.  It follows that the vertices of $\Delta \P$ hold information for necessary and sufficient
conditions for minimality, as shown in Theorem~\ref{minimality-check}. 
The faces of $\Delta \P$ that are contained in $E(\pi)$ are referred to as
\emph{additive faces} and are partially ordered by set inclusion.  The
inclusion-maximal faces are called the \emph{maximal additive faces}. In
section~\ref{section:additivity-discretized}, we show that these maximal additive faces provide a combinatorial description of $E(\pi)$ as the union of certain
polytopes (Lemma~\ref{lemma:covered-by-maximal-valid-triples}). This proves to be a crucial ingredient to show that a piecewise linear function is not extreme.


Further, minimal functions can be classified according to the types of maximal additive faces~$F\in\Delta\P$ that appear. The \emph{generic case} is that in which all maximal additive faces, with the possible exception of those
corresponding to the symmetry condition, are full-dimensional
in~$\R^k\times\R^k$. In this case, the Interval Lemma (for $k=1$) or the
full-dimensional version of the Higher-dimensional Interval Lemma (for
$k\geq2$) are sufficient for proving extremality.  All sufficient conditions
for extremality studied in the previous literature fall into this class. 
\emph{Degenerate cases}, in which some maximal additive faces are 
lower-dimensional, require more machinery.  

In \cite{basu-hildebrand-koeppe:equivariant}, the authors examine functions for $k=1$ with rational breakpoints.  They interpret lower-dimensional maximal additive faces as translation and
reflection operations on the real line.
Using the structure of these
operations, a special class of ``perturbation'' functions is introduced in
\cite{basu-hildebrand-koeppe:equivariant}, which are used as certificates for
the non-extremality of a given minimal function. Understanding the nature of
these lower-dimensional maximal additive faces and their interaction with
these perturbation functions was the key to breaking beyond the existing
arguments from the literature which dealt with only full-dimensional maximal
additive faces for the $k=1$ case. 

For higher dimensions, degenerations of various types are possible and define
a hierarchy of functions. Just like the situation in the $k=1$ case
suggests, as one climbs up in this hierarchy, the extremality
proofs become more and more complex.  In this paper, we initiate this
higher-dimensional theory by studying the $k=2$ case, for piecewise linear
functions over a special triangulation of
$\R^2$ and a particular type of degeneration only. 



\subsubsection{Characterization of extreme piecewise linear functions on a
  standard triangulation of the plane.}


In the present paper, we restrict ourselves to functions on the triangulation $\P_q$
that have a particular type of degeneration of the maximal additive faces
only.  These functions are called \emph{diagonally constrained} functions; the 
definition appears in \autoref{sec:triangulation}.  (The example function
shown in \autoref{figure:diagonallyConstrained-new-figure} is a
diagonally constrained function.)

In the following two theorems, we require that $\f \in \verts(\P_q)$.  This turns out to be a natural assumption because for minimal functions that cannot be viewed as a lower-dimensional function, we must always have $\f \in \verts(\P)$, \autoref{thm:faVertex-gen-k-dim}.  Such functions are called \emph{genuinely $k$-dimensional} and were studied in \cite{bhkm,3slope}.  We detail properties of these functions in \autoref{sec:gen-k-functions}. In particular, we show that the study of continuous piecewise linear extreme functions can, under some mild assumptions, be reduced to the study of genuinely $k$-dimensional functions that are continuous and piecewise linear. 

\begin{theorem}\label{thm:main}
\citedinsurveyas{Theorem 1.8}
Consider the following problem.  
\begin{quote}
  Given a minimal 
  valid function $\pi$ for $R_{\f}(\R^2,\Z^2)$ that is 
  piecewise linear continuous on~$\P_q$ and diagonally constrained with $\f \in \verts(\P_q)$, 
  decide if $\pi$ is extreme.
\end{quote}
There exists an algorithm for this problem that takes a number of elementary operations over the reals that is
bounded by a polynomial in $q$.
\end{theorem}

As a direct corollary of the proof of Theorem~\ref{thm:main}, we obtain the following result relating
the finite and infinite group problems.  

\begin{theorem}\label{thm:1/4q}
\citedinsurveyas{Theorem 1.9}
Let $\pi$ be a minimal continuous piecewise linear function over $\P_q$ that is diagonally constrained and $\f \in \verts(\P_q)$.  Fix $m \in \Z_{\geq 3}$.
Then $\pi$ is extreme for $R_{\ve f}(\R^2, \Z^2)$ if and only if the restriction $\pi\big|_{\tfrac{1}{mq}\Z^2}$ is extreme for $R_{\ve f}(\frac{1}{mq} \Z^2, \Z^2)$.  
\end{theorem}

The two main developments for the proof of Theorem~\ref{thm:main} and Theorem~\ref{thm:1/4q} are an understanding of how additivities combine to imply piecewise linear conditions, such as Theorem~\ref{lem:projection_interval_lemma_fulldim}, and how perturbation functions can imply a function is not extreme.  Specific perturbation functions are described in section~\ref{s:eq-perturb}.  In Appendix~\ref{s:reflection-groups}, we give a more abstract discussion of how perturbation functions can be understood through reflection groups.  The proof of Theorem~\ref{thm:main} and Theorem~\ref{thm:1/4q} is completed in section~\ref{sec:connection-to-finite-group}.

\section{Regular solutions to Cauchy's functional equation on restricted domains of $\R^k$}
\label{s:real-analysis}

\subsection{Cauchy's and Pexider's functional equations}

As mentioned in the introduction, the standard technique for showing
extremality of a minimal valid function~$\pi\colon\R^k\to\R$ is as follows.  Suppose that
$\pi = \frac12(\pi^1 + \pi^2)$, where $\pi^1,\pi^2$ are other (minimal)
valid functions.  One then studies the \emph{additivity domain} $E(\pi)$.
By Lemma~\ref{lemma:tight-implies-tight}, $E(\pi) \subseteq E(\pi^1),
E(\pi^2)$.  One then considers $\pi$, $\pi^1$, $\pi^2$ as solutions to the
\emph{functional equation} 
\begin{equation}\label{eq:cauchy}
  \theta(\u)+\theta(\v) = \theta(\u+\v),\quad (\u,\v) \in F,
\end{equation}
where $F = E(\pi)$.  

This equation is known as the \emph{(additive) Cauchy functional equation}.
Classically (see, e.g., \cite{kuczma09,czerwik02}), it is studied for
functions~$\theta\colon\R^k\to\R$, when the additivity domain~$F$ is the entire
space~$\R^k\times\R^k$.
The solutions to~\eqref{eq:cauchy} with $F=\R^k\times\R^k$ are referred to as
\emph{additive functions}.  
The obvious solutions to~\eqref{eq:cauchy}, 
namely 
the (homogeneous) linear functions $\theta(\ve x) = \ve c\cdot \ve x$, are referred to as the \emph{regular solutions}. In addition, there exist certain
pathological solutions, which are highly discontinuous. In order to rule out
these solutions, one imposes a regularity hypothesis.  Various such regularity
hypotheses have been proposed in the literature.  For example, it is sufficient
to assume that the function~$\theta$ is bounded on bounded intervals, 
or continuous at a point, or bounded below on a finite interval, or locally
Lebesgue integrable;  see~\cite[Theorem 1.2]{kannappan09} for a list of many
more equivalent conditions. Under each of these conditions, one deduces that
the additive function $\theta\colon\R^k\to\R$  is continuous and hence a
(homogeneous) linear function \cite[Theorems 1.1 and 1.2]{kannappan09}. 

A natural and commonly studied generalization of the Cauchy functional equation is the \emph{Pexider equation}
\begin{equation}
\label{eq:pexider}
f(\ve u) + g(\ve v) = h(\ve u + \ve v), \quad (\ve u,\ve v)\in F,
\end{equation}
where $f,g,h \colon \R^k \to \R$.  When $F = \R^k \times \R^k$, it is easily
shown that the solutions to the Pexider equation are $f(\x) = \theta(\x) +
\alpha$, $g(\y) = \theta(\y) + \beta$, $h(\z) = \theta(\z) + \alpha + \beta$
for some additive function $\theta$
satisfying~\eqref{eq:cauchy}~\cite{kannappan09}.   Hence, this equation on the
entire domain reduces to studying the Cauchy functional equation.
Combining this with a regularity condition, we find that 
the regular solutions are affine functions; so we
lose homogeneity of the solutions.

\subsection{Restricted additivity domains}
The additivity domain $E(\pi)$ of a subadditive function~$\pi\colon\R^k \to\R$ can
be a complicated set.  It is convenient to break it into convex sets $F$, which
we then study independently.  

When $F \subsetneq \R^k \times \R^k$, equations~\eqref{eq:cauchy}
and~\eqref{eq:pexider} are referred to as conditional Cauchy and Pexider
equations or as Cauchy and Pexider equations on restricted
domains~\cite{kuczma78, dhombres-ger-78}.  It is clear that the Pexider
equation imposes no conditions on the function values of $f$, $g$, and $h$
outside of the projections $p_1(F)$, $p_2(F)$, and $p_3(F)$, respectively, where the projections are as defined in~\eqref{eq:projections}. 
Baker and Rad{\'o} \cite{baker-rado-1987} show that when the
Pexider equation is satisfied on a restricted open path-connected domain, then the solutions
on each of the projections are constant shifts of the same additive function~$\theta\colon \R^k\to\R$.   We provide a slightly modified version of \cite[Corollary 1]{baker-rado-1987} that removes one assumption.  



\begin{theorem}[\cite{baker-rado-1987}]
\label{thm:baker-rado}
Let $F \subseteq \R^{k} \times \R^{k}$  non-empty, path-connected, and open.  Let $f,g,h\colon \R^k \to \R$ such that~\eqref{eq:pexider} holds for all $(\x,\y) \in F$.  Then there exist an additive function $\theta\colon \R^k \to \R$ and constants $\alpha, \beta\in \R$  such that $f(\x) = \theta(\x) + \alpha$, $g(\y) = \theta(\y) + \beta$, and $h(\z) = \theta(\z) + \alpha + \beta$ for all $\x\in p_1(F)$, $\y \in p_2(F)$, and $\z \in p_3(F)$.

Furthermore, let $D \subseteq \R^{k} \times \R^{k}$ such that $F \subseteq D \subseteq \cl F$, where $\cl F$ denotes the
closure of~$F$.  Suppose that $D$ satisfies the following:
for every $\x \in p_1(D)$ there exist $\y \in p_2(F)$, $\z \in p_3(F)$ such that $\x +\y = \z$ and for every $\y \in p_2(D)$ there exist $\x \in p_1(F)$, $\z \in p_3(F)$ such that $\x + \y = \z$.  
Then $f(\x) = \theta(\x) + \alpha$, $g(\y) = \theta(\y) + \beta$, and $h(\z) = \theta(\z) + \alpha + \beta$ for all $\x\in p_1(D)$, $\y \in p_2(D)$, and $\z \in p_3(D)$.
\end{theorem}
\begin{proof}
By~\cite[Theorem 1]{baker-rado-1987},  there exist an additive function $\theta\colon \R^k \to \R$ and constants $\alpha, \beta\in \R$  such that $f(\x) = \theta(\x) + \alpha$, $g(\y) = \theta(\y) + \beta$, and $h(\z) = \theta(\z) + \alpha + \beta$ for all $\x\in p_1(F)$, $\y \in p_2(F)$, and $\z \in p_3(F)$.  Let $\x \in p_1(D)$, $\y \in p_2(F)$, $\z \in p_3(F)$ such that $\x + \y = \z$.  Then $f(\x) = h(\z) - g(\y) = \theta(\z) - \theta(\y) + \alpha = \theta(\x) + \alpha$.  Similarly, for any $\y \in p_2(D)$, $g(\y) = \theta(\y) + \beta$.  Finally, for any $\z \in p_3(D)$, there exists a preimage $(\x,\y) \in D$ such that $\z = \x + \y$.  Then $h(\z) = f(\x) + g(\y) = \theta(\x) + \theta(\y) + \alpha + \beta = \theta(\z) + \alpha + \beta$.  
\end{proof}

Hence, combined with a regularity condition, affine properties of the
functions on the projections can be deduced.

\subsection{Interval lemma in $\R^1$}
The so-called Interval
Lemma was introduced by Gomory and Johnson to the integer programming
community in~\cite{tspace}.\footnote{Similar results were known independently in the functional equations community.  For
  instance,~\cite{aczel66} states the result for $U = V$.  
 }  
It concerns the Cauchy functional
equation~\eqref{eq:cauchy} on a restricted additivity domain $F$ that is a
rectangle $F=U\times V$, where $U$ and $V$ are bounded intervals.  Then
$p_1(F) = U$, $p_2(F) = V$, and $p_3(F) = U+V$, a Minkowski sum.
We present it here as a corollary of
\autoref{thm:baker-rado}, together with regularity conditions.


The following lemma is stated with the regularity assumption that $f,g,h$ are
bounded functions; but this assumption can be replaced by any of the other
regularity assumptions discussed above.  

\begin{lemma}[Interval lemma]\label{one-dim-interval_lemma}
\citedinsurveyas{Lemma 2.2}
Given real numbers $u_1 < u_2$ and $v_1 < v_2$, let $U = [u_1, u_2]$, $V = [v_1, v_2]$, and $U + V = [u_1 + v_1, u_2 + v_2]$.
Let $f \colon U \rightarrow
\mathbb{R}$, $g \colon  V \rightarrow
\mathbb{R}$, $h \colon U+V  \rightarrow
\mathbb{R}$ be bounded functions. \\ If $f(u)+g(v) = h(u+v)$ for every
$(u, v) \in U \times V$, then there exists $c\in \R$ such that $f(u)=f(u_1)+c(u-u_1)$ for every $u\in U$, $g(v)=g(v_1)+c(v-v_1)$ for every $v\in V$, $h(w)=h(u_1+v_1)+c(w-u_1-v_1)$ for every $w\in U+V$. In other words, $f$, $g$ and $h$ are affine with slope $c$ over $U$, $V$, and $U+V$ respectively.
\end{lemma}
\begin{proof}
Consider the rectangle $D = U \times V \subseteq \R^2$ and let $F = \intr(D)$.  Since $U$ and $V$ are proper intervals, for every $x \in U$, there exists a $y \in \intr(V)$ such that $x + y \in \intr(U+V) = p_3(F)$.  Similarly, for every $y \in V$, there exists a $x \in \intr(U)$ such that $x +y \in p_3(F)$.  Therefore, by~\autoref{thm:baker-rado}, there exists an additive function $\theta \colon \R \to \R$ and constants $\alpha, \beta \in \R$ such that $f(x) = \theta(x) + \alpha$, $g(y) = \theta(y) + \beta$, $h(z) = \theta(z) + \alpha + \beta$ for all $x \in U, y \in V, z \in p_3(D) = U+V$.  

Since $f$ is bounded on $U$, $\theta$ is bounded on $U$.  Therefore, by~\cite[Theorems 1.1 and 1.2]{kannappan09}, $\theta(x) = c x$ for some $c \in \R$.  This competes the proof.
\end{proof}


\subsection{Higher-dimensional Interval Lemma}
\label{s:generalized-interval-lemma}

The generalization of the Interval Lemma for hypercubes $U = V = [a,b]^k$ was stated in~\cite{aczel66}.
The only known generalizations of \autoref{one-dim-interval_lemma} in the integer programming community 
literature appear in \cite{dey3, 3slope} for the case of $k=2$ and in \cite{bhkm} for general $k$.  The results in \cite{3slope, bhkm} are special cases of our Theorem~\ref{lem:generalized_interval_lemma} that require one of the sets to intersect the origin.  The result in \cite{dey3} applies in $k=2$ and allows for so-called \emph{star-shaped} sets that also contain the origin; a similar proof to our generalization also yields a result on star-shaped sets, but we avoid this direction because we do not need this type of result.

In fact, the proof of~\autoref{one-dim-interval_lemma} easily generalizes to the $k$-dimensional setting to prove~\autoref{theorem:generalized_interval_lemma_fulldim}.  This is because the result of Baker and Rad{\'o} (\autoref{thm:baker-rado}) also applies in the $k$-dimensional setting.  Then the problem reduces to $k$ separate one-dimensional problems since any additive function $\theta\colon \R^k \to \R$ can be can be decomposed into $k$ univariate additive functions~\cite[Theorem 1.24]{kannappan09}.  This is under the assumption that the domains of $U$, $V$ of $f$, $g$ are full-dimensional and the additivity domain is the full Cartesian product $U \times V$.  

\begin{figure}[t!]
 \begin{center}
  \includegraphics[scale=0.4]{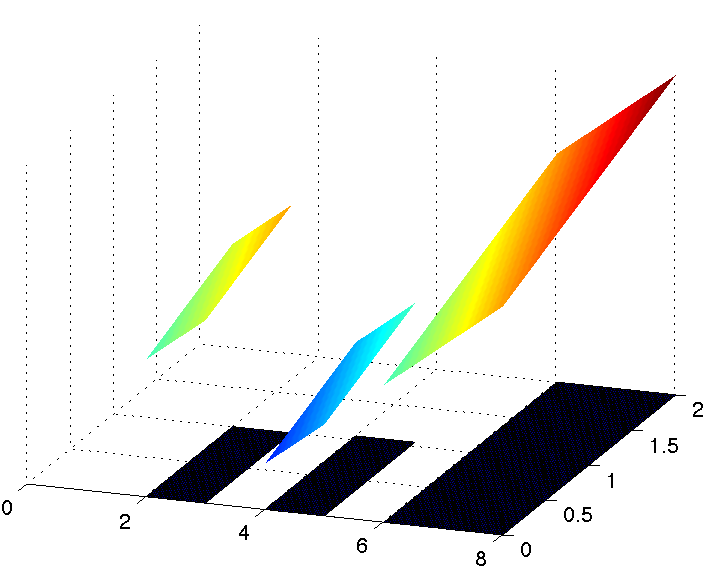}
   \includegraphics[scale=0.4]{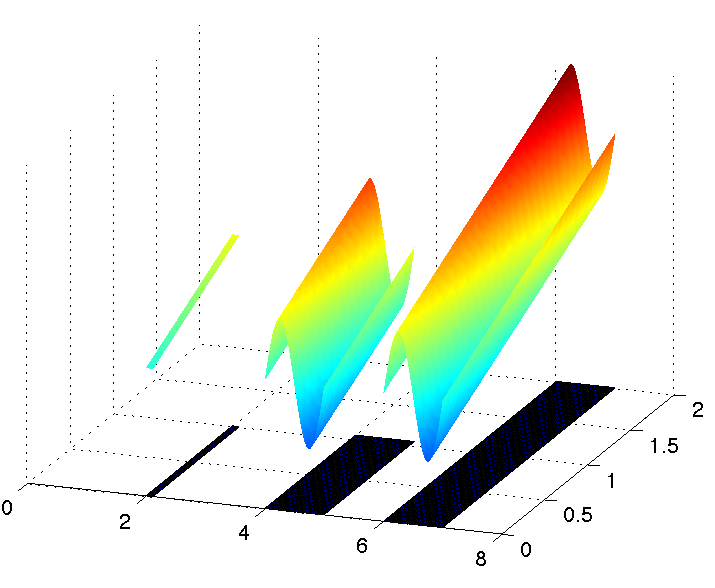}
  \includegraphics[scale=0.4]{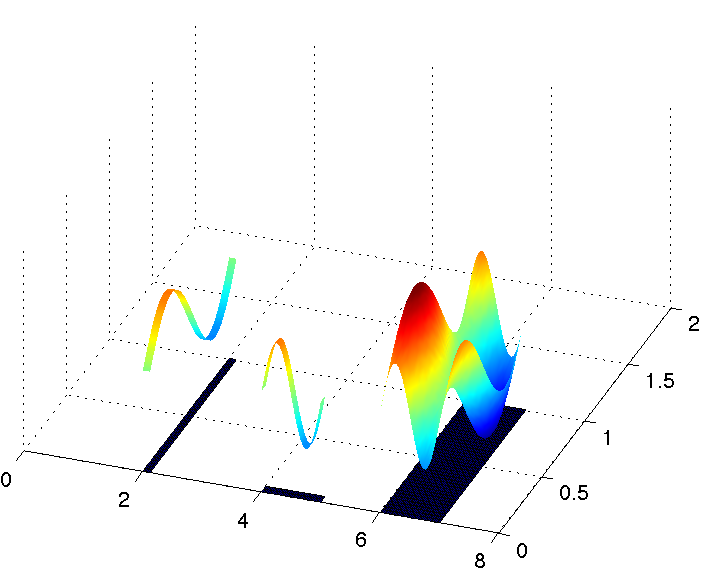} 
     \\
      \end{center}
      \vspace{-0.5cm}
      \begin{flushleft}
   \hspace{2.5cm}   (a) \hspace{4.3cm} (b) \hspace{4.3cm} (c) \hspace{5cm}
      \end{flushleft}
  \caption{Cauchy's functional equation on bounded additivity domains $F = U\times
    V$.  Each diagram shows $p_1(F) = U$ (\emph{left black shadow}), $p_2(F)=V$
    (\emph{middle black shadow}), and $p_3(F) = U + V$ (\emph{right black
      shadow}), and the graph (\emph{colored by function values}) of an
    example function that is additive with respect to this domain. 
    (a)
    Full-dimensional situation. (b) Sum of a one-dimensional and a
    two-dimensional set; not a direct sum.  
   (c) Direct sum of (non-parallel) one-dimensional sets.}
  \label{fig:higher-dim-interval-lemma}
\end{figure}

We prove the result in a more general setting, in which the
additivity domain is $U\times V$ for convex sets $U\subseteq\R^k$ and
$V\subseteq\R^k$, which are not necessarily of the same dimension.  In this
general setting we cannot expect to deduce that the solutions are affine over
$U$, $V$, and $U+V$.  In particular, these results will differ from most literature since the domain of additivity is not full-dimensional.
\begin{remark}
  Indeed, if $U + V$ is a direct sum,
  i.e., for every $\w \in U + V$ there is
  a unique pair $\u \in U$, $\v \in V$ with $\w = \u+\v$, then 
  $f(\u) + g(\v) = h(\u+\v)$ merely expresses a form of separability of~$h$
  with respect to certain subspaces,
  and $f$ and $g$ can be arbitrary functions; 
  see \autoref{fig:higher-dim-interval-lemma}\,(c). 
\end{remark}

\begin{definition}
Let $U \subseteq \R^k$. 
Given a linear subspace $L \subseteq \R^k$, we say $\pi\colon U \to \R$ is
{\em affine with respect to $L$ over $U$} if there exists $\cve \in \R^k$ such
that $\pi(\u^2) - \pi(\u^1) =  \cve\cdot( \u^2 - \u^1)$ for any $\u^1, \u^2
\in U$ such that $\u^2 - \u^1 \in L$.
\end{definition}

\begin{theorem}[Higher-dimensional Interval Lemma]\label{lem:generalized_interval_lemma}
\citedinsurveyas{Theorem 2.5}
Let $f,g,h \colon \R^k \to \R$ be bounded functions. Let $U$ and $V$ be convex subsets of $\R^k$ such that $f(\u) + g(\v) = h(\u+\v)$ for all $(\u,\v) \in F = U\times V$. Let $L$ be a linear subspace of $\R^k$ such that $(L + U) \times (L + V) = (L\times L) + F \subseteq \aff(F) = \aff(U) \times \aff(V)$. Then there exists a vector $\cve\in \R^k$ such that
%
%
%
$f$, $g$ and $h$ are affine with respect to $L$  over $p_1(F) = U$, $p_2(F) = V$ and $p_3(F) = U+V$
respectively, with gradient $\cve$.
\end{theorem}

For the proof, we will only use the machinery of~\autoref{one-dim-interval_lemma}.  We note that certain elements of the proof could also be done using~\autoref{thm:baker-rado}, but there does not seem to be a direct implication.

We will need the following notation and basic result.
For any element $\x \in\R^k$, $k \geq 1$, $\lvert \x \rvert_\infty$ will denote the
standard $\ell^\infty$ norm. We use $B^\infty(\u, r)$ to denote the open
$\ell^\infty$ ball around $\u\in \R^k$ with radius $r\in \R_+$, i.e.,
$B^\infty(\u, r) = \{\, \x \in \R^k \st | \u - \x |_\infty < r \,\}$. 

\begin{lemma}\label{lem:relint-linspace}
Let $U\subseteq \R^k$ be a convex set and let $L$ be a linear space such that $L + U \subseteq \aff(U)$. Then, for any $\u \in \relint(U)$, there exists $r > 0$ such that $B^\infty(\u, r) \cap (\u + L) \subseteq U$.
\end{lemma}
\begin{proof}
It suffices to show that for any $\p \in L$ there exists $\epsilon > 0$ such that $\u + \epsilon\p \in U$. One then can use a basis of $L$ to find the desired $r > 0$. 

Since $L + U \subseteq \aff(U)$, $L$ is a subspace of $\aff(U) - \u$. Thus,
$\p \in \aff(U) - \u$ and therefore, $\u + \p \in \aff(U)$. Since $U$ is
convex and $\u \in \relint(U)$, there exists $\epsilon >0$ such that $\u +
\epsilon\p \in U$.
\end{proof}

\begin{proof}[Proof of Theorem~\ref{lem:generalized_interval_lemma}]
If $m := \dim(L) = 0$, there is nothing to prove. So we assume $m \geq 1$ and let $\p^1, \ldots, \p^m$ be a basis for $L$ (we obviously have $m \leq k$). Since $U$ is convex and $L + U \subseteq \aff(U)$, by Lemma~\ref{lem:relint-linspace} for any vector $\u^0 \in \relint(U)$, there exist real numbers $u^i_1 <0 < u^i_2$ such that the set $U_0 := \{\u^0 + \sum_{i=1}^m \lambda_i\p^i \st  u^i_1 \leq \lambda_i \leq u^i_2 \;\forall i =1, \ldots, m\} \subseteq U$. Similarly, for any vector $\v^0 \in \relint(V)$, there exist real numbers $v^i_1 <0 <  v^i_2$ such that the set $V_0 := \{\v^0 + \sum_{i=1}^m \mu_i\p^i \st  v^i_1 \leq \mu_i \leq v^i_2  \;\forall i =1, \ldots, m\} \subseteq V$.

Fix some $\u^0 \in \relint(U)$, $\v^0\in \relint(V)$ and $i \in \{1, \ldots, m\}$. Let $u^i_1 \leq \bar\lambda_j \leq u^i_2$ and $v^i_1 \leq \bar\mu_j \leq v^i_2$, for $j \neq i$, be real numbers. We consider the two line segments
\begin{align*}
  \textstyle\bigl\{\,\u^0 + \textstyle\sum_{j\neq i}^m \bar\lambda_j\p^j + \lambda_i\p^i \st  u^i_1 \leq
  \lambda_i \leq u^i_2 \,\bigr\} &\subseteq U_0, \\ 
  \textstyle\bigl\{\,\v^0 + \textstyle\sum_{j\neq i}^m \bar\mu_j\p^j + \mu_i\p^i \st  v^i_1 \leq
  \mu_i \leq v^i_2 \,\bigr\}&\subseteq V_0.
\end{align*}

Let $f^i\colon [u^i_1, u^i_2] \to \R$ be defined by $f^i(\lambda) = f(\u^0 + \sum_{j\neq i}^m \bar\lambda_j\p^j + \lambda \p^i)$, $g^i \colon [v^i_1, v^i_2] \to \R$ be defined by $g^i(\lambda) = g(\v^0 + \sum_{j\neq i}^m \bar\mu_j\p^j + \lambda \p^i)$ and $h^i \colon [u^i_1 + v^i_1, u^i_2 + v^i_2] \to \R$ be defined by $h^i(\lambda) = h(\u^0 + \v^0 + \sum_{j\neq i}^m (\bar\lambda_j+ \bar\mu_j)\p^j + \lambda \p^i)$. Applying Lemma~\ref{one-dim-interval_lemma}, there exists a constant $\hat c_i \in \R$ such that 
\begin{equation}\label{eq:affine-prop}
\begin{aligned}
f\bigl(\u^0 + \textstyle\sum_{j\neq i}^m \bar\lambda_j\p^j + \lambda \p^i\bigr) &= f\bigl(\u^0 +
\textstyle\sum_{j\neq i}^m \bar\lambda_j\p^j\bigr) + \hat c_i\cdot \lambda &&\textrm{ for all }
\lambda \in [u^i_1, u^i_2],\\
g\bigl(\v^0 + \textstyle\sum_{j\neq i}^m \bar\mu_j\p^j + \lambda \p^i\bigr) &= g\bigl(\v^0 + \textstyle\sum_{j\neq i}^m \bar\mu_j\p^j\bigr) + \hat c_i\cdot \lambda  &&\textrm{ for all } \lambda \in [v^i_1, v^i_2]. 
\end{aligned}
\end{equation}

Notice that this argument could be made with any other values of $\bar\lambda_j$, $j\neq i$ while using the same $\bar\mu_j$, $j\neq i$. Thus, $\hat c_i$ is independent of the values of $\bar\lambda_j$, $j\neq i$. Thus, we have $m$ real numbers $\hat c_i$, $i=1, \ldots, m$, that only depend on $f,g,h$, $L$ and the two points $\u^0\in\relint(U)$ and $\v^0\in\relint(V)$, and~\eqref{eq:affine-prop} holds for any values of $u^j_1 \leq \bar\lambda_j \leq u^j_2$, $j\neq i$.

We choose $\cve \in \R^k$  satisfying $\cve \cdot \p^i  = \hat c_i$ for all $i=1, \dots, m$ (this can be done since $\p^1, \ldots, \p^m$ are linearly independent).
 Now for any $\p \in L$ such that $\u^0 + \p \in U_0$, we can represent $\p = \sum_{i=1}^m \lambda_i \p^i$ for some $u^i_1 \leq \lambda_i \leq u^i_2$, $i=1, \dots, m$. Thus, $ f(\u^0 + \p)  =  f(\u^0 + \sum_{i=1}^m \lambda_i \p^i)$.

%
%

Now using~\eqref{eq:affine-prop} with $i = m$ we have 
\begin{align*}
  f(\u^0 + \textstyle\sum_{i=1}^m \lambda_i \p^i) 
  & = f(\u^0 + \textstyle\sum_{i=1}^{m-1}\lambda_i \p^i + \lambda_m \p^m) \\
  & = f(\u^0 + \textstyle\sum_{i=1}^{m-1}\lambda_i \p^i) + \hat c_m\cdot\lambda_m,
\end{align*}
which follows because the $\hat c_i$'s do not depend on the particular values $\lambda_i$, $i\neq m$. By applying this argument iteratively, we find that
\begin{align*}
f(\u^0 + \p) 
& = f(\u^0 + \textstyle\sum_{i=1}^m \lambda_i \p^i) \\ 
& = f(\u^0) + \textstyle\sum_{i=1}^{m}\hat c_i\cdot \lambda_i \\
& = f(\u^0) + \textstyle\sum_{i=1}^{m} \lambda_i \cve\cdot \p^i \\
& = f(\u^0) + \cve \cdot \textstyle\sum_{i=1}^{m} \lambda_i \p^i  \\
& = f(\u^0) + \cve\cdot  \p.
\end{align*}

Thus, $f(\u^0 + \p) =  f(\u^0) + \cve\cdot  \p$ for all $\p$ such that $\u^0 + \p \in U_0$, i.e., $f$ is affine with respect to $L$ over $U_0$ with gradient $\cve$. This argument can also be used to show that $g$ is affine with respect to $L$ over $V_0$ with the same gradient $\cve$ (the relations in~\eqref{eq:affine-prop} will now be used on $g$, keeping $\bar \lambda_j, j\neq i$ fixed and allowing $\bar\mu_j, j\neq i$ to vary).
\medskip

Finally, we do one more step to show that $f$ is affine with respect to $L$ over all of $U$ with gradient $\cve$. Let $\u^1, \u^2 \in U $ such that $\u^2 - \u^1 = \p' \in L$. Let $v^0_1< v^0_2 \in \R$, $i = 1, \ldots, m$ be such that $\{\, \v^0 + \lambda \p' \st v^0_1 \leq \lambda \leq v^0_2\,\} \subseteq V_0$.

Let $f^0 \colon [0, 1] \to \R$ be defined by $f^0(\lambda) = f(\u^1 + \lambda \p')$, $g^0 \colon [v^0_1, v^0_2] \to \R$ be defined by $g^0(\lambda) = g(\v^0 + \lambda \p')$ and $h^0 \colon [0 + v^0_1, 1 + v^0_2] \to \R$ be defined by $h^0(\lambda) = h(\u^1 + \v^0 + \lambda \p')$. Applying Lemma~\ref{one-dim-interval_lemma} to $f^0, g^0$ and $h^0$, there exists a constant $\hat c_0 \in \R$ such that 
\begin{subequations}
  \begin{align}
    f(\u^1 + \lambda \p') &= f(\u^1) + \hat c_0\cdot \lambda &&\textrm{for all } \lambda \in [0, 1], \label{eq:affine-prop4}\\
    g(\v^0 + \lambda \p') &= g(\v^0) + \hat c_0\cdot \lambda &&\textrm{for all } \lambda \in [v^0_1, v^0_2].
  \end{align}
\end{subequations}
Since $g$ is affine over $V_0$ with gradient $\cve$, $g(\v^0 + \lambda \p') = g(\v^0) + \lambda (\cve\cdot \p')$  for all $\lambda \in [v^0_1, v^0_2]$. Thus, $\hat c_0 = \cve\cdot \p'$. Using \eqref{eq:affine-prop4}, we get $f(\u^1 + \p') = f(\u^1) + \hat c_0 = f(\u^1) + \cve\cdot \p'$. Therefore, $f(\u^2) - f(\u^1) = \cve\cdot \p'$ as required.
%
%
%
%
%
%
The same argument applies for proving $g$ is affine with respect to $L$ over $V$ with gradient $\cve$. Finally, since $h(\x+\y) = f(\x) + g(\y)$ for all $\x \in U$, $\y\in V$, it follows that $h$ is affine with respect to $L$ over $U+V$ with gradient $\cve$.\end{proof}

\subsection{Pexider functional equation on convex additivity domains in~$\R^k$}

We now prove a technical  lemma which can be used to transfer affine properties using small ``patches'' within a larger domain. This will allow us to connect local applications of the Higher-dimensional Interval Lemma (Theorem~\ref{lem:generalized_interval_lemma}) within convex sets.  This lemma's arguments have been explicitly and implicitly used in the integer programming literature~\cite{3slope, bhkm, basu-hildebrand-koeppe:equivariant,infinite, infinite2, dey1, dey3, tspace,  Miller-Li-Richard2008}, as well as the functional equations literature~\cite{baker-rado-1987,kannappan09}.

\begin{lemma}[Patching lemma]
\label{lem:patching}
Let $U\subseteq \R^k$ be a convex subset. 
Let $\pi \colon U \to \R$ be any function. Suppose $r\colon U \to \R$ is a function such that for every $\u \in U$, 
\begin{itemize}
\item[(i)] $r(\u) > 0$, and
\item[(ii)] 
$\pi$ is affine on $B^\infty(\u, r(\u)) \cap U$.
\end{itemize} 
Then $\pi$ is affine on all of $U$.
\end{lemma}

\begin{proof} If $U$ is empty there is nothing to show. Fix any $\u^0 \in U$. Since $\pi$ is affine on $B^\infty(\u^0, r(\u^0))\cap U$, there exists $\cve \in \R^k$ such that $\pi(\u) - \pi(\u^0) =  \cve \cdot (\u - \u^0)$ for every $\u \in B^\infty(\u^0, r(\u^0))\cap U$. We claim that $\pi(\u) - \pi(\u^0) =  \cve\cdot ( \u - \u^0)$ for every $\u \in U$. This will establish the lemma. Indeed, consider $\u^1, \u^2 \in U$. $\pi(\u^2) - \pi(\u^1) = (\pi(\u^2) - \pi(\u^0)) + (\pi(\u^0) - \pi(\u^1)) =  \cve\cdot ( \u^2 - \u^0 ) -  \cve\cdot ( \u^1 - \u^0 ) =  \cve\cdot ( \u^2 - \u^1 )$.

Consider any arbitrary $\u \in U$ 
and the line segment $[\u, \u^0] \subseteq U$. For every $\x \in [\u, \u^0]$, consider $B^\infty(\x, r(\x))$. Since $r(\x) > 0$ for all $\x \in U$, $\bigcup_{\x \in [\u, \u^0]}B^\infty(\x, r(\x))$ is an open cover of $[\u, \u^0]$. Thus, there exists a finite subcover from this open cover. In particular, there exist points $\x^0, \x^1, \ldots, \x^n \in [\u, \u^0]$ such that the following hold:
\begin{itemize}
\item[(i)] $\u^0 \in B^\infty(\x^0, r(\x^0)) \cap U$, 
\item[(ii)] $\u \in B^\infty(\x^n, r(\x^n))\cap U$, and
\item[(iii)] $(B^\infty(\x^{i-1}, r(\x^{i-1})) \cap U) \cap (B^\infty(\x^i, r(\x^i))\cap U) \neq \emptyset$ for every $i = 1, \ldots, n$.
\end{itemize}

First, because of (i) and the facts that $\pi$ is affine on $B^\infty(\x^0, r(\x^0))\cap U$ and $\pi$ is affine on $B^\infty(\u^0, r(\u^0))\cap U$ with gradient $\cve$, we conclude that $\pi$ is affine with gradient $\cve$ on $B^\infty(\x^0, r(\x^0))\cap U$. From (iii), we know that $(B^\infty(\x^{i-1}, r(\x^{i-1})) \cap U) \cap (B^\infty(\x^i, r(\x^i))\cap U) \neq \emptyset$. Since $\pi$ is affine on $B^\infty(\x^0, r(\x^0))\cap U$ with gradient $\cve$ and $\pi$ is affine over $B^\infty(\x^1, r(\x^1))\cap U$, we conclude $\pi$ is affine over $B^\infty(\x^1, r(\x^1))\cap U$ with gradient $\cve$. Applying this argument repeatedly, we have that $\pi$ is affine on each $B^\infty(\x^i, r(\x^i)) \cap U$ with the same gradient $\cve$. Choose $\y^i$, $i = 1, \ldots, n$ as points in $(B^\infty(\x^{i-1}, r(\x^{i-1})) \cap U) \cap (B^\infty(\x^i, r(\x^i))\cap U)$. Therefore, since $\y^{i+1}, \y^i \in B^\infty(\x^i, r(\x^i))\cap U$ for every $i= 1, \ldots, n-1$, we have $$\pi(\y^{i+1}) - \pi(\y^i) =  \cve\cdot ( \y^{i+1} - \y^i).$$ Also, from (i) and (ii), we have $$\pi(\y^1) - \pi(\u^0) =  \cve \cdot ( \y^1 - \u^0), \qquad \pi(\u) - \pi(\y^n) =  \cve \cdot ( \u - \y^n).$$ Adding these equalities, together, we obtain $\pi(\u) - \pi(\u^0) =  \cve\cdot ( \u - \u^0).$\end{proof}

The Higher-dimensional Interval Lemma will be used to deduce affine
properties from more complicated convex sets.  Since we do not always have additivity on all of $U \times V$, we prove affine properties on smaller cross products  and then patch them together.

We will need the following basic lemma from convex analysis. 

\begin{lemma}[Theorem 6.6 in \cite{rock}]
\label{lem:relintConvex}
Let $C$ be a convex set in $\R^n$ and let $A$ be a linear transformation from $\R^n$ to $\R^m$.   Then 
$$
 A \relint(C) = \relint(A C).
 $$ 
\end{lemma}

\begin{lemma}[Relative interior lemma]
\label{lem:rel-int-7-tuple}
Let $F \subseteq \R^k \times \R^k$ be a convex set.
For any $\x \in \relint(p_1(F))$, there exist $\y \in \relint(p_2(F))$ such that $(\x,\y) \in \relint(F)$ and $p_3(\x,\y) = \x + \y\in\relint(p_3(F))$.  Similarly, for any $\y \in \relint(p_2(F))$, there exist $\x \in \relint(p_1(F))$ such that $(\x,\y) \in \relint(F)$ and $p_3(\x,\y) = \x + \y\in\relint(p_3(F))$. 
\end{lemma}

\begin{proof}
 Since $p_i \colon \R^k \times \R^k \to \R^k$ are linear transformations for $i=1,2,3$, by \autoref{lem:relintConvex}, we have 
 $p_i (\relint(F)) = \relint(p_i (F)).$
Therefore, $p_i\colon\relint(F) \to \relint(p_i(F))$ is a well defined surjective map.  

We only prove the first claim as the second has a similar proof.  Let $\x \in \relint(p_1(F)) = p_1(\relint(F))$.  Hence, there exists a point $\y \in \R^k$ such that $(\x,\y) \in \relint(F)$.  Then, for $i=2,3$, $p_i(\x,\y) \in p_i(\relint(F)) = \relint(p_i(F))$, that is, $\y \in \relint(p_2(F))$ and $\x + \y \in \relint(p_3(F))$. 
\end{proof}



\begin{definition}
  For a linear space $L \subseteq \R^k$ and a set $U \subseteq \R^k$ such that
  for some $\u \in \R^k$ we have $\aff(U) \subseteq L + \u$, we will denote by
  $\intr_L(U)$ the interior of $U$ in the relative topology of $L + \u$.
\end{definition}
Note that $\intr_L(U)$ is well defined because either $\aff(U) = L + \u$, or
$\intr_L(U) = \emptyset$.   We now prove our most general theorem relating to equation~\eqref{eq:cauchy} on a convex domain.

\begin{theorem}[Convex additivity domain lemma]
\label{lem:projection_interval_lemma}
\citedinsurveyas{Theorem 2.11}
Let $f,g,h \colon \R^k \to \R$ be bounded functions. Let $F \subseteq \R^k \times \R^k$ be a convex set
such that $f(\u) + g(\v) = h(\u+\v)$ for all $(\u, \v) \in F$. 
Let $L$ be a linear subspace of $\R^k$ such that $L \times L + F \subseteq \aff(F)$.
Let $(\u^0, \v^0) \in \relint(F)$. 
Then there exists a vector $\cve\in \R^k$ such that $f, g$ and $h$ are affine with gradient $\cve$ over $\intr_L((\u^0 + L) \cap p_1(F))$, $\intr_L((\v^0 + L) \cap p_2(F))$ and $\intr_L((\u^0 + \v^0 + L) \cap p_3(F))$, respectively.
\end{theorem}

\begin{proof}
If $\dim(L)=0$, there is nothing to prove. So we assume $\dim(L) \geq 1$.
Let $I = p_1(F)$, $J = p_2(F)$, $K = p_3(F)$.


For $\u \in \relint(I)$, define 
$$
r(\ve u) = \sup\left\{\frac r2 \in \R \, \bigg|  \, 
\begin{array}{l}
\exists \ve v \in \R^k \text{ such that }
B^\infty((\ve u, \ve v), r) \cap \big((\u,\v) + L\times L\big) \subseteq F
\end{array}
\right\}.
$$

%
%

By Lemma~\ref{lem:rel-int-7-tuple}, for any $\u \in \relint(I)$, 
there exists $\v \in \relint(J)$ such that $(\u,\v) \in\relint(F)$.  
Since $\dim(L) \geq 1$, Lemma~\ref{lem:relint-linspace} implies that $r(\u) > 0$ for every $\u \in \relint(I)$.  Let $\v \in F$ such that $B^\infty((\ve u, \ve v), r(\u)) \cap \big((\u,\v) + L\times L\big) \subseteq F$ 
and let  
$$
\begin{array}{l}
U 
= p_1\left(B^\infty((\ve u, \ve v), r(\u)) \cap \big((\u,\v) + L\times L\big)\right) 
= B^\infty(\ve u, r(\ve u))\cap (\u + L) \text{ and }\\
V 
= p_2\left(B^\infty((\ve u, \ve v), r(\u)) \cap \big((\u,\v) + L\times L\big)\right) 
= B^\infty(\ve v, r(\u))\cap (\v+ L).
\end{array}
$$  
Notice that 
$$U\times V = B^\infty((\ve u, \ve v), r(\u)) \cap \big( (\u,\v) + L\times L\big)\subseteq F.$$     Hence, applying Theorem~\ref{lem:generalized_interval_lemma} with $U$ and $V$, we obtain that $f$ is affine over $U$. 
Thus, we satisfy the hypotheses of Lemma~\ref{lem:patching} and $f$ is affine over $\intr_L((\u + L) \cap I)$ for every $\u \in \relint(I)$.
This argument can be repeated to show that $g$ is affine over $\intr_L((\v + L) \cap J)$ for every $\v \in \relint(J)$.  

For the pair $(\u^0,\v^0) \in \relint(F)$, by Lemma~\ref{lem:relint-linspace}, there exists  $r >0$ such that $B^\infty((\ve u^0, \ve v^0), r) \cap \big( (\u^0,\v^0) + L\times L\big)\subseteq F$.   Then for $U_0=B^\infty(\ve u^0, r(\u)) \cap  (\u^0 + L)$ and $V_0=B^\infty(\ve v^0, r(\u^0))\cap (\v^0 + L)$, we have $U_0 \times V_0 \subseteq F$ and   Theorem~\ref{lem:generalized_interval_lemma} also tells us that $f$ and $g$ have the same gradient $\cve$ in $U_0$ and $V_0$, respectively. Since $f$ and $g$ are affine in $\intr_L((\u^0 + L) \cap I)$ and $\intr_L((\v^0 + L) \cap J)$, respectively, we have that $f$ and $g$ are affine with the same gradient $\cve$ over all $\intr_L((\u^0 + L) \cap I)$ and $\intr_L((\v^0 + L) \cap J)$, respectively. Finally, since $f(\u) + g(\v) = h(\u+\v)$ for all $(\u, \v) \in F$, it follows that $h$ is affine over $\intr_L((\u^0+\v^0 + L) \cap K)$. This finishes the proof.\end{proof}

\begin{remark}[Comparing Theorem~\ref{lem:generalized_interval_lemma} and
  Theorem~\ref{lem:projection_interval_lemma}] 
  \label{remark:projection_interval_lemma_interiors_only}
  \citedinsurveyas{Remark 2.12}
The reader might think that the Higher-dimensional Interval Lemma
(Theorem~\ref{lem:generalized_interval_lemma}) could be obtained as a
corollary of Convex Additivity Domain Lemma 
(Theorem~\ref{lem:projection_interval_lemma}), by setting $F = U\times V$
. 
However, the Higher-dimensional Interval Lemma shows that under the appropriate additivity conditions over $U$ and $V$, we can obtain affine properties over {\em all} of $U$ and $V$ (with respect to $L$); whereas, the Convex Additivity Domain Lemma derives affine properties only over the interiors with respect to~$L$. 
This, however, cannot be avoided.  In particular, there are examples satisfying the hypotheses of Convex Additivity Domain Lemma where the functions are affine over the interiors, but not on the boundaries; see  \cite{baker-rado-1987} for such an example of a $F \subseteq \R \times \R$ and bounded functions $f,g,h$ that satisfy~\eqref{eq:pexider}, but are not affine.  

\end{remark}

\begin{remark}[Extension not valid even with all additive relations]
  \label{remark:projection_interval_lemma_interiors_only2}
The example in~\cite{baker-rado-1987} mentioned above is obtained by choosing
a subset $F\subseteq\R\times\R$ such that $F\subsetneq \{(x,y) \st x \in p_1(F), y
\in p_2(F), x +y \in p_3(F)\}$. The strict containment means that additivity
does not hold for all possible pairs $(x,y) \in p_1(F)\times p_2(F)$ such that $x + y\in p_3(F)$. We now give a similar example where the set containment is not strict, meaning that all possible additive relations from the projections are allowed. In particular, we construct an $F\subseteq\R^2\times\R^2$ such that $F  = \{\,(\x,\y) \st \x \in p_1(F),\, \y \in p_2(F),\, \x +\y \in p_3(F)\,\}$.
Let
\begin{align*}
F = \conv\Bigg(  &\begin{pmatrix} \u \\ \v \end{pmatrix} = 
  \begin{pmatrix}   0       \\      4      \\        5     \\         1       \end{pmatrix},
   \begin{pmatrix}    0     \\         3     \\         5       \\       2       \end{pmatrix},
   \begin{pmatrix}    2      \\       11/3     \\       4       \\       1       \end{pmatrix},
    \begin{pmatrix}   1           \\   4      \\        4    \\          1       \end{pmatrix},
    \begin{pmatrix}   8/3     \\      35/9     \\       4       \\       1       \end{pmatrix},
    \begin{pmatrix}   5/2      \\      4      \\        4    \\          1       \end{pmatrix},
    \begin{pmatrix}   1      \\       10/3    \\        5     \\         2       \end{pmatrix},
    \begin{pmatrix}   0       \\       4  \\            5      \\        2       \end{pmatrix},
   \begin{pmatrix}    3      \\        4      \\        5     \\         0    \end{pmatrix}
\Bigg).
\end{align*}


\begin{figure}[t!]  
\begin{center}
\ifpdf
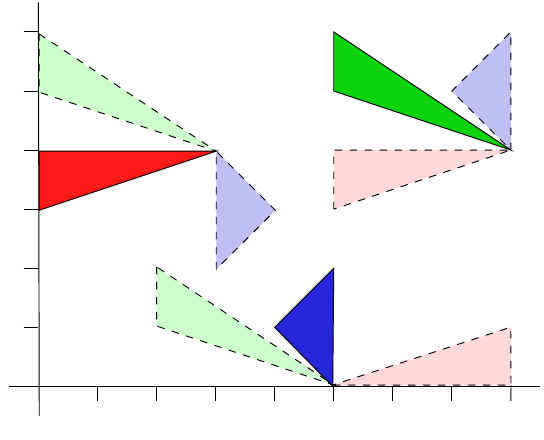
\else
\input{figureValidTripleCounterExample.pstex_t}
\fi
\end{center}
\caption{An illustration of the counterexample of Remark~\ref{remark:projection_interval_lemma_interiors_only2}.  
  The $4$-dimensional simplex~$F$ projects to the three closed triangles $U =
  p_1(F)$, $V=p_2(F)$, $W=p_3(F)$.  The points $\u,\v,\w$ are additive, i.e.,
  $\u + \v = \w$, but none of them is additive with any other points.  To see
  this, we plot the sums $U + \v$, $V + \u$, $\w + (-U)$, $W - \u$,  $\w + (-V)$, and
  $W - \v$ and show that these sets intersect $U$, $V$, and $W$ only at the
  points $\u$, $\v$, and $\w$. } 
\label{figureValidTripleCounterExample}
\end{figure}
This is a full-dimensional set of~$\R^2 \times \R^2$, which has the projections
\begin{align*}
  U = p_1(F) &= \conv\bigl(\u = \ColVec{3}{4}, \ColVec{0}{4}, \ColVec{0}{3}\bigr),\\
  V = p_2(F) &= \conv\bigl(\v=\ColVec{5}{0},\ColVec{5}{2},\ColVec{4}{1}\bigr),\\
  W = p_3(F) &= \conv\bigl(\w = \ColVec{8}{4},\ColVec{5}{6},\ColVec{5}{5}\bigr).
\end{align*}
We refer to \autoref{figureValidTripleCounterExample} for an illustration.  Furthermore, it can be shown that $F = \{\,(\x,\y) \st \x \in U,  \;\y \in V,\;  \x + \y \in W\,\}$.
Now define $f,g,h \colon \R^2 \to \R$ in the following way:
$$
f(\x) = \begin{cases}
1 &\text{ if } \x = \u,\\
0 &\text{ otherwise, } 
\end{cases}\qquad
g(\x) = \begin{cases}
2 &\text{ if } \x = \v,\\
0 &\text{ otherwise, }
\end{cases}\qquad
h(\x) = \begin{cases}
3 &\text{ if } \x = \w,\\
0 &\text{ otherwise. } 
\end{cases}
$$
\bigskip

{\em Claim 1. 
$ f(\x) + g(\y) = h(\x + \y)$ for all $(\x,\y) \in F$.}\medskip

Clearly this equation holds whenever $\x \neq \u$, $\y \neq \v$, $\x + \y \neq
\w$.   So suppose $\x = \u$.  Since $(\u + V) \cap W = \{\w\}$ and $(W - \u)\cap V = \{\v\}$, the only choice for $\y$ is $\v$.  Similarly, if we choose
$\y = \v$, the only choice for $\x$ is $\u$ or if we choose $\x + \y = \w$,
the only choices for $\x$ and $\y$ are $\u$ and $\v$.   We refer the reader to
\autoref{figureValidTripleCounterExample} to see these arguments illustrated.   
Therefore, the claim holds if and only if
$$
f(\x) + g(\y) = h(\x + \y) \text{ for all } \x \in U\setminus \{\u\},\ \y \in V\setminus \{\v\},\ \x + \y \in W\setminus \{\w\},
$$
and
$$
f(\u) + g(\v) = h(\w).
$$
Since all these equations hold, the claim is proved.\bigskip

Observe that, since $F$ is full-dimensional,  Theorem~\ref{lem:projection_interval_lemma} applies with 
$L=\R^2$. 
We deduce affine properties over the interiors of $p_1(F) = U$, $p_2(F) = V$ and $p_3(F) = W$
. This shows that Theorem~\ref{lem:projection_interval_lemma} cannot be extended to deduce
affine properties on all of $U, V, W$, unless we require further restrictions
on the types of convex sets~$F$ that we consider.  
\end{remark}


Of course, if we use the stronger regularity assumption that $f$, $g$, and $h$
are continuous functions (rather than merely bounded functions), then the
affine properties extend to the boundary as well.

\begin{corollary}[Convex additivity domain lemma for continuous functions]
\label{lem:projection_interval_lemma-corollary}
\citedinsurveyas{Corollary 2.14}
Let $f,g,h \colon \R^k \to \R$ be continuous functions. Let $F \subseteq \R^k \times \R^k$ be a convex set
such that $f(\u) + g(\v) = h(\u+\v)$ for all $(\u, \v) \in F$. 
Let $L$ be a linear subspace of $\R^k$ such that $L \times L + F \subseteq \aff(F)$.
Let $(\u^0, \v^0) \in \relint(F)$. 
Then there exists a vector $\cve\in \R^k$ such that $f, g$ and $h$ are affine with gradient $\cve$ over $(\u^0 + L) \cap p_1(F)$, $(\v^0 + L) \cap p_2(F)$ and $(\u^0 + \v^0 + L) \cap p_3(F)$, respectively.
\end{corollary}

\section[Discrete geometry of piecewise linear minimal valid functions and
  their additivity domains]{Discrete geometry of piecewise linear minimal valid functions\\ and
  their additivity domains}\label{s:prelim} 


\subsection{Polyhedral complexes and piecewise linear functions}

We introduce the notion of polyhedral complexes, which serves two purposes in
our paper.  First, it provides a framework to define piecewise linear functions.
Second, it is a tool for studying subadditivity and additivity relations of these functions.

\begin{definition}
\label{def:polyhedralComplex}
A {\em polyhedral complex} is a collection $\P$ of polyhedra in $\R^k$ such that:
\begin{enumerate}[\rm(i)]
\item $\emptyset \in \P$,
\item if $I \in \P$, then all faces of the polyhedron $I$ are in $\P$,
\item the intersection $I \cap J$ of two polyhedra $I,J \in \P$ is a face of both $I$ and $J$,
\item $\P$ is \emph{locally finite}, i.e., any compact subset of $\R^k$ intersects only finitely many faces in $\P$.
\end{enumerate}
\end{definition}
A polyhedron $I$ from $\P$ is called a {\em face} of the complex.
A polyhedral complex~$\P$ is called {\em pure} if all its maximal faces (with
respect to set inclusion) have the same
dimension. In this case, we call the maximal faces of $\P$ the {\em cells}
of~$\P$. 
A polyhedral complex~$\P$ is {\em complete} if the union of all faces of
the complex is~$\R^k$. 
The reader can find examples illustrating this and the following definitions in
\autoref{sec:triangulation}.


Given a pure and complete polyhedral complex $\P$, we call a function $\pi\colon \R^k\to \R$ \emph{continuous piecewise linear over
  $\P$} if it is affine over each of the cells of~$\P$. We introduce the following notation for a continuous piecewise linear function $\pi$ over $\P$.  
%
%

Motivated by Gomory--Johnson's characterization of minimal valid functions
(\autoref{thm:minimal}), we are interested in functions~$\pi\colon \R^k\to \R$
that are periodic modulo~$\Z^k$, i.e., for all $\ve x\in \R^k$ and all vectors 
$\ve t\in\Z^k$, we have $\pi(\ve x+\ve t) = \pi(\ve x)$.  If $\pi$ is periodic modulo~$\Z^k$
and continuous piecewise linear over a pure and complete complex~$\P$, then we
will usually assume that $\P$ is
also \emph{periodic modulo~$\Z^k$}, i.e., for all $I\in \P$ and all vectors
$\ve t\in\Z^k$, the translated polyhedron~$I + \ve t$ also is a face of~$\P$.

\begin{remark}
  Under these assumptions it is clear that there are various ways to make the
  description of~$\pi$ finite.  For example, $\tilde D := [0,1)^k$ is a
  fundamental domain (system of unique representatives) of~$\R^k$ with respect
  to the natural action of~$\Z^k$, and so it suffices to know the values of
  $\pi$ on~$\tilde D$.  However, it is inconvenient that $\tilde D$ is not closed. On
  the other hand, if we use instead its closure, $D := [0,1]^k$, we lose
  uniqueness since not every point $\ve x \in \R^k$ would have a unique decomposition as $\ve x = \ve d + \ve z$ for some $\ve d \in D$ and $\ve z \in \Z^k$.  Another viewpoint, considering polyhedral complexes of the torus  $\R^k/\Z^k$, would require more complicated definitions.  Thus, in most of this
  paper, we find it most convenient and natural to work with periodic
  functions and infinite periodic complexes.
\end{remark}

\subsection{The extended complex $\Delta \P$}
\label{section:delta-p-definition}

Let $\P$ be a pure, complete polyhedral complex of $\R^k$.
 For any $I,J,K\subseteq \R^k$, we define the set
$$F(I,J,K) = \{\,(\x,\y) \in \R^k \times \R^k \st \x \in I,\, \y \in J,\, \x + \y \in K\,\}.$$   In the specific case where $I,J,K$ are polyhedra, $F(I,J,K)$ is also a polyhedron. 
In order to study the additivity domain of a piecewise linear function over $\P$, 
we define the following family of polyhedra in $\R^k \times \R^k$,
$$\Delta\P =
\{\, F(I,J,K) \st I, J, K \in \P\, \}.$$

First, we present formulas for the projections $p_1,p_2,p_3$ of $F(I,J,K)$, as
defined in \eqref{eq:projections}, in
terms of $I,J$ and $K$.  The proofs of the simpler results of this section can be found in Appendix~\ref{s:additional-proofs:polyhedral-complexes}.
\begin{prop}\label{prop:projection}
\citedinsurveyas{Proposition 3.3}
Let $I,J,K \subseteq \R^k$.  Then
\begin{align*}
  p_1(F(I,J,K)) &= (K + (-J)) \cap I, \\
  p_2(F(I,J,K)) &= (K+ (-I)) \cap J,\\
  p_3(F(I,J,K)) & = (I+J) \cap K.
\end{align*}
\end{prop}
\begin{remark}
  Note that in general, $p_1(F(I,J,K)) \subsetneq I$, $p_2(F(I,J,K))
  \subsetneq J$, and $p_3(F(I,J,K)) \subsetneq K$. Consider $I = [0,1], J = [0,1],
  K = [1.5,2.5]$.  Then $F(I,J,K)$ is the triangle $\conv\{(1,0.5), (1,1),
  (0.5,1)\}$, so $p_1(F(I,J,K)) = [0.5, 1]$, $p_2(F(I,J,K)) = [0.5, 1]$ and
  $p_3(F(I,J,K))=[1.5,2]$.
\end{remark}
The next lemma explains the tight relation between $F$ and its projections $p_1(F), p_2(F)$ and $p_3(F)$.

\begin{lemma}
\label{lemma:setsIJKF}
\citedinsurveyas{Lemma 3.5}
Let $I,J,K \subseteq \R^k$ and let $F = F(I,J,K)$.  
Let $I' = p_1(F)$, $J' = p_2(F)$, and $K' = p_3(F)$.  Then $F =  F(I', J', K')$.
\end{lemma}
\begin{proof}
By definition of $I', J', K'$ it follows that $I' \subseteq I$, $J' \subseteq J$, $K' \subseteq K$.  
 Therefore $F(I',J',K') \subseteq F(I,J,K)$. 
 
Observe that for any $\bar F \subseteq \R^k \times\R^k$ and $(\bar{\ve x}, \bar{\ve y}) \in \bar F$, by definition we have $\bar{\ve x} \in p_1(\bar F)$, $\bar{\ve y} \in p_2(\bar F)$, and $\bar{\ve x} + \bar{\ve y} \in p_3(\bar F)$.  
Therefore $(\bar{\ve x},\bar{\ve y}) \in \{ (\ve x, \ve y) \st \ve x \in p_1(\bar F), \ve y \in p_2(\bar F), \ve x + \ve y \in p_3(\bar F)\} = F(p_1(\bar F), p_2(\bar F), p_3(\bar F))$.  Hence, $\bar F \subseteq F(p_1(\bar F), p_2(\bar F), p_3(\bar F))$.
Thus, 
\begin{displaymath}
  F(I,J,K) \subseteq F\bigl(p_1\bigl(F(I,J,K)\bigr), p_2\bigl(F(I,J,K)\bigr), p_3\bigl(F(I,J,K)\bigr)\bigr) = F(I',J',K').
\end{displaymath}
Therefore, $F(I,J,K) = F(I',J', K')$.
\end{proof}

The next lemma shows that $\Delta\P$ is a polyhedral complex, which
follows from the fact that $\P$ is a polyhedral complex. 
\begin{lemma}
  \label{lemma:delta-p-is-complex}
  \citedinsurveyas{Lemma 3.6}
  If $\P$ is a pure, complete polyhedral complex in $\R^k$, then $\Delta\P$ is a pure, complete polyhedral complex in $\R^k \times \R^k$.
\end{lemma}

Let $\pi$ be a continuous piecewise linear function over $\P$.  
We will study the function $\Delta\pi\colon \R^k \times \R^k \to \R$, as defined in Lemma~\ref{lemma:tight-implies-tight}, which measures the slack in the subadditivity constraints.

\begin{lemma}
\label{lem:delta-pi-cts}
\citedinsurveyas{Lemma 3.7}
$\Delta\pi$ is continuous piecewise linear over $\Delta\P$.
\end{lemma}

\begin{proof}
First,  $\Delta \pi$ is continuous since it  is the sum of continuous functions.

For any $F(I,J,K) \in \Delta\P$, $\Delta \pi|_{F(I,J,K)}(\x,\y) = \pi|_I(\x) + \pi|_J(\y) - \pi|_K(\x + \y)$.  Since $\pi|_I$, $\pi|_J$, $\pi|_K$ are all affine, it follows that $\Delta \pi|_{F(I,J,K)}$ is affine.  Therefore $\Delta\pi$ is affine over every face in $\Delta \P$, i.e., $\Delta \pi$ continuous piecewise linear over $\Delta \P$. 
\end{proof}

\begin{remark}\label{rem:delta-p-finite}
  If $\pi$ and $\P$ are periodic modulo~$\Z^k$, then 
  $\Delta\pi$ and $\Delta\P$ are periodic modulo~$\Z^k\times\Z^k$.  
  Indeed, let $F\in\Delta\P$, so $F = F(I,J,K)$ for some $I, J, K \in \P$. 
  Then for $(\ve u, \ve v) \in \Z^k\times\Z^k$ we have $F + (\ve u, \ve v) = 
  F(I+\ve u, J+\ve v, K+\ve u+\ve v) \in \Delta\P$. 
  In order to make the description of $\Delta\pi$ finite, 
  we can choose a fundamental domain (system of unique representatives)
  of~$\R^k\times\R^k$ with respect to the action of~$\Z^k\times\Z^k$, for
  example $\Delta \tilde D := [0,1)^k \times [0,1)^k$.
\end{remark}

\begin{remark}\label{rem:delta-p-swap}
  We remark that $\Delta\pi(\x,\y)$ is also invariant under exchanging $\x$
  and $\y$.  This can be expressed as an action of the symmetric group~$S_2$. 
  Together we obtain the action of the group~$\Z^k \wr S_2$, a wreath product, 
  and so we would be able to choose a smaller fundamental domain,
  corresponding to the action of this group.  Thus, in a practical implementation of our algorithms, this allows us to store less information when handling $\Delta \pi$ and hence improve the running time of our algorithms.\end{remark}

\subsection{Finite test for minimality of piecewise linear functions}\label{section:minimalityTest}

By Theorem~\ref{thm:minimal}, we can test whether a function is minimal by
testing subadditivity and the symmetry condition. These properties are easy to test when the function is continuous piecewise linear.  The first of such tests came from Gomory and Johnson \cite[Theorem
7]{tspace} for the case $k=1$.\footnote{Note that in \cite{tspace}, the word ``minimal'' needs to be replaced by
  ``satisfies the symmetry condition'' throughout the statement of their
  theorem and its proof.} 
Richard, Li, and Miller \cite[Theorem
22]{Richard-Li-Miller-2009:Approximate-Liftings} gave a similar
superadditivity test for discontinuous piecewise linear functions.
In~\cite{basu-hildebrand-koeppe:equivariant}, the authors gave a minimality
test for discontinuous piecewise linear functions for the $k=1$ case.  In the present paper, we give a similar test for continuous piecewise linear functions for general $k$.  As in~\cite{basu-hildebrand-koeppe:equivariant}, we do not claim novelty for these ideas.  Since our focus of this paper is classifying extreme functions and our theorems only consider minimal functions, we present these minimality tests to give a complete picture.

We assume that the function given to us is periodic and is described by a pure
and complete polyhedral complex $\P$ where every cell in $\P$ is bounded and
therefore each cell is the convex hull of its vertices.  
As we explain in section~\ref{sec:bounded-cells} of the Appendix, the assumption that every cell is bounded is not very restrictive.  
In particular, we show that every continuous minimal piecewise linear function $\pi\colon \R^k \to \R$ that satisfies a certain regularity condition called \emph{genuinely $k$-dimensional} (see Definition~\ref{def:genk}) has the property that if $\P$ is periodic modulo $\Z^k$, then every cell of $\P$ is bounded (Lemma~\ref{lemma:genKvertices}). 
Furthermore, if the function $\pi$ is not genuinely $k$-dimensional, then we can project it into a lower dimension and study it there (\autoref{prop:dim-reduction}, \autoref{remark:dimension-reduction}).

We use $\verts(\cdot)$ to denote the set of vertices of a polyhedron or
polyhedral complex.  For a polyhedral complex $\P$ in $\R^k$ and a set $S \subseteq \R^k$, we define  $S \cap \P := \{\, S \cap F \st F \in \P\,\}$.  When $S$ is a polyhedron in $\R^k$, the collection $S \cap \P$ is again a polyhedral complex.
We write $\ve 1$ to denote the vector with all entries as one, and$\pmod{\ve 1}$ to denote componentwise equivalence modulo 1. 


\begin{theorem}[Minimality test]
\label{minimality-check}
\citedinsurveyas{Theorem 3.10}
Let $\P$ be a pure, complete, polyhedral complex  in $\R^k$ that is periodic modulo $\Z^k$ and every cell of $\P$ is bounded.
Let $\pi\colon \R^k
\to\R$ be a continuous function that is periodic modulo $\Z^k$ and that is
piecewise linear function over $\P$.  Let $\Delta\tilde D = [0,1)^k \times [0,1)^k$ or another fundamental
domain as described in Remarks~\ref{rem:delta-p-finite} and \ref{rem:delta-p-swap}.
Then $\pi$ is minimal for $R_\f(\R^k, \Z^k)$  if and only if the following conditions hold:
\begin{enumerate}
\item $\pi(\0) = 0$,
\item Subadditivity test: $\Delta\pi(\u,\v) \geq 0$ for all $(\u,\v) \in \Delta\tilde D \cap \verts(\Delta \P)$.
\item Symmetry test: $\pi(\f) = 1$ and 
  $$\Delta\pi(\u,\v) = 0 \quad\text{for all}\quad
  (\u,\v)\in \Delta\tilde D \cap \verts(\Delta \P \cap\{\,(\u,\v)\st \u + \v
  \equiv \f \textstyle\pmod{\ve1}\,\}).$$\label{symmetry-test}
\end{enumerate}
\end{theorem}

\begin{proof}
  We use the characterization of minimal functions given by
  Theorem~\ref{thm:minimal}. 
  Clearly these conditions are necessary.  We will show that they are sufficient.

Since every cell of $\P$ is bounded, the cells of $\Delta \P$ are also bounded.  
By Lemma~\ref{lem:delta-pi-cts}, $\Delta \pi$ is continuous piecewise linear over $\Delta \P$. Therefore $\Delta \pi$ is completely determined by the values on $\verts(\Delta \P)$. 

Let $(\x, \y) \in \R^k \times \R^k$.  
For subadditivity, we need to show that $\Delta\pi(\x,\y) \geq 0$.  Let $F \in \Delta \P$ be such that $(\x,\y) \in F$.  Consider any vertex $(\u, \v) \in \verts(F)$.  Since $\Delta\tilde D$ is a fundamental domain for $\Z^k \times \Z^k$, and $\Delta \P$ is periodic modulo $\Z^k \times \Z^k$, there exists a point $(\w,\z) \in \Z^k \times \Z^k$ such that $(\u + \w,\v + \z)  \in \Delta\tilde D \cap \verts(\Delta \P)$.  Since $\Delta \pi$ is periodic modulo $\Z^k \times \Z^k$ and is nonnegative on $(\u + \w,\v + \z) $, we have that $\Delta \pi$ is also nonnegative on $(\u,\v)$.  Therefore $\Delta \pi$ is nonnegative on all of $\verts(F)$, and since $\Delta \pi|_F$ is affine, by convexity it follows that $\Delta \pi(\x,\y) \geq 0$. Therefore $\pi$ is subadditive.  

 Similarly, to show symmetry, we need to show that $\Delta\pi(\x,\y) = 0$ for all $\x,\y \in \R^k$ such that $\x + \y \equiv \f\pmod{\ve1}$. Observe that $\Delta\P\cap \{\,(\u,\v)\st \u + \v \equiv \f\pmod{\ve1}\,\}$ is a polyhedral complex.  Let $(\x,\y) \in \R^k$ such that $\x + \y \equiv \f\pmod{\ve1}$.  By letting $F \in \Delta\P\cap \{\,(\u,\v)\st \u + \v \equiv \f\pmod{\ve1}\,\}$ such that $(\x,\y) \in F$, the same argument as above shows that $\Delta\pi = 0$ for all vertices of $F$, and by convexity, $\Delta \pi|_F = 0$.  Therefore $\Delta \pi(\x,\y) = 0$ and we conclude that $\pi$ is symmetric.
%

Finally, we show that $\pi$ is nonnegative.  First, since $\pi$ is continuous on the compact set $[0,1]^k$, and is periodic, $\pi$ is bounded. Suppose for the sake of contradiction that $\pi(\x) < 0$ for some $\x \neq \0$.  Since $\pi$ is subadditive, $\pi(n \x) \leq n \pi(\x)$.  But since $n \pi(\x) \to -\infty$ as $n\to \infty$, this shows that $\pi$ is unbounded, which is a contradiction.
\end{proof}

\begin{remark}[Symmetry test simplified]\label{rem:f-break-point}
\citedinsurveyas{Remark 3.11}
Suppose $\P$ is pure, complete polyhedral complex that is periodic modulo $\Z^k$ and contains $\{\f\} \in \P$.        
   Then $\Delta\tilde D \cap\verts(\Delta\P \cap \{\,(\u,\v)\st \u + \v \equiv
\f\pmod{\ve 1}\,\}) \subseteq \Delta\tilde D \cap \verts(\Delta\P)$.  In particular,  the symmetry test (\ref{symmetry-test} in Theorem~\ref{minimality-check}) then reduces to checking on vertices $(\u, \v)\in \Delta\tilde D$ of $\Delta \P$ such that $\u + \v \equiv \f \pmod{\ve1}$.

To see this, consider any face $F(I,J,K) \in \Delta \P$, where $I,J,K \in \P$, and any $\z \in \Z^k$.  Then  $F(I,J,K) \cap \{(\x,\y) \in \R^k \times \R^k \st \x + \y = \f+ \z\} = F(I,J, \{\f + \z\}\cap K)$.  Since $\P$ is periodic modulo $\Z^k$ and $\{\f\}\in \P$, we have $\{\f + \z\} \in \P$.  Since $\P$ is a polyhedral complex, $\{\f + \z\} \cap K \in \P$.  Therefore, $F(I,J, \{\f + \z\}\cap K) \in \Delta \P$.  Therefore, $\verts(\Delta\P \cap \{\,(\u,\v)\st \u + \v \equiv
\f\pmod{\ve 1}\,\}) \subseteq  \verts(\Delta\P)$.  Intersecting both sides with $\Delta \tilde D$ maintains the containment relationship.

\end{remark}

\subsection{Combinatorializing the additivity domain}
\label{section:additivity-discretized}

Let $\pi \colon \R^k \to \R$ be a continuous piecewise linear function over a  pure, complete polyhedral complex $\P$.  
Recall the definition of  the \emph{additivity domain} of~$\pi$,
\begin{displaymath}
  E(\pi) = \bigl\{\,(\x, \y) \bigst \Delta \pi(\x, \y) = 0\,\bigr\}.
\end{displaymath}
We now give a combinatorial representation of this set
using the faces of~$\P$; this extends a technique
in~\cite{basu-hildebrand-koeppe:equivariant}.  
Let 
$$
E(\pi, \P) = \bigl\{\, F \in \Delta \P \bigst \Delta\pi|_F = 0\,\bigr\}.
$$
We consider $E(\pi, \P)$ to include $F=\emptyset$, on which $\Delta\pi|_F = 0$
holds trivially.  Then $E(\pi, \P)$ is another polyhedral complex, a
subcomplex of $\Delta\P$.   
As mentioned, if $\pi$ is continuous, then $\Delta \pi$ is continuous.
Under this continuity assumption, we can consider only the set of maximal faces in $E(\pi, \P)$.  We define 
$$
E_{\max{}}(\pi,\P) = \bigl\{\,F\in E(\pi, \P)\bigst F \text{ is a maximal face by set inclusion in } E(\pi, \P)\,\bigr\}. 
$$
  
\begin{lemma}\label{lemma:covered-by-maximal-valid-triples}
\citedinsurveyas{Lemma 3.12}
Suppose that $\pi$ is subadditive.  Then
\begin{equation*}
E(\pi) = \bigcup \{ F  \in E(\pi, \P)  \}
= \bigcup \{ F \in E_{\max{}}(\pi, \P)  \}.
\end{equation*}
\end{lemma}

\begin{proof}
Clearly $E(\pi) \supseteq \bigcup \{ F  \in E(\pi, \P)  \}
\supseteq \bigcup \{ F \in E_{\max{}}(\pi, \P)  \}$. 
We show the reverse inclusions.
Suppose $(\x,\y) \in E(\pi)$.  Since $\Delta \P$ is a polyhedral complex that covers all of $\R^k\times \R^k$, there exists a face $F \in \Delta \P$ such that $(\x,\y) \in \relint(F)$.  Note that if $(\x,\y) \in \verts(\Delta \P)$, then $F = \{ (\x,\y)\}$ is $0$-dimensional face of $\Delta \P$.
Suppose that $(\x,\y) \not\in \verts(\Delta\P)$.  Since $\pi$ is subadditive, $\Delta \pi \geq 0$.  
Further, since $\Delta \pi$ is affine in $F$, $(\x,\y) \in \relint(F)$, and $\Delta\pi(\x,\y) = 0$,  we have that $\Delta\pi|_{F} = 0$.  Therefore, $F \in E(\pi,\P)$ and $(\x,\y) \in F$ is contained in the first right hand side.   Clearly, if $F$ is not maximal in $E(\pi, \P)$, then it is contained in a maximal face $F' \in E_{\max{}}(\pi,\P)$, and hence the reverse inclusions also hold.
\end{proof}
This combinatorial
representation can then be made finite by choosing representatives under the action
of~$\Z^k\times\Z^k$, which leaves $E(\pi)$ and thus $E(\pi, \P)$ and
$E_{\max}(\pi, \P)$ invariant, as in \autoref{rem:delta-p-finite}.

\subsection{Non-extremality via perturbation functions}
We now give a method of showing $\pi$ is not extreme when we are given a
certain piecewise linear perturbation function $\bar \pi$. 

\begin{theorem}[Perturbation]
\label{corPerturb}
\citedinsurveyas{Theorem 3.13}
Let $\P$ be a pure, complete, polyhedral complex  in $\R^k$ that is periodic modulo $\Z^k$ and every cell of $\P$ is bounded.
Suppose $\pi$ is minimal and continuous piecewise linear over $\P$.  Suppose
$\bar \pi \not\equiv 0$ is continuous piecewise linear over $\P$, is periodic
modulo $\Z^k$ and satisfies $E(\pi) \subseteq E(\bar \pi)$ and $\bar
\pi(\f) = 0$. Then $\pi$ is not extreme. Furthermore, given $\bar\pi$, there
exists an $\epsilon > 0$ such that $\pi^1 = \pi + \epsilon \bar \pi$ and
$\pi^2 = \pi - \epsilon \bar \pi$ are distinct minimal functions that are
continuous piecewise linear over $\P$ such that $\pi = \tfrac12(\pi^1 +
\pi^2)$. 
\end{theorem}

  \begin{proof}

Let $\Delta\tilde D = [0,1)^k \times [0,1)^k$ (or any other fundamental domain
as in Remarks~\ref{rem:delta-p-finite} and \ref{rem:delta-p-swap}).
Let 
$$
\epsilon = \frac{1}{2} \frac{\min(\Delta \pi(\x,\y) \st (\x, \y) \in \Delta\tilde D \cap \verts(\Delta \P) , \Delta \pi(\x, \y) \neq 0)}{\max( |\Delta \bar \pi(\u, \v)| \st (\u, \v) \in \Delta\tilde D \cap\verts(\Delta \P), \Delta \bar \pi(\u, \v) \neq 0)}.
$$
Note that $\epsilon$ exists and $\epsilon > 0$  since $\Delta \pi$ and $\Delta \bar \pi$ are non-zero somewhere, $\Delta \pi$ is a nonnegative function because $\pi$ is minimal, and $\verts(\Delta \P) \neq \emptyset$ since $\Delta\P$ is a collection of bounded polyhedra. 

Setting $\pi^1 = \pi + \epsilon \bar \pi, \pi^2 = \pi - \epsilon \bar \pi$, we see that $\pi^1, \pi^2$ are piecewise linear and periodic modulo $\Z^k$.  We show that $\pi^1, \pi^2$ satisfy conditions (1), (2), and (3) of  Theorem~\ref{minimality-check} to show that $\pi^1, \pi^2$ are minimal functions.  Since $\epsilon >0$ and $\bpi \not\equiv 0$, $\pi^1, \pi^2$ are then distinct minimal functions that show that $\pi$ is not extreme.

We use the assumption that $E(\pi) \subseteq E(\bpi)$, which implies that $\Delta \bpi(\x,\y) = 0$ whenever $\Delta \pi(\x,\y) = 0$.  

First, $\Delta \pi(\0,\0) = \pi(\0) + \pi(\0) - \pi(\0) = \pi(\0) = 0$, therefore $0 = \Delta\bpi(\0,\0) = \bpi(\0)$.  Therefore $\pi^1(\0) = \pi^2(\0) = 0$.  Since $\bpi(\f) =0$ and $\pi(\f) = 1$, it follows that $\pi^1(\f) = \pi^2(\f) =1$.  These results along with $E(\pi) \subseteq E(\bpi)$ satisfy conditions (1) and (3) and Theorem~\ref{minimality-check}.  

Next, for any $(\x, \y) \in \Delta\tilde D \cap \verts(\Delta \P)$, from the definition of $\epsilon$ and the fact that $E(\pi) \subseteq E(\bpi)$, which implies that $\Delta \bpi(\x,\y) = 0$ whenever $\Delta \pi(\x,\y) = 0$,  we have 
$$
\Delta \pi(\x, \y) \pm \epsilon \Delta \bar \pi(\x, \y) \geq \Delta \pi(\x, \y) - \epsilon |\Delta \bar\pi(\x,\y)| \geq \tfrac{1}{2} \Delta \pi(\x, \y) \geq 0.
$$
Therefore $\pi^1, \pi^2$ satisfy also condition (2) of Theorem~\ref{minimality-check}, and we are done.
\end{proof}

\section{A class of minimal valid functions defined over $\R^2$}\label{sec:triangulation}

We now define the class of {\em diagonally constrained} functions $\pi \colon
\R^2 \to \R$. 
We first introduce a special two-dimensional polyhedral complex. The functions will be continuous piecewise linear over this complex.

\subsection{The standard triangulations $\P_q$ of~$\R^2$ and their geometry}\label{s:geometry-Pq}
\citedinsurveyas{Section 4.1}

 Let $q$ be a positive integer.  Consider the
arrangement~$\mathcal H_q$ of all hyperplanes (lines) of~$\R^2$ of the form
$\ColVec{0}{1}\cdot \x = b$, $\ColVec{1}{0}\cdot \x = b$, and $\ColVec{1}{1}\cdot\x  = b$,
where $b \in \tfrac{1}{q}\Z$.  The complement of the arrangement~$\mathcal
H_q$ consists of two-dimensional cells, whose closures are the triangles
$$\FundaTriangleLower = \tfrac1q \conv(\{ \ColVec{0}{0}
, \ColVec{1}{0}
, \ColVec{0}{1}
\})\qquad\text{and}\qquad \FundaTriangleUpper = \tfrac1q \conv(\{\ColVec{1}{0}
, \ColVec{0}{1}
,
\ColVec{1}{1}
\})$$ and their translates by elements of the lattice $\smash[t]{\frac1q\Z^2}$. 

We denote by $\P_q$ the collection of these triangles and the vertices and
edges that arise as intersections of the triangles, and the empty set.  Thus $\P_q$ is a locally finite
polyhedral complex that is periodic modulo $\Z^2$.  Since all
nonempty faces of~$\P_q$ are simplices, it is a triangulation of the
space~$\R^2$. 

\begin{example}\label{ex:diag-constrained-function}
  \autoref{figure:diagonallyConstrained-new-figure}, which appeared in the
  introduction,
  shows the complex~$\P_5$
  with an example of a minimal valid continuous piecewise linear function
  on~$\P_5$ with $\f = \ColVec[5]{2}{2}$ that is periodic modulo~$\Z^2$.  The function is uniquely
  determined by its values on the vertices of $\P_5$ that lie within the
  fundamental domain $\tilde D = [0,1)^2$.  Note that, due the
  periodicity of the function modulo~$\Z^2$, the values of the function on the
  left and the right edge (and likewise on the bottom and the top edge) of $D = [0,1]^2$ match. 
\end{example}

There is a partial ordering structure on the family of triangulations~$\P_q$,
whose importance to us will become clear later: For every $m>1$, the
triangulation $\P_{mq}$ is a subtriangulation (refinement) of $\P_q$, i.e.,
every face of $\P_q$ is a union of faces of~$\P_{mq}$.

Within the polyhedral complex $\P_q$, let $\Ipoint$ be the set
of 0-faces (vertices), 
$\Iedge$ be the set of 1-faces (edges), and $\Itri$ be the set of 2-faces
(triangles).  The sets of diagonal, vertical, and horizontal edges will be denoted by
$\Idiag$, $\Ivert$, and $\Ihor$, respectively. We also use abbreviations such as $\Ipointdiag =
\Ipoint\cup\Idiag$, $\Ipointdiagtri = \Ipoint\cup\Idiag\cup\Itri$, etc. 

\begin{remark}\label{rem:matrixA}
Let 
$$
A =  \begin{bmatrix}  1 & -1 &  0 & 0 &  1 & -1 \\
 0 & 0 & 1 & -1  & 1 & -1  \end{bmatrix}^T.
$$
Then for every face $I \in \P_q$, there exists a vector $\ve b \in \frac1q\Z^6$ such that $I = \{\, \x \st A \x \leq \ve b\,\}.$ Furthermore, for every vector $\ve b \in \frac1q\Z^6$, the set $\{\, \x \st A \x \leq \ve b\,\}$ is a union of faces of $\P_q$ (possibly empty), since each inequality corresponds to a hyperplane in the arrangement $\mathcal{H}_q$.
\end{remark}


The matrix $A$ is totally unimodular.  Thus 
the specific choice of the triangulation $\P_q$ lends itself to strong
unimodularity properties that reveal structure in the complex.
More importantly, they allow us to develop a simple theory of extremality, 
in which all relevant properties of the function can be expressed using the
faces of the original complex~$\P_q$.

The following lemma can be shown by enumerating cases and using simple 2-dimensional geometry.  We give an alternate proof that utilizes the total unimodularity of $A$ and avoids case analysis.

\begin{lemma}\label{lemma:I+J}
Let $I, J \in \P_q$.  Then $-I$ and $I + J$ are unions of faces in $\P_q$.
\end{lemma}
\begin{proof}
If $I = \{\,\x \in \R^2 \st A\x \leq \b\,\}$ for some $\b \in \frac1q\Z^6$, then $-I = \{\,\x \in \R^2 \st -A\x \leq \b\,\}$. Since $-A$ has the same rows as $A$ (with a permutation), by Remark~\ref{rem:matrixA}, $-I$ is a union of faces of $\P_q$.

We now 
show that
the Minkowski sum $I+J$ is a union of faces in $\P_q$. 

Let $\ve a^i$ be the $i^\text{th}$ row vector of $A$. 
Then there exist vectors $\ve b^1, \ve b^2$ such that $I = \{\, \x \st A \x
\leq \ve b^1\,\}$, $J = \{\, \y \st A \y \leq \ve b^2\,\}$.
Moreover, due to the total unimodularity of the matrix~$A$, the right-hand
side vectors
$\ve b^1, \ve b^2$ can be chosen so that $\ve b^1, \ve b^2$ are tight, i.e., 
\begin{equation}
\max_{\x \in I} \ve a^i \cdot \x = \ve b^1_i, \quad    \max_{\y \in J}  \ve a^i \cdot \y = \ve b^2_i,
\label{eq:support-vector}
\end{equation}
and $\ve b^1, \ve b^2 \in \frac1q\Z^6$.

We claim that $I + J = \{\, \x \st A \x \leq \ve b^1 + \ve b^2 \,\}$.  Clearly
$I+J \subseteq \{\, \x \st A \x \leq \ve b^1 + \ve b^2 \,\}$.  We show the
reverse direction.    
Let $K'$ be a facet (edge) of $I+J$.  Then $K'=I'+J'$, where $I'$ is a face of
$I$ and $J'$ is a face of $J$.  Without loss of generality, assume that $I'$
is an edge; then $J'$ is either a vertex or an edge.  
By well-known properties of Minkowski sums, the normal cone of $K'$ is the
intersection of the normal cones of $I'$ in $I$ 
and $J'$ in $J$.  Thus $K'$ has the same normal direction as the facet
(edge)~$I'$.  (This argument relied on the fact that we are in dimension two.)
This proves that $I + J = \{\, \x \st A \x \leq \ve b\,\}$ for
some vector $\ve b$. 

Let $\x^*$, $\ve y^*$ be maximizers in~\eqref{eq:support-vector}.  Then $\x^* + \y^* \in I + J$, and thus
$$
\ve b^1_i  + \ve b^2_i  = \ve a^i \cdot \x^* +  \ve a^i \cdot \y^* \leq  \max_{\ve z \in I + J}  \ve a^i \cdot \ve z \leq \max_{\x \in I} \ve a^i \cdot \x  + \max_{\y \in J}  \ve a^i \cdot \y = \ve b^1_i + \ve b^2_i.
$$
Therefore, $\max_{\ve z \in I + J}  \ve a^i \cdot \ve z = \ve b^1_i + \ve
b^2_i$,  which shows that every constraint $\ve a^i \cdot \ve z \leq \ve b^1_i$ is
met at equality, and therefore $I+J = \{\, \x \st A \x \leq \ve b^1 + \ve b^2\,\}$
and we conclude that $I + J$ is a union of subsets in $\P_q$. 
\end{proof}

This result has an important consequence for the complex $\Delta \P_q$, allowing component
projections of the faces of $\Delta \P_q$ to be faces of $\P_q$.

\begin{lemma}
  \label{lemma:vertices}
  \begin{enumerate}[\rm(i)]
  \item Let $F \in \Delta \P_q$. Then the projections $p_1(F)$, $p_2(F)$, and
    $p_3(F)$ are faces in the complex~$\P_q$.
  \item In particular, let $(\x,\y)$ be a vertex of $\Delta\P_q$.  Then
    $\x,\y$ are vertices of the complex $\P_q$, i.e., $\x,\y\in
    \frac{1}{q}\Z^2$.
  \end{enumerate}
\end{lemma}


\begin{proof}
  By definition of $\Delta \P_q$, there exist $I,J,K\in\P_q$ such that $F =
  F(I,J,K)$.  Let $I' = p_1(F)$, $J' = p_2(F)$, and $K' = p_3(F)$.  
  By Proposition~\ref{prop:projection},
  \begin{align*}
    I' = p_1(F) &= (K + (-J)) \cap I, \\
    J' = p_2(F) &= (K+ (-I)) \cap J,\\
    K' = p_3(F) & = (I+J) \cap K,
  \end{align*}
  and thus, by Lemma~\ref{lemma:I+J}, $I'$, $J'$, and $K'$ are faces of~$\P_q$.
\end{proof}

\begin{theorem}[Simplified minimality test]
\label{minimality-check-2d}
\citedinsurveyas{Theorem 4.5}
Let $\pi \colon \R^2 \to \R$ be a continuous piecewise linear function  over $\P_q$ that is periodic modulo $\Z^2$.  Suppose $\f \in \verts(\P_q)$.  Then $\pi$ is minimal for $R_\f(\R^2, \Z^2)$ if
and only if the following conditions hold.
\begin{enumerate}
\item $\pi(\0) = 0$.
\item Subadditivity test: $\pi(\x) + \pi(\y) \geq \pi(\x+ \y)$ for all $\x,\y \in \frac{1}{q}\Z^2\cap [0,1)^2$.
\item  Symmetry test: $\pi(\x) + \pi(\f- \x) = 1$ for all $\x \in \frac{1}{q} \Z^2\cap [0,1)^2$.
\end{enumerate}
\end{theorem}

\begin{proof}
Since $\{\f\}\in  \P_q$, the result follows by applying Theorem~\ref{minimality-check} with $\Delta\tilde D = [0,1)^2\times [0,1)^2$ and using Remark~\ref{rem:f-break-point} and Lemma~\ref{lemma:vertices} to show that the vertices that need to be considered are vertices $(\x,\y) \in \tfrac{1}{q}\Z^2 \times \tfrac{1}{q}\Z^2 \cap \Delta \tilde D$.
\end{proof}

\begin{example}[\autoref{ex:diag-constrained-function}, continued]\label{ex:diag-constrained-function-continued}
  We now visualize the additive faces $F \in E(\pi,\P_q)$ 
  (\autoref{figure:diagonallyConstrained-new-figure}); following
  \autoref{lemma:covered-by-maximal-valid-triples}, we are particularly
  interested in the maximal additive faces $\bar F \in E_{\max}(\pi,\P_q)$. 
  Following \autoref{rem:delta-p-finite}, $E(\pi,\P_q)$ is invariant
  under the action of $\Z^k\times\Z^k$. By the construction of~$\P_q$, we can always choose a representative  
  $\bar F \in E_{\max}(\pi,\P_q)$ that is a subset of the closure $\Delta D = [0,1]^2 \times
  [0,1]^2$ of the fundamental domain.  Then all faces $F \in E(\pi,\P_q)$ with
  $F \subseteq \bar F$ also are subsets of $\Delta D$. 
  
  By \autoref{lemma:setsIJKF}, each $F \in E(\pi,\P_q)$ is determined by its
  projections $I = p_1(F)$, $J = p_2(F)$, $K = p_3(F)$ as $F=F(I,J,K)$.  Due
  to the choice of triangulation~$\P_q$, by \autoref{lemma:vertices}, $I$,
  $J$, and $K$ are faces of~$\P_q$. When $F \subseteq \Delta D = [0,1]^2\times
  [0,1]^2$, we have $I,J \subseteq D = [0,1]^2$ and $K \subseteq 2D =
  [0,2]^2$. 

  Thus we can visualize faces $F\subseteq \Delta D$ by showing three diagrams,
  corresponding to its projections $p_i(F) \in \P_q$, where $p_1(F), p_2(F)
  \subseteq D$ and $p_3(F) \subseteq 2D$ as follows.  For example, consider
  the face $\bar F$ with
  $$
  p_1(\bar F) = \myCVthree{3}{2}{3}{1}{2}{2}{CornflowerBlue}\,, \qquad 
  p_2(\bar F) = \myCVthree{4}{0}{5}{0}{4}{1}{CornflowerBlue}\,, \qquad
  p_3(\bar F) = \myCVthree[1201]{3}{2}{3}{1}{2}{2}{CornflowerBlue}.
  $$
  It is a maximal additive face.  It has of course many smaller included
  faces, for example $F$ given by
  $$
  p_1(F) = \myCVtwo{3}{1}{3}{2}\,, \qquad 
  p_2(F) = \myCVthree{4}{0}{5}{0}{4}{1}{CornflowerBlue}\,, \qquad 
  p_3(F) = \myCVthree[1201]{3}{2}{3}{1}{2}{2}{CornflowerBlue}.
  $$
  Here, $p_1(F) \in \Ivert$, $p_2(F) \in \P_{q,\triup}$, and $p_3(F) \in \P_{q,\tridown}$.  

  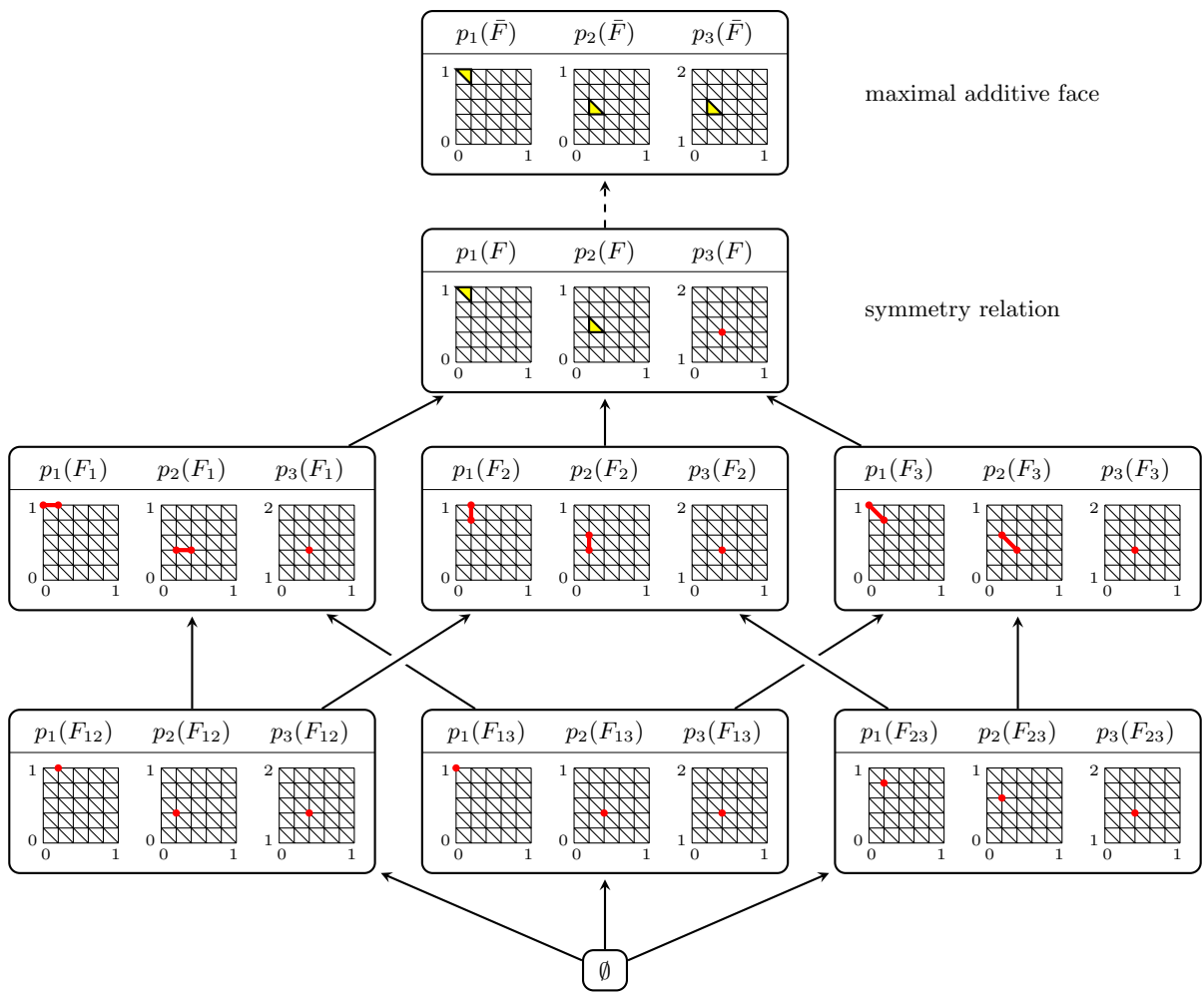
\begin{figure}[tp]
    \centering

\newcommand\FigureFace[4]
{\begin{tabular}{c@{\hspace{0.6em}}c@{\hspace{0.6em}}c}
    $p_1(#4)$ & $p_2(#4)$ & $p_3(#4)$ \\
    \cmidrule{1-3}
    #1 & #2 & #3
  \end{tabular}}
  
\begin{tikzpicture} [ font = \small, align = flush center, >=stealth, thick, node distance = 2.3cm, 
  every node/.style={draw, inner xsep=0pt, rounded corners}]
  \node [fill=white] (Fmax) 
  {\FigureFace{\myCVthree{1}{5}{1}{4}{0}{5}{Yellow}}
    {\myCVthree{1}{2}{2}{2}{1}{3}{Yellow}}
    {\myCVthree[0112]{1}{2}{2}{2}{1}{3}{Yellow}}{\bar F}};
  \node[draw=none, right = 1cm of Fmax] {maximal additive face};
  \node [fill=white, below=0.7cm of Fmax] (F) 
  {\FigureFace{\myCVthree{1}{5}{1}{4}{0}{5}{Yellow}}
    {\myCVthree{1}{2}{2}{2}{1}{3}{Yellow}}
    {\myCVone[0112]{2}{2}}{F}};
  \node[draw=none, right = 1cm of F] {symmetry relation};
  \node [fill=white, below = 0.7cm of F] (F2) 
  {\FigureFace{\myCVtwo{1}{4}{1}{5}}{\myCVtwo{1}{2}{1}{3}}{\myCVone[0112]{2}{2}}{F_2}}; 
  \node [fill=white, left = 0.6cm of F2] (F1) 
  {\FigureFace{\myCVtwo{0}{5}{1}{5}}{\myCVtwo{1}{2}{2}{2}}{\myCVone[0112]{2}{2}}{F_1}};
  \node [fill=white, right = 0.6cm of F2] (F3) 
  {\FigureFace{\myCVtwo{1}{4}{0}{5}}{\myCVtwo{2}{2}{1}{3}}{\myCVone[0112]{2}{2}}{F_3}};
  \node [fill=white, below = 1.3cm of F2] (F13) 
  {\FigureFace{\myCVone{0}{5}}{\myCVone{2}{2}}{\myCVone[0112]{2}{2}}{F_{13}}};
  \node [fill=white, below = 1.3cm of F3] (F23) 
  {\FigureFace{\myCVone{1}{4}}{\myCVone{1}{3}}{\myCVone[0112]{2}{2}}{F_{23}}};
  \node [fill=white, below = 1.3cm of F1] (F12) 
  {\FigureFace{\myCVone{1}{5}}{\myCVone{1}{2}}{\myCVone[0112]{2}{2}}{F_{12}}};
  \node [fill=white, below = 1cm of F13, inner ysep=4pt, inner xsep=6pt] (Empty) 
  {$\emptyset$};
  \draw[->,shorten >= 2pt] (Empty) -- (F13);
  \draw[->,shorten >= 2pt] (Empty) -- (F23);
  \draw[->,shorten >= 2pt] (Empty) -- (F12);
  \draw[->,shorten >= 2pt] (F13) -- (F1);
  \draw[->,shorten >= 2pt] (F13) -- (F3);
  \draw[->,shorten >= 2pt] (F12) -- (F1);
  \draw[->,shorten >= 2pt,preaction={draw=white, -,line width=6pt,shorten >=
    .6cm, shorten <= 1cm}] (F12) -- (F2);
  \draw[->,shorten >= 2pt,preaction={draw=white, -,line width=6pt,shorten >=
    .6cm, shorten <= 1cm}] (F23) -- (F2);
  \draw[->,shorten >= 2pt] (F23) -- (F3);
  \draw[->,shorten >= 2pt] (F1) -- (F);
  \draw[->,shorten >= 2pt] (F2) -- (F);
  \draw[->,shorten >= 2pt] (F3) -- (F);
  \draw[->,dashed,shorten >= 2pt] (F) -- (Fmax); 
\end{tikzpicture}

    \caption{A non-maximal additive face~$F\in E(\pi, \P_q)$ corresponding to a symmetry
      relation, the poset of its faces, and a maximal additive face $\bar F \in
      E_{\max}(\pi,\P_q)$ with $F \subset \bar F$. The
      triangles in these diagrams are colored yellow, matching
      Figure~\ref{figure:diagonallyConstrained-new-figure}, while points and
      edges are colored red.   
      This figure reveals that there are many additive faces $F \in E(\pi, \P_q)$ that do not appear in $E_{\max}(\pi,\P_q)$ and hence are not recorded in \autoref{tab:example-emax}.  Furthermore, notice that face $F_1$, with its projections described above, is not a valid maximal additive face for a diagonally constrained function.  This diagram explains that $F_1$ is not maximal in $E(\pi, \P_q)$, and hence does not contradict the fact that $\pi$ is diagonally constrained.
      }  
    \label{fig:symmetry-poset-example}
  \end{figure}
  Since $\pi$ is a minimal valid function, the symmetry condition implies that for any
  face $I \in \P_q$, we have $F(I,\f - I, \{\f\}) \in E(\pi, \P_q)$; but these
  are not necessarily maximal additive faces, even when $I \in \Itri$.  
  We illustrate this in \autoref{fig:symmetry-poset-example}, which shows a
  face $F = F(I,\f - I, \{\f\})$ with $I = p_1(F) \in \Itri$ with a
  containing maximal additive face $\bar F$ and the poset of the faces of~$F$.


  \autoref{tab:example-emax} shows all maximal additive faces 
  $F \in E_{\max}(\pi, \P_q)$ after all the faces arising from the symmetry
  condition have been removed.  Following \autoref{rem:delta-p-swap},
  $F(I,J,K) \in E_{\max}(\pi, \P_q)$ if and only if $F(J,I,K) \in E_{\max}(\pi,
  \P_q)$, so we have also removed the redundancy of swapping $I$ and
  $J$ by choosing either one of the two
  representatives arbitrarily. 

\begin{table}[tp]
\caption{All maximal faces $F \in E_{\max}(\pi, \P_q)$ of the function~$\pi$
  from \autoref{ex:diag-constrained-function}, except for the faces
  corresponding to the symmetry condition
  . The triangles in these diagrams are colored to match
  Figure~\ref{figure:diagonallyConstrained-new-figure}, while points and edges
  are just colored red.  Type numbers refer to \autoref{lemma:cases}.
  Notice that none of the light green triangles, for instance, the triangle with vertices $\ColVec{1/5}{0}, \ColVec{1/5}{1/5}, \ColVec{2/5}{0}$ appear in the table below.  This is because the only additive relations these triangles satisfy are from the symmetry condition, which we do not list below.
  } 
\label{tab:example-emax}
\centering
\begin{tabular}{c@{\hspace{0.6em}}c@{\hspace{0.6em}}clc@{\hspace{0.6em}}c@{\hspace{0.6em}}clc@{\hspace{0.6em}}c@{\hspace{0.6em}}c}
  \toprule
  $p_1(F)$ & $p_2(F)$ & $p_3(F)$ &&
  $p_1(F)$ & $p_2(F)$ & $p_3(F)$ &&
  $p_1(F)$ & $p_2(F)$ & $p_3(F)$ \\
  \cmidrule{1-3}  \cmidrule{5-7} \cmidrule{9-11}   
  \multicolumn{11}{c}{Three triangles (Type 2)}\\
  \midrule
  \myCVthree{0}{0}{1}{0}{0}{1}{Sepia} & \myCVthree{0}{0}{1}{0}{0}{1}{Sepia} & \myCVthree{0}{0}{1}{0}{0}{1}{Sepia} &&%
  \myCVthree{0}{0}{1}{0}{0}{1}{Sepia} & \myCVthree{2}{2}{2}{1}{1}{2}{Sepia} & \myCVthree{2}{2}{2}{1}{1}{2}{Sepia}&&%
    \myCVthree{1}{5}{1}{4}{0}{5}{Yellow} & \myCVthree{1}{5}{1}{4}{0}{5}{Yellow} & \myCVthree[0112]{1}{5}{1}{4}{0}{5}{Yellow}\\[5.0ex]
  \myCVthree{0}{4}{1}{4}{0}{5}{RoyalBlue} & \myCVthree{0}{4}{1}{4}{0}{5}{RoyalBlue} & \myCVthree[0112]{0}{4}{1}{4}{0}{5}{RoyalBlue}&&
  \myCVthree{0}{4}{1}{4}{0}{5}{RoyalBlue} & \myCVthree{2}{3}{2}{2}{1}{3}{RoyalBlue} & \myCVthree[0112]{2}{3}{2}{2}{1}{3}{RoyalBlue}&&%
  \myCVthree{1}{5}{1}{4}{0}{5}{Yellow} & \myCVthree{1}{2}{2}{2}{1}{3}{Yellow} & \myCVthree[0112]{1}{2}{2}{2}{1}{3}{Yellow}\\[5.0ex]
  \myCVthree{4}{0}{5}{0}{4}{1}{CornflowerBlue} & \myCVthree{4}{0}{5}{0}{4}{1}{CornflowerBlue} & \myCVthree[1201]{4}{0}{5}{0}{4}{1}{CornflowerBlue}&&
  \myCVthree{4}{0}{5}{0}{4}{1}{CornflowerBlue} & \myCVthree{4}{0}{5}{0}{4}{1}{CornflowerBlue} & \myCVthree[1201]{4}{1}{4}{0}{3}{1}{CornflowerBlue}&&%
  \myCVthree{4}{0}{5}{0}{4}{1}{CornflowerBlue} & \myCVthree{4}{0}{5}{0}{4}{1}{CornflowerBlue} & \myCVthree[1201]{3}{1}{4}{1}{3}{2}{CornflowerBlue}\\[5.0ex]%
  \myCVthree{4}{0}{5}{0}{4}{1}{CornflowerBlue} &  \myCVthree{4}{1}{4}{0}{3}{1}{CornflowerBlue} &  \myCVthree[1201]{4}{1}{4}{0}{3}{1}{CornflowerBlue}&&
   \myCVthree{4}{0}{5}{0}{4}{1}{CornflowerBlue} & \myCVthree{3}{2}{3}{1}{2}{2}{CornflowerBlue} &  \myCVthree[1201]{3}{2}{3}{1}{2}{2}{CornflowerBlue}&&%
  \myCVthree{4}{0}{5}{0}{4}{1}{CornflowerBlue} &\myCVthree{4}{1}{4}{0}{3}{1}{CornflowerBlue} &  \myCVthree[1201]{3}{2}{3}{1}{2}{2}{CornflowerBlue}\\[5.0ex]%
  \myCVthree{4}{0}{5}{0}{4}{1}{CornflowerBlue} & \myCVthree{4}{1}{4}{0}{3}{1}{CornflowerBlue} &  \myCVthree[1201]{3}{1}{4}{1}{3}{2}{CornflowerBlue}&&%
  \myCVthree{4}{0}{5}{0}{4}{1}{CornflowerBlue} &  \myCVthree{3}{1}{4}{1}{3}{2}{CornflowerBlue} &  \myCVthree[1201]{3}{2}{3}{1}{2}{2}{CornflowerBlue}&&%
  \myCVthree{4}{0}{5}{0}{4}{1}{CornflowerBlue} & \myCVthree{3}{1}{4}{1}{3}{2}{CornflowerBlue} &  \myCVthree[1201]{3}{1}{4}{1}{3}{2}{CornflowerBlue}\\[5.0ex]
   \myCVthree{4}{1}{4}{0}{3}{1}{CornflowerBlue} & \myCVthree{4}{1}{4}{0}{3}{1}{CornflowerBlue} & \myCVthree[1201]{3}{2}{3}{1}{2}{2}{CornflowerBlue}&&%
  \myCVthree{4}{1}{4}{0}{3}{1}{CornflowerBlue} & \myCVthree{3}{1}{4}{1}{3}{2}{CornflowerBlue} & \myCVthree[1201]{3}{2}{3}{1}{2}{2}{CornflowerBlue}&&
    \myCVthree{5}{5}{5}{4}{4}{5}{Blue} & \myCVthree{5}{5}{5}{4}{4}{5}{Blue} & \myCVthree[1212]{5}{5}{5}{4}{4}{5}{Blue}\\[5.0ex]%
  \myCVthree{5}{1}{5}{0}{4}{1}{BrickRed} & \myCVthree{5}{1}{5}{0}{4}{1}{BrickRed} & \myCVthree[1201]{5}{1}{5}{0}{4}{1}{BrickRed}&&%
  \myCVthree{5}{1}{5}{0}{4}{1}{BrickRed} &  \myCVthree{2}{1}{3}{1}{2}{2}{BrickRed}   & \myCVthree[1201]{2}{1}{3}{1}{2}{2}{BrickRed}&&%
    \myCVthree{5}{5}{5}{4}{4}{5}{Blue} &\myCVthree{2}{2}{3}{2}{2}{3}{Blue} &  \myCVthree[1212]{2}{2}{3}{2}{2}{3}{Blue}\\
  \midrule
  \multicolumn{11}{c}{Two triangles, one edge (Type 4)}\\
  \midrule
  \myCVthree{2}{3}{3}{3}{2}{4}{SeaGreen} & \myCVthree{5}{4}{5}{3}{4}{4}{SeaGreen} & \myCVtwo[1212]{3}{1}{2}{2}&&
  \myCVtwo{5}{0}{4}{1} & \myCVthree{2}{3}{3}{3}{2}{4}{SeaGreen} & \myCVthree[1212]{2}{3}{3}{3}{2}{4}{SeaGreen}&&
  \myCVtwo{5}{0}{4}{1} & \myCVthree{4}{3}{5}{3}{4}{4}{YellowOrange} & \myCVthree[1201]{4}{3}{5}{3}{4}{4}{YellowOrange}\\[5.0ex]
  \myCVtwo{5}{0}{4}{1} & \myCVthree{3}{4}{3}{3}{2}{4}{YellowOrange} & \myCVthree[1201]{3}{4}{3}{3}{2}{4}{YellowOrange}&&
  \myCVtwo{5}{0}{4}{1} & \myCVthree{5}{4}{5}{3}{4}{4}{SeaGreen} & \myCVthree[1201]{5}{4}{5}{3}{4}{4}{SeaGreen}\\
  \midrule
  \multicolumn{11}{c}{Three points (Type 1)}
  \\
  \midrule
   \myCVone{3}{3} & \myCVone{3}{3} & \myCVone[1212]{1}{1}\\
   \bottomrule
 \end{tabular}
 \vspace{5ex}

  \end{table}
\end{example}

\subsection{Diagonally constrained functions on $\P_q$}
\label{subsection:maximal-faces}
\citedinsurveyas{Section 4.2}

\begin{figure}[htbp]
\bgroup
\newcommand\DefinePosetNode[4]{
  \global\def#1{#2 & #3 & #4}}
\newcommand\PosetNode[2]
{#1\\
  \begin{tabular}{c@{\hspace{0.75em}}c@{\hspace{0.75em}}c}
    \midrule
    $p_1(F)$ & $p_2(F)$ & $p_3(F)$ \\
    \cmidrule{1-3}
    #2
  \end{tabular}}
\newcommand\PosetNodeDistance{2.8cm}
\newcommand\PosetInnerXsep{0pt}
\DefinePosetNode\onedallexample\EquiOneDimPointWithEnclosingEdge\EquiOneDimEdge\EquiOneDimEdge
\DefinePosetNode\onedgenericexample\EquiOneDimEdge\EquiOneDimEdge\EquiOneDimEdge
\DefinePosetNode\allexample\hor\diag\triup
\DefinePosetNode\diagexample\diag\triup\tridown
\DefinePosetNode\horizexample\hor\tridown\triup
\DefinePosetNode\vertexample\ver\triup\tridown
\DefinePosetNode\pointexample\triup\tridown\EquiPointNE
\DefinePosetNode\fullexample\triup\tridown\triup
\begin{tikzpicture} [ font = \small, align = flush center, >=stealth, thick,
  node distance = \PosetNodeDistance, 
  inner xsep = \PosetInnerXsep,
  every node/.style={draw, rounded corners}]
  \node [fill=red!30] (all) 
  {\PosetNode{arbitrarily\\constrained}\allexample};
  \node [fill=orange!30, below of = all] (diag) 
  {\PosetNode{full-dimensionally, \\ diagonally, and \\ point constrained}\diagexample}; 
  \node [fill=orange!30, left = 0.6cm of diag] (horiz) 
  {\PosetNode{full-dimensionally, \\ horizontally, and \\ point constrained}\horizexample};
  \node [fill=orange!30, right = 0.6cm of diag] (vert) 
  {\PosetNode{full-dimensionally, \\ vertically, and \\ point constrained}\vertexample};
  \node [fill=yellow!20, below of = diag] (point) 
  {\PosetNode{full-dimensionally\\ and point constrained}\pointexample};
  \node [fill=green!20, below = 0.6cm of point] (full) 
  {\PosetNode{full-dimensionally\\ constrained}\fullexample};
  \draw[->] (full) -- (point);
  \draw[->] (point) -- (diag);
  \draw[->] (diag) -- (all);
  \draw[->, shorten >= 0.5pt, shorten <= -1.5pt] (point) -- (horiz);
  \draw[->, shorten >= 0.5pt, shorten <= -1.5pt] (horiz) -- (all);
  \draw[->, shorten >= 0.5pt, shorten <= -1.5pt] (point) -- (vert);
  \draw[->, shorten >= 0.5pt, shorten <= -1.5pt] (vert) -- (all);
  \node [fill=yellow!20, left = 4cm of point] (onedall)
  {\PosetNode{full-dimensionally\\ and point constrained}\onedallexample};
  \node [fill=green!20, below = 0.6cm of onedall] (onedgeneric) 
  {\PosetNode{full-dimensionally\\ constrained}\onedgenericexample};
  \draw[->] (onedgeneric) -- (onedall);
\end{tikzpicture}
\egroup

  \caption{A hierarchy of minimal valid functions. 
    At the top is the most general, at the bottom the least general class of
    functions.  Each class is illustrated with the projections $p_1(F)$,
    $p_2(F)$, $p_3(F)$ of a maximal face
    $F\in\Delta\P$ with $F\subseteq E(\pi)$ that is
    allowed in this class, but not in the classes below.
    \emph{Left}, case $k=1$. 
    \emph{Right}, case $k=2$ for the standard triangulation~$\P_q$.
    }
  \label{fig:constrained-functions-poset}
\end{figure}
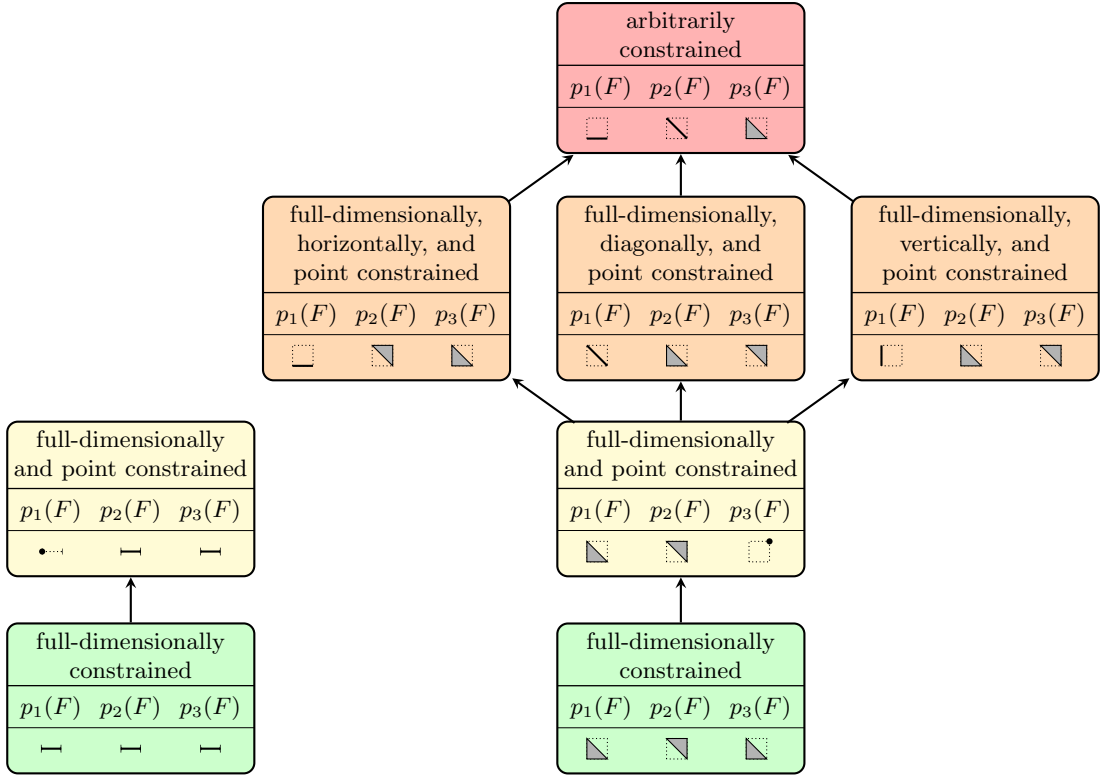

\autoref{fig:constrained-functions-poset} (on the right) shows a hierarchy of
minimal valid functions $\pi$ depending on the type of the possible projections $p_i(F)$ for maximal additive faces $F \in E_{\max{}}(\pi, \P_q)$.  
%
The labeling of the class is meant to be self-explanatory in \autoref{fig:constrained-functions-poset}; for example, in the lowest class ``full-dimensionally constrained'', all projections are 2-dimensional (triangles), in the class ``full-dimensionally and point constrained'', the projections are either 2-dimensional (triangles) or 0-dimensional (points), in the class ``full-dimensionally, horizontally and point constrained'' means the projections are either triangles, or horizontal edges, or points.

In this paper, we study the family of minimal valid functions that allows for two types of degenerations of the maximal additive faces, and characterize (in the sense of Theorems~\ref{thm:main} and~\ref{thm:1/4q}) the extreme functions within this family. Specifically, we assume that the maximal additive faces $F\in E_{\max{}}(\pi, \P_q)$ are so that its projections $p_i(F)$ are either full-dimensional (triangles $\triup, \tridown$), points
($\point$), or diagonal edges $(\diag)$, but not horizontal or vertical
edges.  These full-dimensionally, diagonally, and point constrained minimal
valid functions (\autoref{fig:constrained-functions-poset})  will be called
{\em diagonally constrained} minimal valid functions for brevity. 

\begin{definition}  A continuous piecewise linear function $\pi$ on $\P_q$ is
  called \emph{diagonally constrained} if whenever $F \in E_{\max{}}(\pi, \P_q)$, then $p_i(F) \in \Ipointdiagtri$ for $i=1,2,3$.
  %
\end{definition}

There are many examples of diagonally constrained functions. The unimodular properties of $\Delta \P_q$ provide an easy method to compute $E(\pi, \P_q)$ and test if
a function is diagonally constrained by using simple arithmetic and set
membership operations on vertices of~$\P_q$; see~\cite{hildebrand-thesis} for details. This can be done in polynomial time in $q$. 

%
%


\begin{example}[Example~\ref{ex:diag-constrained-function}, continued] \label{ex:diag-constrained-function-continued2}
  Since no relation appearing in the list of all maximal additive faces
  (\autoref{tab:example-emax}) involve a vertical or
  horizontal edge, the function is diagonally constrained.  Note that there
  are relations derived from two triangles and one diagonal edge.  These
  relations create affine properties as described in
  Figure~\ref{fig:higher-dim-interval-lemma}\,(b), and makes the analysis of this
  function more complicated than full-dimensionally constrained functions.
  \end{example}

\label{s:properties-valid-triples}

The following lemma characterizes the types of possible maximal additive faces that
can exist for a valid function that is diagonally constrained.  
\begin{lemma}\label{lemma:cases}
Suppose $\pi$ is continuous piecewise linear over $\P_q$ and is diagonally constrained.  Suppose that $F \in E_{\max}(\pi, \P_q)$. Let $I = p_1(F), J = p_2(F), K = p_3(F)$. Then one of the following is true.
\begin{enumerate}[(Type 1)]
\item\label{type:all-point/diag} $I,J,K \in \Ipointdiag$.
\item\label{type:all-tri}  $I, J,K \in \Itri$.
\item\label{type:tri,tri,point} One of $I,J,K$ is in $\Ipoint$, while the other two are in $\Itri$.
\item\label{type:tri,tri,edge} One of $I,J,K$ is in $\Idiag$, while the other two are in $\Itri$.
\end{enumerate}
\end{lemma}
All of these types of maximal additive faces appear in the function from
\autoref{ex:diag-constrained-function-continued2}: Maximal faces corresponding
to the symmetry condition are of Type~3, whereas Types 1, 2, and~4 appear in
Table~\ref{tab:example-emax}.
\begin{proof}
By definition of diagonally constrained functions, $I,J,K \in \Ipoint\cup
\Idiag \cup \Itri$. 
Elementary counting reveals that there are 27 possible ways to put $I,J,K$ into those three sets, whereas 15 possibilities are described above.  We will show that the 12 remaining cases not listed above are not possible because $I,J,K$ are projections of $F$. 
\begin{enumerate}
\item Suppose $I, J \in \Ipointdiag$, $K \in \Itri$.   By Proposition~\ref{prop:projection},  $K \subseteq I+J$.  But this is not possible because $I+J$ is one-dimensional while $K$ is two-dimensional. 
\item Suppose $I, K\in \Ipointdiag$, $J \in \Itri$.  By Proposition~\ref{prop:projection}, $J \subseteq K + (- I)$.  But again, this is not possible because $K + (-I)$ is one-dimensional while $J$ is two-dimensional.
\item Suppose $J, K\in \Ipointdiag$, $I \in \Itri$. This is similar to the last case.\qedhere
\end{enumerate}
\end{proof}

\subsection{Affine properties of $\pi^i$ on projections of faces in $E(\pi, \P_q)$}\label{sec:affine-lin}
Let $\pi$ be a minimal valid function that is continuous piecewise linear over~$\P_q$.
The lemmas of this subsection will be used to deduce affine properties of valid functions $\pi^1, \pi^2$ when $\pi = \tfrac12 (\pi^1 + \pi^2)$ by using Lemma~\ref{lemma:tight-implies-tight}.  Here we will apply Corollary~\ref{lem:projection_interval_lemma-corollary} to conclude affine properties on faces of $\P_q$.  By using Corollary~\ref{lem:projection_interval_lemma-corollary}, we are using the continuity of the function to extend affine properties to the boundaries of faces.  

\begin{lemma}\label{cor:triangle+triangle}
Suppose $\theta\colon\R^2 \to \R$ is a continuous function and let $F \in E(\theta, \P_q)$ such that $p_i(F) \in \Itri$ for $i=1,2,3$. Then $\theta$ is affine in $p_i(F)$ for $i=1,2,3$ with the same gradient.
\end{lemma}
\begin{proof}
We apply Corollary~\ref{lem:projection_interval_lemma-corollary} to $F$ with $f,g,h = \theta$ and $L =
\R^2$.  Since $p_1(F), p_2(F) \in \Itri$, and triangles are two-dimensional objects, we have $L\times L + F \subseteq \aff(F)$. The conclusion of the corollary then says
that $\theta$ is affine over $p_i(F)$  for $i=1,2,3$ with the same gradient.  
\end{proof}

\begin{lemma}\label{cor:triangle+diagonal}
Let $\theta\colon \R^2 \to \R$ be a continuous function. Let $F \in E(\theta,\P_q)$ such that  $p_1(F), p_2(F) \in \Itri$ and $p_3(F) \in \Iedge$ (resp., $p_1(F), p_3(F) \in \Itri$ and $p_2(F) \in \Iedge$).
Let $L$ be the linear space such that $\aff(p_3(F))$ (resp., $\aff(p_2(F))$) is a translate of $L$. 
Then for some $\cve \in \R^2$, $\theta$ is affine with respect to $L$ over $p_1(F), p_2(F), p_3(F)$ with the gradient $\cve$.
\end{lemma}

\begin{proof}

We only give the proof for $p_1(F), p_2(F) \in \Itri$ and $p_3(F) \in \Iedge$. The other case is similar. 

Consider any $(\u^1, \v^1), (\u^2, \v^2) \in \relint(F)$. By applying
Corollary~\ref{lem:projection_interval_lemma-corollary} with $F$ and $L$ we see that there
exist vectors $\cve^1, \cve^2\in \R^k$ such that $\theta$ is affine with
gradient $\cve^i$ over $(\u^i + L) \cap p_1(F)$, $(\v^i + L)
\cap p_2(F)$ and $(\u^i + \v^i + L) \cap p_3(F)$ for $i = 1,2$.  Let $\bar \cve^1, \bar \cve^2$ be the orthogonal projections of $\cve^1$ and $\cve^2$, respectively, onto the linear space $L$.  Therefore, $\theta$ is affine with
gradient $\bar \cve^i$ over $(\u^i + L) \cap p_1(F)$, $(\v^i + L)
\cap p_2(F)$ and $(\u^i + \v^i + L) \cap p_3(F)$ for $i = 1,2$.
Then, since $(\u^1 + \v^1 + L) \cap p_3(F) = (\u^2 + \v^2 + L) \cap
p_3(F) = p_3(F)$, we have $\bar \cve^1 = \bar \cve^2$.  
Therefore, we obtain that $\theta$ is affine with respect to $L$ with gradient $\bar \cve = \bar \cve^1 = \bar \cve^2$ over $p_i(F)$ for $i=1,2,3$. 
\end{proof}

\begin{definition}
\label{def:Ldiag}
Define 
      \begin{displaymath}
    L_{\diag}
    = \{\, \ve
    x\in\R^2\st \ve 1\cdot\ve x = 0\,\} = \{\,\lambda \ColVec{-1}{1} \st \lambda \in
    \R\,\}.
  \end{displaymath}
\end{definition}

\begin{lemma}[Geometric adjacent transference] \label{obs:adjacent}
Let $I,J \in \Itri$ be triangles such that $I\cap J \in \Iverthor$.  Let $\pi$ be a continuous function defined on $I\cup J$ satisfying the following properties:
\begin{itemize}
\item[(i)] $\pi$ is affine on $I$.
\item[(ii)] $\pi$ is affine with respect to the linear space  $L_{\diag}$
(the diagonal direction) on $J$. 
\end{itemize}
Then $\pi$ is affine on $J$.
\end{lemma}
\begin{proof}
Let $e=I\cap J\in \Iverthor$ be the common edge of $I$ and $J$. We assume that $e$ is
vertical (the argument for horizontal edges is exactly the same) and let $\v^0
\in \R^2$ be the vertex of $e$ such that the other vertex is $\v^0 +
\ColVec[q]{0}{1}$. Since $\pi$ is affine with respect to the linear space
$ L_{\diag}$ on $J$, there exists $c \in \R$ such that $\pi(\x + \lambda \ColVec{-1}{1}) = \pi(\x) + c\cdot \lambda$ for all $\x \in J$ and $\lambda \in \R$ such that $\x + \lambda \ColVec{-1}{1} \in J$. Since $\pi$ is affine on $I$, there exists $c' \in \R$ such that $\pi(\v^0 + \lambda \ColVec{0}{1}) = \pi(\v^0) + c' \cdot\lambda$ for all $0 \leq \lambda \leq \tfrac{1}{q}$. 

 Now observe that any point in $J$ can be written as $\v^0 + \mu_1  \ColVec{0}{1} + \mu_2  \ColVec{-1}{1}$ with $0 \leq \mu_1, \mu_2 \leq \tfrac{1}{q}$ and therefore, $\pi(\v^0 + \mu_1  \ColVec{0}{1} + \mu_2  \ColVec{-1}{1}) = \pi(\v^0 + \mu_1  \ColVec{0}{1}) + c\cdot\mu_2$ (using (ii) in the hypothesis) and $\pi(\v^0 + \mu_1  \ColVec{0}{1}) + c\cdot\mu_2 = \pi(\v^0) + c'\cdot\mu_1 + c\cdot\mu_2$. Thus, $\pi$ is affine over $J$.
\end{proof}

\section{Proof of the main results for the two-dimensional case}
\label{sec:main-results}
In this section we prove our main results for continuous piecewise linear
functions over $\P_{q}$.  

\begin{assumption}
For the remainder of the paper, we assume that $\pi$ is a minimal valid function that is continuous piecewise linear over $\P_q$.  
\end{assumption}



\begin{definition}
\label{def:affine-imposing}
  \begin{enumerate}[\rm(a)]
  \item For any $I \in \P_q$, if $\pi$ is affine in
    $I$ and if for all valid functions $\pi^1, \pi^2$ such that $\pi =
    \tfrac{1}{2}( \pi^1 +  \pi^2)$ we have that $\pi^1, \pi^2$ are
    affine in $I$, then we say that $\pi$ is \emph{affine imposing in
      $I$}.
\item      For any $I \in \P_q$, if $\pi$ is affine with respect to $L_{\diag}$ over $I$ and if for all valid functions $\pi^1, \pi^2$ such that $\pi =
    \tfrac{1}{2}( \pi^1 + \pi^2)$ we have that $\pi^1$, $\pi^2$ are both
    affine with respect to $L_{\diag}$ over $I$, then we say that $\pi$ is \emph{diagonally affine imposing in
      $I$}.
  \item For a collection $\I\subseteq \P_q$, if for all
    $I \in \I$, $\pi$ is affine imposing (or diagonally affine imposing) in $I$, then we say that $\pi$ is
    \emph{affine imposing (diagonally affine imposing) in $\I$}.
  \end{enumerate}
\end{definition}

\paragraph{\bf Section outline.} We either show that $\pi$ is affine imposing in $\P_q$
(subsection~\ref{section:AI}) or construct a continuous piecewise linear perturbation (subsection~\ref{s:eq-perturb}) that proves $\pi$ is not extreme (subsection \ref{sec:non-extreme-by-perturbation}).
If $\pi$ is affine imposing in $\P_q$, we set up a system of linear equations to decide if~$\pi$ is extreme or not (subsection~\ref{section:system}). This implies Theorem~\ref{thm:main} stated in the introduction.

\subsection{Imposing affine linear properties on faces of~$\P_q$}
\label{section:AI}
As briefly discussed in~\autoref{sec:affine-lin}, the set $E(\pi,\P_q)$ helps one to deduce affine linear properties of $\pi^1, \pi^2$. There are essentially three types of such deductions: (i) Lemma~\ref{cor:triangle+triangle} and Lemma~\ref{cor:triangle+diagonal} (deducing affine linear properties by Interval Lemma type arguments), (ii) transferring affine linear properties through lower dimensional faces, and (iii) Lemma~\ref{obs:adjacent}, which transfers affine linear properties via adjacency between cells of $\P_q$. We build a finite graph to formally record these interactions in section~\ref{sec:finite-graph}.

\subsubsection{Covered triangles.}
We now consider faces of $\Itri$ on which we can deduce affine properties.  
$$
\begin{array}{r}
\Itri^1 = \{\,  I,J \in \Itri \st \exists K \in \Idiag, F \in E(\pi, \P_q) \text{ with } (I,J,K) = (p_1(F), p_2(F), p_3(F)) \;\;\;\;\\
\text{ or  } (I,K,J) =(p_1(F), p_2(F), p_3(F)) \,\},
\end{array}
$$
$$
\Itri^2 = \{\, I,J,K\in \Itri \st \exists F \in E(\pi, \P_q) \text{ with } (I,J,K) = (p_1(F), p_2(F), p_3(F)) \,\}.
$$
It follows from Lemma~\ref{cor:triangle+triangle} that $\pi$ is affine imposing in $\Itri^2$ and from Lemma~\ref{cor:triangle+diagonal} that $\pi$ is diagonally affine imposing in $\Itri^1$.  The superscripts here correspond to the dimension of the linear space on which $\pi$ is affine imposing on a face.

\subsubsection{Finite graph.}\label{sec:finite-graph}
Next we will define a finite graph~$\G$ whose nodes correspond to the
two-dimensional faces (triangles) in $\Itri$.  To make this graph finite, we
will use the periodicity of the function~$\pi$
and of the complex~$\P_q$ modulo~$\Z^2$.  
By $\Itri/\Z^2$ we denote the set of
equivalence classes 
\begin{equation}
  \EqClass{I} = \{\, \tau_{\ve s}(I) = I + \ve s \st \ve s\in\Z^2\,\}
  \label{eq:equivalence-class}
\end{equation}
of two-dimensional faces (triangles) $I \in\Itri$ modulo translations by
integer vectors~$\ve s\in\Z^2$.  We can identify an 
equivalence class with its unique representative that is a triangle contained
in~$[0,1]^2$. 

\begin{definition}
\label{def:graph}
  Let $\G = \G(\Itri/\Z^2,\E)$ be the finite undirected graph with node set
  $\Itri$ and edge set $\E = \E_\point \cup \E_\diag$ where $\{\EqClass{I},
  \EqClass{J}\} \in \E_\point$ (resp., $\{\EqClass{I}, \EqClass{J}\} \in
  \E_\diag$) if and only if $[I]\neq[J]$ and for some $K \in \Ipoint$ (resp.,
  $ K \in \Idiag$) and $F \in E(\pi, \P_q)$, 
  we have one of the following cases:
  \begin{enumerate}[({Case} a.)]
  \item  \label{case-a} $(I,J,K) = (p_1(F), p_2(F), p_3(F))$, which implies $F' := F(J,I,K) \in
    E(\pi,\P_q)$ with $(J,I,K) = (p_1(F'), p_2(F'), p_3(F'))$, or 
  \item  \label{case-b} $(I,K,J) = (p_1(F), p_2(F), p_3(F))$, or
  \item \label{case-c} $(J,K,I) = (p_1(F), p_2(F), p_3(F))$.
  \end{enumerate}
\end{definition}
Therefore we record an edge between two cells in $\Itri$ whenever there is an $F \in E(\pi, \P_q)$ such that two of the projections $p_i(F)$, $i=1,2,3$, are these two cells and the third projection is in $\Ipointdiag$. By the symmetry between $p_1$ and $p_2$ and the symmetry in the definition of $E(\pi, \P_q)$, for every $F \in E(\pi, \P_q)$ there exists an $F' \in E(\pi, \P_q)$ such that $p_1(F) = p_2(F')$, $p_2(F) = p_1(F')$, and $p_3(F) = p_3(F')$.  Therefore, when considering an $F \in E(\pi, \P_q)$ with two projections in $I,J\in \Itri$ and a third projection $K \in \Ipointdiag$, we can always assume that either $p_2(F) = K$ or $p_3(F)= K$.

Some faces in $\E_\diag$ are inherently also in $\E_\point$.  \autoref{fig:why-diagonals} depicts how this can happen and also shows an edge in $\E_\diag$ that is not necessarily in $\E_\point$.  Thus, $\E_\point$ alone is not sufficient to describe all the relations in the graph that we need to consider.

\begin{figure}
 $$
 \myCVfigSevenOne \qquad\qquad 
 \myCVfigSevenTwo \qquad\qquad 
 \myCVfigSevenThree
$$

%
%
%

  \caption{An example of an important edge connection in $\E_\diag$ that is
    not captured with $\E_\point$.  
    For a given minimal valid
    function~$\pi$, we could have $F(I,J,K_1), F(I,J,K_2), F(I,J,K_3) \in
    E(\pi, \P_q)$.  Therefore,
    $\{\EqClass{J},\EqClass{K_1}\},\{\EqClass{J},\EqClass{K_2}\},\{\EqClass{J},\EqClass{K_3}\}
    \in \E_\diag$. Thus, $[J], [K_1], [K_2], [K_3]$ are connected in the graph
    $\G(\Itri/\Z^2,\E)$.  Notice however that 
    already the smaller faces
    $F(I_1, J, K_1)\subseteq F(I,J,K_1)$ and
    $F(I_2, J,K_2) \subseteq F(I,J,K_2)$ (corresponding to the vertices $I_1$
    and $I_2$ of the one-dimensional face~$I$) 
    induce $F(I_1, J,K_1), F(I_2, J,K_2) \in E(\pi,\P_q)$.  Thus
    $\{\EqClass{J},\EqClass{K_1}\},\{\EqClass{J},\EqClass{K_2}\} \in
    \E_\point$, and so we can realize these two edges in the graph without
    using $\E_\diag$.  But the third edge, $\{\EqClass{J},\EqClass{K_3}\}$, is
    only realized in $\E_\diag$.  Therefore, recording $\E_\diag$ is crucial
    to the construction of $\G(\Itri/\Z^2, \E)$. 
}\label{fig:why-diagonals}
\end{figure}

The functions $\pi, \pi^1, \pi^2$ have related slopes on faces that are connected in the graph.

\begin{lemma}
\label{lemma:point-and-line-Lemma}
Let $L \subseteq \R^2$ be a linear subspace. Let $I, J \in \Itri$ with $\{\EqClass{I},\EqClass{J}\} \in \E$. Suppose $\pi^1, \pi^2$ are valid functions with $\pi = \tfrac12(\pi^1 + \pi^2)$.  For $\theta = \pi, \pi^1$, or $\pi^2$, 
if $\theta$ is affine with respect to $L$ over $I$, 
then $\theta$ is affine with respect to $L$ over~$J$ as well.
\end{lemma}
\begin{proof}
By Lemma \ref{lem:tightness} and Assumption 1, $\pi^1, \pi^2$ are minimal and continuous and $E(\pi, \P_q)\subseteq E(\theta, \P_q)$ for $\theta = \pi,\pi^1, \pi^2$.

\emph{Case (i).}  Suppose $\{\EqClass{I},\EqClass{J}\}\in \E_\point$.  
Then there exists $\a \in \frac{1}{q}\Z^2$ such that, setting $K = \{\a\} \in \Ipoint$, there exists $F \in E(\pi, \P_q)$ such that either $(I,J,K)= (p_1(F), p_2(F), p_3(F))$, $(I,K,J) = (p_1(F), p_2(F), p_3(F))$, or $(J,K,I) = (p_1(F), p_2(F), p_3(F))$; these are cases \ref{case-a}, \ref{case-b}, and \ref{case-c} from Definition~\ref{def:graph}, respectively. We only consider the case $(I,J,K) = (p_1(F), p_2(F), p_3(F))$; the other cases are similar. Then $\theta|_I(\x) + \theta|_J(\y) = \theta|_K(\a)$ for all $\x \in I$, $\y \in J$, $\x + \y = \a$. Consider any $\y^1, \y^2 \in J$ such that $\y^2 - \y^1 \in L$. Set $\x^i = \a - \y^i \in I$ for $i=1,2$. Thus, $\theta|_J(\y^2) - \theta|_J(\y^1) = \theta|_I(\x^1) - \theta|_I(\x^2)$ and $\x^1 - \x^2 = \y^2 - \y^1 \in L$. Since $\theta$ is affine with respect to $L$ over $I$, $\theta$ is affine with respect to $L$ over $J$.

\emph{Case (ii).} Suppose $\{\EqClass{I},\EqClass{J}\}\in \E_\diag$.   We show that for any $\y \in \relint(J)$, and any $\p \in L$ there exists $\epsilon > 0$ such that $\theta$ is affine over $\{\,\y + \lambda\p\st -\epsilon \leq \lambda \leq \epsilon\,\}$. Using Lemma~\ref{lem:patching}, this will then imply that $\theta$ is affine with respect to $L$ over $J$. We only consider the case when there exists $F \in E(\pi, \P_q)$ with $(I,J,K) =(p_1(F), p_2(F), p_3(F))$ (case \ref{case-a} from Definition~\ref{def:graph}); the cases~\ref{case-b} and \ref{case-c} are similar.  Thus, $I,J \in \Itri$ and $K \in \Idiag$.  Then $\theta|_I(\x) + \theta|_J(\y) = \theta|_K(\a)$ for all $\x \in I$, $\y \in J$, $\x + \y = \a\in K$.

Let $\y\in \relint(J)$.  Using \autoref{lem:rel-int-7-tuple}, there exists $\x \in \relint(I)$ and $\a \in \relint(K)$ such that $\x + \y = \a$. Since $\y \in \relint(J)$, there exists $\epsilon > 0$ such that $\{\,\x + \lambda\p\st -\epsilon \leq \lambda \leq \epsilon\,\} \subseteq I$ and $\{\,\y + \lambda\p\st -\epsilon \leq \lambda \leq \epsilon\,\} \subseteq J$. Since $\theta$ is affine over $\{\,\x + \lambda\p\st -\epsilon \leq \lambda \leq \epsilon\,\} \subseteq I$ and $\theta|_I(\x) + \theta|_J(\y) = \theta|_K(\a)$ for all $\x \in I$, $\y \in J$, $\x + \y = \a$, a similar argument as case (i) proves $\theta$ is affine over $\{\,\y + \lambda\p\st -\epsilon \leq \lambda \leq \epsilon\,\}$.\end{proof}

With this in mind, 
for each $I \in \Itri$, let $\G_I$ be the connected component of $\G$
containing $\EqClass{I}$.  We define the two sets of faces that contain complete connected components in the graph $\G$,
$$
\Stri^1 =\bigl\{\,J \in \Itri \bigst \EqClass{J} \in \G_I \text{ for some } I \in \Itri^1\,\bigr\},
$$
$$
\Stri^2 =\bigl\{\,J \in \Itri \bigst \EqClass{J} \in \G_I \text{ for some } I \in \Itri^2\,\bigr\}.
$$
\begin{obs}
\label{obs:affine-imposing}
It follows from Lemma~\ref{lemma:point-and-line-Lemma}, Lemma~\ref{cor:triangle+triangle} and the periodicity of $\pi$, $\pi^1$, and $\pi^2$ that $\pi$ is affine imposing in $\Stri^2$.   Similarly, it follows from 
Lemma~\ref{lemma:point-and-line-Lemma}, Lemma~\ref{cor:triangle+diagonal} and the periodicity of $\pi$, $\pi^1$, and $\pi^2$
that $\pi$ is diagonally affine imposing in $\Stri^1$.
\end{obs}  
\begin{obs}[Geometrically adjacent triangles]
\label{obs:adjacent-apply-lemma}
From Lemma~\ref{obs:adjacent}, it follows that if $I \in \Stri^2$, $J \in \Stri^1$ and $I\cap J \in \Iverthor$, then $\pi$ is affine imposing in $J$.  Furthermore, by periodicity of $\pi$, $\pi^1$, and $\pi^2$, $\pi$ is affine imposing in all $J' \in \EqClass{J}$.
\end{obs}
This observation motivates the following graph definition that is a super-graph of $\G$.
\begin{definition}
  Let $\bar\G = \bar\G(\Itri/\Z^2,\bar\E)$ be the finite undirected graph with node set
  $\Itri/\Z^2$ and edge set $\bar \E = \E_\point \cup \E_\diag \cup \E_{\ver\hor}$ where $\E_\point$ and $\E_\diag$ are defined in Definition~\ref{def:graph} and where $\{\EqClass{I},
  \EqClass{J}\} \in \E_{\ver \hor}$  if and only if $[I]\neq[J]$ and for some $I' \in \EqClass{I}$, $J' \in \EqClass{J}$ we have $I',J' \in \Stri^1 \cup \Stri^2$ and $I \cap J \in \Iverthor$.
\end{definition}

In contrast to the graph $\G$ and Lemma~\ref{lemma:point-and-line-Lemma}, faces in $\bar\G$ connected by edges from $\E_{\ver\hor}$ do not necessarily have related slopes, even if $\pi$ is affine imposing on these faces.

For each $I \in \Itri$, let $\bar \G_I$ be the connected component of $\bar \G$ containing $\EqClass{I}$.
 Let 
 \begin{equation}
\barStri^2 = \bigl\{\, K \in \Itri \bigst \EqClass{K} \in \bar \G_I \text{ for some }  I \in \Stri^2 \,\bigr\}.
\end{equation}

Note that $\barStri^2 \subseteq \Stri^1 \cup \Stri^2$.
Let \begin{equation} 
\barStri^1 = \Stri^1 \setminus \barStri^2.
\end{equation}

\begin{theorem}
\label{theorem:AI}
If $\barStri^2 = \Itri$, then $\pi$ is affine imposing in $\Itri$, and therefore $\theta$ is continuous piecewise linear over $\P_q$ for $\theta = \pi^1, \pi^2$ whenever we have that $\pi^1, \pi^2$  are valid functions and $\pi = \tfrac12(\pi^1 + \pi^2)$.  
\end{theorem}
\begin{proof}
By Lemma \ref{lem:tightness}, $\pi^1, \pi^2$ are minimal and continuous.
Since they are minimal, they are also periodic.  From
Observation~\ref{obs:affine-imposing}, $\pi$ is affine imposing in $\Stri^2$
and diagonally affine imposing in $\Stri^1$.   By
Observation~\ref{obs:adjacent-apply-lemma}, $\pi$ is affine imposing in any
$J'$ such that there exists a $J$ with $J \in [J']$ and $I \cap J \in
\P_{q,\ver\hor}$ for some $I$ such that $\pi$ is affine imposing in $I$.  In
particular, this holds for all $I \in \Stri^2$.   Consider any  $K\in \Stri^1$
where $\EqClass{K}$ is connected by a path to $\EqClass{I}$ in the graph
$\bar\G$.  By induction on the number of edges in the path from $\EqClass{K}$ to $\EqClass{I}$ and using Lemma~\ref{lemma:point-and-line-Lemma},  $\pi$ is affine imposing in $K$.  Therefore, $\pi$ is affine imposing in $\Stri^1 \cap \barStri^2$.  Since $\Itri = \barStri^2  \subseteq \Stri^1 \cup \Stri^2 \subseteq \Itri$, it follows that $\pi$ is affine imposing in all of $\Itri$.  
\end{proof}


\subsection{Perturbation functions}\label{s:eq-perturb}

In this section we study functions $\psi \colon \R^2 \to \R$ that satisfy entire classes of additivity relations that appear in $E(\pi, \P_q)$. These will be used to construct perturbation functions $\bpi$ such that $E(\pi) \subseteq E(\bpi)$. We may then leverage \autoref{corPerturb} to show that $\pi$ is not extreme.  

\newcommand\DrawAxisLabels{%
  \node[anchor=north] at (-1cm, -1cm) {$-\tfrac1q$};
  \node[anchor=north] at (0cm, -1cm) {$0\vphantom{\tfrac1q}$};
  \node[anchor=north] at (1cm, -1cm) {$\tfrac1q$};
  \node[anchor=north] at (2cm, -1cm) {$\tfrac2q$};
  \node[anchor=east] at (-1cm, -1cm) {$-\tfrac1q$};
  \node[anchor=east] at (-1cm, 0cm) {$0$};
  \node[anchor=east] at (-1cm, 1cm) {$\tfrac1q$};
  \node[anchor=east] at (-1cm, 2cm) {$\tfrac2q$};
}

\newcommand\DrawTriangle{
  \fill[ultra thick] (0,0) -- (0,1) -- (1,0) -- cycle;  
  \draw[ultra thick] (1,0) -- (0,0) -- (0,1);
  \draw[ultra thick] (0,1) -- (1,0);
  \visible<8-> {
    \node at (0.25,0.25) {\color{magenta}$+$};
  }
}

\newcommand\DrawUpperTriangle{
  \fill[ultra thick] (1,1) -- (0,1) -- (1,0) -- cycle;  
  \draw[ultra thick] (0,1) -- (1,1) -- (1,0);
  \draw[ultra thick] (0,1) -- (1,0);
  \visible<8-> {
    \node at (0.75,0.75) {\color{magenta}$-$};
  }
}

\newcommand\DrawTwodimGrid{
\begin{tikzpicture}[fill=white,>=stealth]
        \visible<2->{
          \DrawUpperTriangle
        }
        \visible<3->{
          \begin{scope}[yshift=1cm] 
            \DrawTriangle
          \end{scope}
        }
        \visible<5->{
          \foreach \x in {-1cm, 0cm, 1cm} {
            \foreach \y in {-1cm, 0cm, 1cm} {
              \begin{scope}[xshift=\x,yshift=\y] 
                \DrawTriangle
                \DrawUpperTriangle
              \end{scope}
            }
          }
        }
        \visible<1->{
          \begin{scope}[fill=cyan!50,draw=cyan]
            \DrawTriangle
          \end{scope}
        }
        \visible<2-6>{
        }
        %
        
        
        \DrawAxisLabels
      \end{tikzpicture}      
}


\newcommand\DrawTrianglePoint{
  \visible<3->{%
%

\fill[ fill=white] (0,0) -- (1/3,0) -- (0,1/3) -- cycle; 
  \fill[ fill=white] (1,0) -- (2/3,0) -- (2/3,1/3) -- cycle; 
  \fill[ fill=white] (0,1) -- (0,2/3) -- (1/3,2/3) -- cycle; 
  \shade[top color=green, bottom color=white] (1/3,0)  -- (1/3,1/3) --
  (2/3,0) -- cycle;
  \shade[ right color=green, left color=white] (0,1/3) -- (1/3,1/3) -- (1/3,2/3) --
  (0,2/3) -- cycle;
  \shade[left color=green, right color=white] 
  (1/3,1/3) -- (2/3,0) -- (2/3,1/3) -- cycle;
  \shade[ bottom color=green, top color=white] 
  (1/3,1/3) -- (0,2/3) -- (1/3,2/3) -- cycle;
  
   \draw[thin, draw=black] (1/3,0) -- (0,1/3); 
   \draw[thin, draw=black] (2/3,0) -- (0,2/3);
   \draw[thin, draw=black] (1/3,0) -- (1/3,2/3);
   \draw[thin, draw=black] (2/3,0) -- (2/3,1/3);
   \draw[thin, draw=black] (0,1/3) -- (2/3,1/3);
   \draw[thin, draw=black] (0,2/3) -- (1/3,2/3);
   
      \draw[ultra thick, draw=black] (0,0) -- (0,1) -- (1,0) -- cycle; 
  \draw[ultra thick, draw=black] (0,0) -- (0,1) -- (1,0) -- cycle; 

  }  
}

\newcommand\DrawUpperTrianglePoint{
  \visible<3->{%
    \fill[ fill=white] (1,1) -- (2/3,1) -- (1,2/3) -- cycle; 
  \fill[ fill=white] (0,1) -- (1/3,1) -- (1/3,2/3) -- cycle; 
  \fill[ fill=white] (1,0) -- (1,1/3) -- (2/3,1/3) -- cycle; 

\shade[ bottom color=red, top color=white] (2/3,1) -- (2/3,2/3) -- (.66,2/3) --
  (1/3,1) -- cycle;
  \shade[ left color=red, right color=white] (1,2/3) -- (2/3,2/3) -- (2/3,.66) --
  (1,1/3) -- cycle;
  \shade[ right color=red, left color=white] 
  (0.66,2/3) -- (1/3,1) -- (1/3,02/3) -- cycle;
  \shade[ top color=red, bottom color=white] 
  (2/3,0.66) -- (1,1/3) -- (2/3,1/3) -- cycle;

   \draw[thin, draw=black] (1/3,1) -- (1,1/3); 
   \draw[thin, draw=black] (2/3,1) -- (1,2/3);
   \draw[thin, draw=black] (1/3,1) -- (1/3,2/3);
   \draw[thin, draw=black] (2/3,1) -- (2/3,1/3);
   \draw[thin, draw=black] (1,1/3) -- (2/3,1/3);
   \draw[thin, draw=black] (1,2/3) -- (1/3,2/3);
   
  
      \draw[ultra thick, draw=black] (1,1) -- (0,1) -- (1,0) -- cycle; 
  \draw[ultra thick, draw=black] (0,0) -- (0,1) -- (1,0) -- cycle;

}
}

\newcommand\DrawTwodimGridPointPerturbation{

    \begin{tikzpicture}[transform canvas={rotate=-45}, xshift = -1cm, yshift= -1cm]

    { 

    \begin{scope}[rotate = 45]
    \begin{scope}[yshift = 1.7cm]
     \foreach \x in {-1,0, 1} {
       \foreach \y in {-1,0, 1} {
   \shade [top color=white,bottom color=green ] (1/3+\x,1/3+\y) -- (1/3+\x ,2/3+\y) -- (2/3+\x,1/3+\y)  -- cycle;
   \shade [top color=green,bottom color=white ] (1/3+\x,1/3+\y) -- (1/3+\x ,0+\y) -- (0+\x,1/3+\y)  -- cycle;
     \shade [top color=white,bottom color=red ] (2/3+\x,2/3+\y) -- (2/3+\x ,1+\y) -- (1+\x,2/3+\y)  -- cycle;
   \shade [top color=red,bottom color=white ] (2/3+\x,2/3+\y) -- (2/3+\x ,1/3+\y) -- (1/3+\x,2/3+\y)  -- cycle;
   }
   }
   \end{scope}
   \end{scope}

  }
\end{tikzpicture}
\hspace{-3.39cm}
\begin{tikzpicture}[fill=white,>=stealth]
        \visible<2->{
          \foreach \x in {-1cm, 0cm, 1cm} {
            \foreach \y in {-1cm, 0cm, 1cm} {
              \begin{scope}[xshift=\x,yshift=\y] 
                \DrawTrianglePoint
                \DrawUpperTrianglePoint
              \end{scope}
            }
          }
        }

        \DrawAxisLabels

      \end{tikzpicture}      
}

\newcommand\DrawTriangleDiag{
  \filldraw[thick, draw=black] (0,0) -- (0,1) -- (1,0) -- cycle;  
  \visible<5->{
    \shadedraw[draw=black, right color=white, left color=red,shading angle=-45]
    (0.5,0) -- (0,0.5) -- (0,2/3) -- (2/3,0) -- cycle;
    \shadedraw[draw=black, right color=red, left color=white,shading angle=-45]
    (0,2/3) -- (2/3,0) -- (1,0) -- (0,1) -- cycle;
  }
  \visible<4->{
    \draw[ultra thick,draw=green] (1/3,0) -- (0,1/3);
    \draw[ultra thick,draw=white] (0.5,0) -- (0,0.5);
    \draw[ultra thick,draw=red] (0,2/3) -- (2/3,0);
    \draw[ultra thick,draw=white] (1,0) -- (0,1);
  }

}

\newcommand\DrawUpperTriangleDiag{
  \filldraw[thick, draw=black] (1,1) -- (1,0) -- (0,1) -- cycle;  
  \visible<5->{
    \shadedraw[draw=black, left color=white, right color=red,shading angle=-45]
    (1,1) -- (2/3,1) -- (1,2/3) -- cycle;
    \shadedraw[draw=black, left color=red, right color=white,shading angle=-45]
    (2/3,1) -- (1,2/3) -- (1,0.5) -- (0.5,1) -- cycle;
    \shadedraw[draw=black, left color=white, right color=green,shading angle=-45]
    (0.5,1) -- (1,0.5) -- (1,1/3) -- (1/3,1) -- cycle;
    \shadedraw[draw=black, left color=green, right color=white,shading angle=-45]
    (1,1/3) -- (1/3,1) -- (0,1) -- (1,0) -- cycle;
  }
   \visible<4->{
     \draw[ultra thick,draw=red] (2/3,1) -- (1,2/3);
     \draw[ultra thick,draw=white] (0.5,1) -- (1,0.5);
     \draw[ultra thick,draw=green] (1,1/3) -- (1/3,1);
     \draw[ultra thick,draw=white] (1,0) -- (0,1);
   }
}

\newcommand\DrawTwodimGridDiag{
  \begin{tikzpicture}[transform canvas={rotate=45},fill=white,>=stealth]   

  \visible<4->{

    \begin{scope}[rotate = -45, xshift =0 cm, yshift = 0.7cm]

    \begin{scope}[xshift = 1cm, yshift = 1cm]
             \node[anchor=north, rotate = -45] at (-1cm, -1cm) {$-\tfrac{1}{q}$};
  \node[anchor=north, rotate = -45] at (0cm, -1cm) {$0\vphantom{\tfrac1q}$};
  \node[anchor=north, rotate = -45] at (1cm, -1cm) {$\tfrac1q$};
  \node[anchor=north, rotate = -45] at (2cm, -1cm) {$\tfrac2q$};
  \node[anchor=east, rotate = -45] at (-1cm, -1cm) {$-\tfrac1q$};
  \node[anchor=east, rotate = -45] at (-1cm, 0cm) {$0$};
  \node[anchor=east, rotate = -45] at (-1cm, 1cm) {$\tfrac1q$};
  \node[anchor=east, rotate = -45] at (-1cm, 2cm) {$\tfrac2q$};
  \end{scope}

      \foreach \x in {0, 1, 2} {
    \shade[ right color=green, left color=white]
   (0+\x,0) -- (0,0+\x) -- (0,1/3+\x) -- (1/3+\x,0) -- cycle;
    \shade[ right color=white, left color=green]
    (1/3+\x,0) -- (0,1/3+\x) -- (0,0.5+\x) -- (0.5+\x,0) -- cycle;
    \shade[ right color=red, left color=white]
    (0.5+\x,0) -- (0,0.5+\x) -- (0,2/3+\x) -- (2/3+\x,0) -- cycle;
    \shade[right color=white, left color=red]
    (0,2/3+\x) -- (2/3+\x,0) -- (1+\x,0) -- (0,1+\x) -- cycle;
    
    \shade[ right color=green, left color=white]
   (0+\x,3) -- (3,0+\x) -- (3,1/3+\x) -- (1/3+\x,3) -- cycle;
    \shade[ right color=white, left color=green]
    (1/3+\x,3) -- (3,1/3+\x) -- (3,0.5+\x) -- (0.5+\x,3) -- cycle;
    \shade[ right color=red, left color=white]
    (0.5+\x,3) -- (3,0.5+\x) -- (3,2/3+\x) -- (2/3+\x,3) -- cycle;
    \shade[ right color=white, left color=red]
    (3,2/3+\x) -- (2/3+\x,3) -- (1+\x,3) -- (3,1+\x) -- cycle;

  }

    
    
   
     \foreach \x in {0, 1/3, 2/3, 1, 4/3, 5/3, 2, 7/3, 8/3} {
        \draw[thin, draw = black] (\x,0) -- (\x,3);
          \draw[thin , draw = black] (0,\x) -- (3,\x);    
         \draw[thin , draw = black] (\x,0) -- (0,\x);
         \draw[thin , draw = black] (\x,3) -- (3,\x);
         \draw[thin, draw = black] (\x,0) -- (\x,3);
          \draw[thin , draw = black] (0,\x) -- (3,\x);

          }

          \foreach \x in {0,1,2,3} {
         \draw[ultra thick, draw = black] (\x,0) -- (\x,3);
          \draw[ultra thick , draw = black] (0,\x) -- (3,\x);    
         \draw[ultra thick , draw = black] (\x,0) -- (0,\x);
         \draw[ultra thick , draw = black] (\x,3) -- (3,\x);
         \draw[ultra thick, draw = black] (\x,0) -- (\x,3);
          \draw[ultra thick , draw = black] (0,\x) -- (3,\x);    
          }

   \end{scope}
  }
\end{tikzpicture}

\begin{tikzpicture}
\end{tikzpicture}

}


\newcommand\DrawTwodimGridDiagFundDomain{
\begin{tikzpicture}[fill=white,>=stealth]
  \visible<4->{
    \begin{scope}[xshift = -1cm, yshift = -1cm]

      \foreach \x in {0, 0.5,1,1.5 ,2, 2.5} {
     \node at (0.365+\x,0.365+\x) {\color{magenta}$-$};
     \node at (0.125+\x,0.125+\x) {\color{magenta}$+$};

  }

          \foreach \x in {0,.5,1,1.5,2,2.5,3} {
         \draw[ultra thick , draw = black] (\x,0) -- (0,\x);
         \draw[ultra thick , draw = black] (\x,3) -- (3,\x);
          }

        \draw [ultra thick, cyan, >=serif cm,<->] (1,1) -- (1.5,1);

     
   \end{scope}
  }
  \DrawAxisLabels
\end{tikzpicture}

\begin{tikzpicture}
\end{tikzpicture}

}

\begin{figure}[t!]
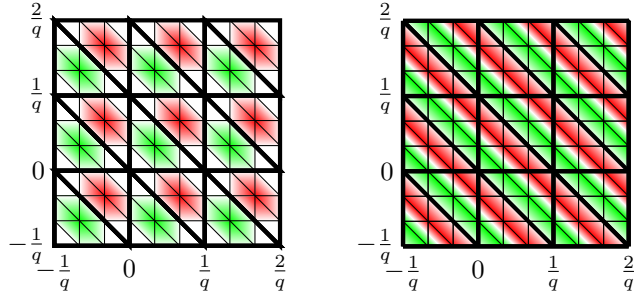

  \centering
  \begin{tabular}{ccc}
    \DrawTwodimGridPointPerturbation
    &\qquad \qquad& 
    \DrawTwodimGridDiag
  \end{tabular}
  \caption{Perturbation functions $\psi^m_{q,\point}$ (\emph{left}) and
    $\varphiD^m$ (\emph{right}) for $m=3$. 
    Colors indicate whether the value of the function is negative
    (\emph{red}), positive (\emph{green}), zero (\emph{white}).
    Two polyhedral complexes are drawn: $\P_q$ (\emph{thick lines})
    and its refinement $\P_{mq}$ (\emph{thin lines}).}
  \label{fig:2d-perturbation-functions}
\end{figure}

For $m\geq 3$, we will use the subtriangulation (refinement) $\P_{mq}$ of $\P_{q}$.  We define $\psiPoint\colon\R^2 \to \R$ that is piecewise linear over $\P_{mq}$ as follows: at  all vertices of $\P_{mq}$ that lie on the boundary of~$\FundaTriangleLower$, let $\psiPoint$ take the value $0$, and at all vertices of $\P_{mq}$ that lie in the interior of~$\FundaTriangleLower$, we assign
$\psiPoint$ to have the value $1$. Interpolate these values linearly to define $\psiPoint$ on all of $\FundaTriangleLower$. For every point $\x$ in $\FundaTriangleUpper$ define $\psiPoint(\x) = -\psiPoint\bigl(\ColVec[q]{1}{1}-\x\bigr)$. Finally, for any $\y \in \R^2$, let $\x \in [0,\frac{1}{q}]^2$ and $\t \in \frac{1}{q}\Z^2$ be vectors such that $\y = \x + \t$; define $\psiPoint(\y) = \psiPoint(\x)$. Since $\psiPoint$ vanishes on the boundary of $[0,\frac{1}{q}]^2$, this periodic extension is well-defined. The function for $m=3$ is shown in Figure~\ref{fig:2d-perturbation-functions} (left).  


The following result is quite easy to verify from the definition of $\psiPoint$. Formally, the assertions follow from (i), (iv) and (v) in Lemma~\ref{lemma:our-equivariant-psi}, in \autoref{s:reflection-groups},
 whose proof uses more general tools which, in our opinion, are of independent interest.

\begin{lemma}\label{lem:our-equivariant-psi}
For every $m\geq 3$, the function $\psiPoint\colon\R^2\to\R$ constructed above has the following properties:
  \begin{enumerate}[\rm(i)]
  \item  $\psiPoint\big|_I = 0$ on all edges and vertices $I \in \Ipointedge$.
      \item Let $i=1,2$ or $3$, and let $F\in \Delta \P_q$ be such that $p_i(F) \in \Ipoint$. Then, $F \subseteq E(\psiPoint)$.
  \item $\psiPoint$ is continuous piecewise linear over $\P_{mq}$.
  \end{enumerate}
\end{lemma}

We will also need another class of functions $\varphiD^m\colon\R^2 \to \R$
parametrized by $m \geq 3$. Let $\varphiD^m\colon \R^2 \to \R$ be the piecewise linear function over $\P_{mq}$ defined in the following way. The values on the vertices of $\P_{mq}$ are given as follows: for any $\x \in \verts(\P_{mq})$, 
$$
\varphiD^m(\x) = 
\begin{cases}
1 & \text{if } \ve1\cdot\x \equiv \frac{i}{mq} \pmod{\frac{1}{q}}  \text{ for any } 1 \leq i < \tfrac{m}{2}, i \in \Z,\\
-1 & \text{if }\ve1\cdot\x \equiv \frac{i}{mq} \pmod{\frac{1}{q}} \text{ for any } \tfrac{m}{2} < i \leq m-1, i \in \Z,\\
0 &  \text{if }\ve1\cdot\x \equiv  0 \text{ or } \frac{1}{2q} \pmod{\frac{1}{q}}.
\end{cases}
$$
The function $\varphiD^m$ is then uniquely extended to $\R^2$ continuously by interpolation on the faces of $\P_{mq}$. The function is shown for $m=3$ in Figure~\ref{fig:2d-perturbation-functions} (right). 

The next result can also be easily verified from the definition of the function $\varphiD^m$.  Formally, we again the assertions follow from (i), (iv) and (v) in Lemma~\ref{lemma:our-equivariant-varphiD} in \autoref{s:reflection-groups}.

\begin{lemma}\label{lem:our-equivariant-varphiD}
  The function $\varphiD^m\colon\R^2\to\R$ constructed above is well-defined and has the following properties:
  \begin{enumerate}[\rm(i)]
  \item  $\varphiD^m\big|_I = 0$ on all edges and vertices $I \in \Ipointdiag$.
      \item Let $i=1,2,$ or $3$ and let $F \in \Delta \P_q$ be such that $p_i(F) \in \Ipointdiag$. Then, $F \subseteq E(\varphiD^m)$.
  \item $\varphiD^m$ is continuous piecewise linear over $\P_{mq}$.
  \end{enumerate}
\end{lemma}

\subsection{Non-extremality by equivariant perturbation}
\label{sec:non-extreme-by-perturbation}
\label{sec:non-extreme-by-diag-perturbation}
In this subsection, we will prove the following lemma that shows that when $\pi$ is piecewise linear over $\P_q$, it must be affine imposing for it to be extreme.  This is done by defining specific perturbations that can be used to show $\pi$ is not extreme.  

We will derive sufficient conditions for extremality in the subsequent subsection.

\begin{lemma}
  \label{lemma:not-extreme}
  Let $\pi$ be a minimal, continuous piecewise linear function over $\P_q$ that is diagonally constrained.  If $\barStri^2 \neq \Itri$, then $\pi$ is not extreme.
\end{lemma}

In the proof, we will need $\psiPoint$ and $\varphiD^m$, as constructed in \autoref{s:eq-perturb}.  We first will analyze a case that uses $\psiPoint$.

Recall that $\Delta \pi(\x,\y) := \pi(\x) + \pi(\y) - \pi(\x + \y)$ and that when $\pi$ is piecewise linear over $\P_q$, we have that $\Delta \pi$ is piecewise linear over $\Delta \P_q$, as explained in section~\ref{sec:triangulation}.

\begin{lemma}[Perturbation only on interior of triangles]\label{lemma:not-extreme1}
Let $\pi$ be a minimal, continuous piecewise linear function over $\P_q$ with $\f \in \verts(\P_q)$ that is diagonally constrained. Suppose there exists $I^* \in \Itri\setminus ( \Stri^2 \cup \Stri^1)$.  Then $\pi$ is not extreme.  

Furthermore, for any $m \in \Z_{\geq 3}$, there exist distinct minimal valid functions $\pi^1, \pi^2$ that are continuous piecewise linear over $\P_{mq}$ such that $\pi = \tfrac{1}{2}(\pi^1 + \pi^2)$. 
\end{lemma}

\begin{proof}
Fix $m \in \Z_{\geq 3}$.  Let $R = \bigcup\{\,J \st \EqClass{J} \in \G_{I^*}\,\}$.    
Let $\psiPoint\colon \R^2 \to \R$ be the function constructed in \autoref{s:eq-perturb}.
Let $\bar \pi = \delta_R  \cdot\psiPoint$ where $\delta_R$ is the indicator function for the set $R$. 
Since $\{\0\}, \{\f\} \in \Ipoint$, by Lemma~\ref{lem:our-equivariant-psi}\,(i), we have $\psiPoint(\0) =  0$ and $\psiPoint(\f) = 0$.  Hence, $\bar \pi(\0) =  0$ and $\bar \pi(\f)= 0$.

Since $\psiPoint\big|_I = 0$ for all $I \in \Iedge$ and $R$ is a union of faces in $\Itri$, we find that $\bar \pi$ is continuous.  Since $\psiPoint$ is piecewise linear over $\P_{mq}$, $\bar \pi$ is also piecewise
linear over $\P_{mq}$. 
  Finally, notice that $\bar \pi$ is periodic modulo $\Z^2$ since $\psiPoint$ and $\delta_R$ are both periodic modulo $\Z^2$.

We will show that $E(\pi) \subseteq E(\bpi)$.  Since $I^* \in R$ and $\psiPoint \neq 0$ on $\intr(I^*)$, we have that $\bar \pi \not\equiv 0$.  Since $\bar \pi(\f) = 0$ and $\bar \pi \not\equiv 0$, by \autoref{corPerturb}, this will show that $\pi$ is not extreme.
By Lemma~\ref{lemma:covered-by-maximal-valid-triples}, we only need to consider maximal faces in the  complex $\Delta \P_q$.  Let $F \in E_{\max}(\pi, \P_q)$.  Define $\Delta \bpi(\x,\y) := \bpi(\x) + \bpi(\y) - \bpi(\x + \y)$.    
We will show that $\Delta\bpi|_F = 0$.  Note that $\bpi$ is defined over the finer complex $\P_{mq}$.  Therefore $\Delta \bpi$ is piecewise linear over $\Delta \P_{mq}$.  Since $F \in \Delta\P_q$, the function $\Delta \bpi$ is not necessarily affine over $F$.

Let $I=p_1(F)$, $J = p_2(F)$, and $K=p_3(F)$.  By Lemma~\ref{lemma:setsIJKF}, $F = F(I,J,K)$.
Since $\pi$ is diagonally constrained, we enumerate the possible cases for $I,J,K$ as listed in Lemma~\ref{lemma:cases} and show that $F = F(I,J,K) \subseteq E(\bpi)$.  Observe that we can write $\Delta \bpi|_F(\x,\y) = \bpi|_I(\x) + \bpi|_J(\y) - \bpi|_K(\x + \y)$ and that $F \subseteq E(\bpi)$ if and only if $\Delta \bpi|_F = 0$.
\begin{enumerate}[(Type 1)]
\item\label{type:all-point/diag-proof} $I,J,K \in \Ipointdiag$.   By Lemma~\ref{lem:our-equivariant-psi}\,(i), $\psiPoint = 0  = \bpi$ on the faces $I,J,K$ and thus we  have $\Delta \bpi|_F = 0$.
\item\label{type:all-tri-proof}  $I, J,K \in \Itri$.  By definition of $\Stri^2$, we have $I,J,K \in \Stri^2$.  Therefore $I\cap R, J\cap R, K \cap R \in \Ipointedge$. By Lemma~\ref{lem:our-equivariant-psi}\,(i), $\psiPoint = 0$ on $I\cap R$, $J\cap R$, $K \cap R$.  Since $\delta_R = 0$ on $I \setminus R$, $J \setminus R$, $K \setminus R$, we have $\bpi = 0$ on $I,J,K$ and  thus $\Delta \bar \pi|_F = 0$.
\item\label{type:tri,tri,point-proof} One of $I,J,K$ is in $\Ipoint$, while the other two are in $\Itri$.   Label $I,J,K$ as $I',J',K'$ where $I' \in \Ipoint$ and $J', K' \in \Itri$.  By Lemma~\ref{lem:our-equivariant-psi}\,(i), $\psiPoint = 0  = \bpi$ on $I'$.  We consider four cases.
\begin{enumerate}[{Case} i.]
\item $\EqClass{J'}, \EqClass{K'} \notin \G_{I^*}$.  Then $J' \cap R, K' \cap R \in \Ipointedge$.  By Lemma~\ref{lem:our-equivariant-psi}\,(i), $\psiPoint = 0 = \bpi$ on $J' \cap R$ and $K'\cap R$.  Furthermore, $\delta_R = 0$ on $J' \setminus R$ and $K' \setminus R$.  Hence, $\bpi = 0$ on $I',J',K'$ and hence $\Delta \pi|_F = 0$.
\item $\EqClass{J'}, \EqClass{K'} \in \G_{I^*}$.  By the relations in Lemma~\ref{lem:our-equivariant-psi}\,(ii) and the fact that $\delta_R = 1$ on $J', K'$, we have that $\Delta \bpi|_F = 0$.
\item $\EqClass{J'} \in \G_{I^*}$, $\EqClass{K'} \notin \G_{I^*}$.  We show that this case cannot happen.  Since $F \in E(\pi)$ and $I' \in \Ipoint$, we have that $\{\EqClass{J'}, \EqClass{K'}\} \in \E_\point$.  Therefore, $\EqClass{K'} \in \G_{J'}$.    Since $\EqClass{J'} \in \G_{I^*}$, we have that $\G_{I^*} = \G_{J'}$, which is a contradiction because then $\EqClass{K'} \in \G_{I^*}$.
\item $\EqClass{K'} \in \G_{I^*}$, $\EqClass{J'} \notin \G_{I^*}$.  This is similar to the previous case.
\end{enumerate}
\item\label{type:tri,tri,edge-proof} One of $I,J,K$ is in $\Idiag$, while the other two are in $\Itri$.
In this case, by definition, the two triangles are in $\Stri^1$.  Since triangles in $\Stri^1$ only intersect $R$ on lower-dimensional faces $\Ipointedge$, we have that $I\cap R, J\cap R, K\cap R \in \Ipointedge$.  By Lemma~\ref{lem:our-equivariant-psi}\,(i), $\psiPoint = 0 = \bpi$ on $I \cap R$, $J \cap R$, and $K\cap R$.  Since $\delta_R = 0$ on $I\setminus R$, $J\setminus R$ and $K \setminus R$, we have $\bpi = 0$ on $I,J,K$ and  we  have $\Delta \bar \pi|_F = 0$.
\end{enumerate}
We conclude that $E(\pi) \subseteq E(\bar \pi)$, $\bar \pi(\f) = 0$, and $\pi$ and $\bpi$ are both piecewise linear over $\P_{mq}$.  Therefore, by \autoref{corPerturb}, $\pi$ is not extreme and there exist distinct minimal functions $\pi^1, \pi^2$ that are continuously piecewise linear over $\P_{mq}$.
\end{proof}

We next use the function $\varphiD^m$, as defined in subsection~\ref{s:eq-perturb}, as the basis for a perturbation function $\bpi$.  

As in Lemma~\ref{lemma:not-extreme1}, we will allow the perturbation $\bar \pi$ to only apply to a subset of the triangles, this time corresponding to a connected component in the graph $\bar \G$.  Since $\varphiD^m$ is non-zero on the vertical and horizontal faces $\Iverthor$, we must be careful about geometrically adjacent triangles.


To handle the geometrically adjacent triangles easier, we consider the case where $\Itri =  \Stri^2 \cup \Stri^1$. 

\begin{figure}
\begin{tikzpicture}[scale=0.4]
\draw[top color=pink,bottom color=pink] (0,0) -- (5,0) -- (5,5) -- (0,5)-- cycle; 
\draw[top color=white,bottom color=white] (3,0) -- (0,3) -- (0,4) -- (4,0) -- cycle; 
\draw[thin, blue] (3.25,0) -- (0,3.25) -- (0,3.75) -- (3.75,0) -- cycle; 
\draw[thin, blue] (3.5,0) -- (0,3.5) -- (0,3.75) -- (3.75,0) -- cycle; 
\draw[top color=white,bottom color=white] (3,5) -- (5,3) -- (5,4) -- (4,5) -- cycle; 
\draw[thin, blue] (3.25,5) -- (5,3.25) -- (5,3.75) -- (3.75,5) -- cycle; 
\draw[thin, blue] (3.5,5) -- (5,3.5) -- (5,3.5) -- (3.5,5) -- cycle; 
\diagGrid
\end{tikzpicture}
\caption{A case where $\Itri = \Stri^2 \cup \Stri^1=
  \barStri^2 \cup \barStri^1$.  Therefore, on every triangle, $\pi$ is either
  affine imposing (\emph{shaded triangles}), or only diagonally affine
  imposing (\emph{striped triangles}). 
  Observation~\ref{obs:adjacent-apply-lemma} shows that a shaded triangle that
  is geometrically adjacent to a striped triangle along a vertical or
  horizontal face in $\Iverthor$ forces the striped triangle to become shaded.
  Therefore, no striped triangle can be geometrically adjacent to a
  shaded triangle along a vertical or horizontal face in $\Iverthor$.  In this
  example, every striped triangle is connected in the graph $\bar \G$ by a
  path with edges in $\E_{\ver \hor}$.  Therefore, all of these triangles form
  a connected component in the graph.  We can choose any one of these
  triangles as $I^*$ in Lemma~\ref{lemma:not-extreme2} to perturb on this
  connected component.  } 

\label{fig:diagonal-stripe}
\end{figure}
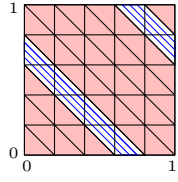

\begin{obs}
\label{obs:connections}
Suppose $\Itri = \Stri^2 \cup \Stri^1$ and let $I^* \in \barStri^1$.  Let $J,K \in \Itri$ such that $\EqClass{J} \in \bar\G_{I^*}$ and $J \cap K \in \P_{q,\ver\hor}$.  Then $\EqClass{K} \in \bar\G_{I^*}$ as well.  This is because  $\Itri = \Stri^2 \cup \Stri^1$ and therefore $\{\EqClass{J},\EqClass{K}\} \in \E_{\ver \hor} \subseteq \bar\E$;  see Figure~\ref{fig:diagonal-stripe}. 
\end{obs}

\begin{lemma}[Diagonal perturbation touching vertical and horizontal boundaries of triangles]\label{lemma:not-extreme2}
Suppose $\pi$ is continuous piecewise linear over $\P_q$ with $\f \in \verts(\P_q)$ and is diagonally constrained. 
Suppose further that  $\Itri = \Stri^2 \cup \Stri^1$ and there exists $I^* \in \barStri^1$.  Then $\pi$ is not extreme.

Furthermore, for any $m \in \Z_{\geq 3}$, there exist distinct minimal valid functions $\pi^1, \pi^2$ that are continuous piecewise linear over $\P_{mq}$ such that $\pi = \tfrac{1}{2}(\pi^1 + \pi^2)$. 
\end{lemma}
\begin{proof}
 Let $R = \bigcup\{\,J\in \Itri \st \EqClass{J} \in \bar\G_{I^*}\,\}$.
 Note that $\bar\G_{I^*} \subseteq \barStri^1$ and recall that $\barStri^1 \cap \barStri^2 = \emptyset$.  Furthermore, $\barStri^1 \cap \Stri^2 = \emptyset$.   Let $\varphiD^m\colon \R^2 \to \R$ be the function constructed in~\autoref{s:eq-perturb}.

Let $\bar \pi = \delta_{R}(\x)  \cdot\varphiD^m(\x)$.  First recognize that $\bar \pi$ is a continuous function. To see this, note that $\varphiD^m(\x)$ is continuous and $\delta_{R}$ is continuous on $R$ and $\R^2 \setminus R$.  By Observation~\ref{obs:connections}, it follows that $\partial R \subseteq \bigcup \{\, I \st I \in \P_{q,\diag}\, \}$.  By \autoref{lem:our-equivariant-varphiD}\,(i), $\varphiD^m$ vanishes on $\partial R$, that is $\varphiD^m = 0$ on $\partial R$.  These together imply that $\bpi$ is continuous.

Since $\varphiD^m$ is piecewise linear over $\P_{mq}$, $\bar \pi$ is also piecewise
linear over $\P_{mq}$.  Also, since $\varphiD^m\big|_I = 0$ for all $I \in \Idiag$, we find $\bar \pi$ is also continuous. 
Finally, notice that $\bar \pi$ is periodic modulo $\Z^2$ since $\varphiD^m$ and $\delta_R$ are both periodic modulo $\Z^2$.

We will show that $E(\pi) \subseteq E(\bpi)$.  Since $I^* \in R$ and $\varphiD^m \not\equiv 0$ on $\intr(I^*)$, we have $\bar \pi \not\equiv 0$.  Since $\bar \pi(\f) = 0$ and $\bar \pi \not\equiv 0$, by \autoref{corPerturb}, this will show that $\pi$ is not extreme.
By Lemma~\ref{lemma:covered-by-maximal-valid-triples}, we only need to consider maximal faces in the complex $\Delta \P_q$.  Let $F \in E_{\max}(\pi, \P_q)$.  

Define $\Delta \bpi(\x,\y) := \bpi(\x) + \bpi(\y) - \bpi(\x + \y)$.    
We will show that $\Delta\bpi|_F = 0$.  Note that $\bpi$ is defined over the finer complex $\P_{mq}$.  Therefore $\Delta \bpi$ is piecewise linear over $\Delta \P_{mq}$.  Since $F \in \Delta \P_q$, the function $\Delta \bpi$ is not necessarily affine over $F$.

Let $I=p_1(F)$, $J = p_2(F)$, and $K=p_3(F)$.  By Lemma~\ref{lemma:setsIJKF}, $F = F(I,J,K)$.
Since $\pi$ is diagonally constrained, we enumerate the possible cases for $I,J,K$ as listed in Lemma~\ref{lemma:cases} and show that $F = F(I,J,K) \subseteq E(\bpi)$.  Observe that we can write $\Delta \bpi|_F(\x,\y) = \bpi|_I(\x) + \bpi|_J(\y) - \bpi|_K(\x + \y)$ and that $F \subseteq E(\bpi)$ if and only if $\Delta \bpi|_F = 0$.

\begin{enumerate}[(Types 3 $\&$ 4)]
\item[(Type 1)]\label{type:all-point/diag-proof} $I,J,K \in \Ipointdiag$.   By Lemma~\ref{lem:our-equivariant-varphiD}\,(i), $\varphiD^m = 0 = \bar\pi$ on all faces $I,J,K$ and thus we  have $\Delta \bar \pi|_F = 0$.

\item[(Type 2)]\label{type:all-tri-proof}  $I, J,K \in \Itri$.  By definition of $\barStri^2$, we have $I,J,K \in \barStri^2$.  By Observation~\ref{obs:connections}, we must have $I\cap R,J\cap R,K\cap R \in \Ipointdiag$ and hence $\varphiD^m = 0$ on $I\cap R,J\cap R,K\cap R$ by \autoref{lem:our-equivariant-varphiD}\,(i).  Therefore, $\bpi = 0$ on $I,J,K$ and  we  have $\Delta \bar \pi|_F = 0$.

\item[(Types 3 $\&$ 4)]\label{type:tri,tri,point-proof} One of $I,J,K$ is in $\Ipointdiag$, while the other two are in $\Itri$.  Label $I,J,K$ as $I',J',K'$ where $I' \in \Ipointdiag$ and $J', K' \in \Itri$.  By Lemma~\ref{lem:our-equivariant-varphiD}\,(i), $\varphiD^m = 0$ on $I'$.  We consider four cases.
\begin{enumerate}[{Case} i.]
\item $\EqClass{J'}, \EqClass{K'} \notin \bar\G_{I^*}$.   By Observation~\ref{obs:connections}, we must have $J' \cap R, K' \cap R \in \Ipointdiag$.  Therefore, by Lemma~\ref{lem:our-equivariant-varphiD}\,(i), $\varphiD^m = 0$ on $J' \cap R, K' \cap R$, while $\delta_R = 0$ on $J' \setminus R$, and $K' \setminus R$.  Therefore $\bpi = 0$ on $I,J,K$ and hence $\Delta \pi_F = 0$.
\item $\EqClass{J'}, \EqClass{K'} \in \bar\G_{I^*}$.  By Lemma~\ref{lem:our-equivariant-varphiD}\,(ii) and the fact that $\delta_R = 1$ on $J', K'$,  we have that $\Delta \bpi|_F = 0$.
\item  $\EqClass{J'} \in \bar\G_{I^*}$, $\EqClass{K'} \notin \bar\G_{I^*}$.  We show that this case cannot happen.  Since $F \subseteq E(\pi)$ and $I' \in \Ipointdiag$, we have that $\{\EqClass{J'}, \EqClass{K'}\} \in \E\subseteq \bar \E$.  Therefore, $\EqClass{K'} \in \bar\G_{J'}$.    Since $\EqClass{J'} \in \bar\G_{I^*}$, we have that $\bar\G_{I^*} = \bar\G_{J'}$ which is a contradiction because then $\EqClass{K'} \in \bar\G_{I^*}$.
\item $\EqClass{K'} \in \bar\G_{I^*}$, $\EqClass{J'} \notin \bar\G_{I^*}$.  This is similar to the previous case.
\end{enumerate}

\end{enumerate}
We conclude that $E(\pi) \subseteq E(\bar \pi)$, $\bar \pi(\f) = 0$, and $\pi$ and $\bar \pi$ are both continuous piecewise linear over $\P_{mq}$.  Therefore, by \autoref{corPerturb}, $\pi$ is not extreme and there exist distinct minimal functions $\pi^1, \pi^2$ that are continuously piecewise linear over $\P_{mq}$.
\end{proof}

\begin{proof}[Proof of Lemma~\ref{lemma:not-extreme}]
This follows directly from Lemmas \ref{lemma:not-extreme1} and \ref{lemma:not-extreme2}.
\end{proof}

The specific form of our perturbations as continuous piecewise linear
functions over $\P_{mq}$ implies the following corollary.

\begin{corollary}
\label{corollary:AIG4}
Fix $m \in \Z_{\geq 3}$.  Suppose $\pi$ is a continuous piecewise linear function over $\P_q$ and is diagonally constrained.  If $\pi$ is not affine imposing over $\Itri$, then there exist distinct minimal
$\pi^1, \pi^2$ that are continuous piecewise linear over $\P_{mq}$ such that $\pi = \frac12(\pi^1 + \pi^2)$.  
\end{corollary}

\subsection{Extremality and non-extremality by linear algebra}
\label{section:system}

In this section we suppose $\pi$ is a minimal  continuous piecewise linear function over $\P_{q}$ that is affine imposing in  $\Itri$.  Therefore, by Lemma~\ref{lem:tightness} and Definition~\ref{def:affine-imposing}, $\pi^1$ and $\pi^2$ must also be minimal continuous piecewise linear functions over $\P_{q}$.   Recall from  Lemma~\ref{lemma:tight-implies-tight} that $E(\pi) \subseteq E(\pi^1), E(\pi^2)$.

We now set up a system of linear equations that $\pi$ satisfies and that
$\pi_1$ and $\pi_2$ must also satisfy.
Let $\varphi\colon \frac{1}{q}\Z^2 \to \R$ be a periodic function modulo $\Z^2$.  
Suppose $\varphi$ satisfies the following system of linear equations:
\begin{equation}\tag{$\mathrm E_q(\pi)$}
\begin{cases}
\varphi(\0) = 0,\\
\varphi(\f) = 1,\\
\varphi(\u) + \varphi(\v)  =  \varphi(\u + \v) 
&\text{ for all $\u,\v \in \frac{1}{q} \Z^2$ with } 
\pi(\u) + \pi(\v)  =  \pi(\u + \v).
\end{cases}
\end{equation}

 Since $\pi$ exists and satisfies ($\mathrm E_q(\pi)$), we know that the
 system has a solution.  Since $\varphi$ and $\pi$ are periodic, we can identify
 variables $\varphi(\x)$ and $\varphi(\x + \t)$ for $\x\in \frac1q\Z^2$ and $\t \in
 \Z^2$, and thus the system can be
 represented with finitely many variables and finitely many equations.

\begin{theorem}
\label{theorem:systemNotUnique}
\citedinsurveyas{Theorem 5.16}
Let $\pi\colon \R^2 \to \R$ be a continuous piecewise linear valid function
over $\P_q$.
\begin{enumerate}[\rm(i)]
\item If the system $(\mathrm E_q(\pi))$ does not have a
  unique solution, then $\pi$ is not extreme. 
\item Suppose $\pi$ is minimal and affine imposing in $\Itri$.
  Then $\pi$ is extreme if and only if the system of equations
  $(\mathrm E_q(\pi))$ has a unique solution.  
\end{enumerate}
\end{theorem}

The proof is similar to the proof of~\cite[Theorem 4.11]{basu-hildebrand-koeppe:equivariant}.  

\begin{proof}
  \emph{Part (i).}
Suppose ($\mathrm E_q(\pi)$) does not have a unique solution. Let $\bar \varphi\colon \tfrac{1}{q} \Z^2\to\R$ be a non-trivial element in the kernel of the system above.  Then for any $\epsilon$, $\pi\big|_{\tfrac{1}{q} \Z^2} + \epsilon \bar\varphi$ also satisfies the system of equations. 
Let $\bar\pi\colon \R^2 \to \R$ be the continuous piecewise linear extension of $\bar \varphi$ over $\P_q$.  Therefore $\bar \pi(\f) = 0$ and $\bar \pi \not\equiv 0$.  
Let $\u,\v \in \frac{1}{q}\Z^2$.
If $\Delta\pi(\u,\v) = 0$, then $\Delta\varphi(\u,\v) = 0$, as implied by the system of equations.  Since $\verts(\Delta \P_q) \subseteq \frac{1}{q}\Z^2$, this shows that for any $\x,\y \in \R^2$, $\Delta \pi(\x,\y) = 0$ implies that $\Delta \bar \pi(\x,\y) = 0$.   Therefore $E(\pi) \subseteq E(\bar \pi)$.
  Therefore, by \autoref{corPerturb}, $\pi$ is not extreme.
\smallbreak

\emph{Part (ii).}
Suppose there exist distinct, valid functions $\pi^1, \pi^2$ such that $\pi =
\tfrac12(\pi^1 + \pi^2)$.  Since $\pi$ is minimal and affine
imposing in $\Itri$, $\pi^1,\pi^2$ are minimal continuous piecewise linear functions over $\P_q$.  Furthermore, $\pi\big|_{\tfrac{1}{q} \Z^2}$ and, also $\pi^1\big|_{\tfrac{1}{q} \Z^2}, \pi^2\big|_{\tfrac{1}{q} \Z^2}$ satisfy the system
of equations ($\mathrm E_q(\pi)$).  If this system has a unique solution,
then $\pi = \pi^1  = \pi^2$, which is a contradiction since $\pi^1, \pi^2$
were assumed distinct.  Therefore $\pi$ is extreme. 

On the other hand, if the system ($\mathrm E_q(\pi)$) does not have a unique solution, then by Part (i), $\pi$\ is not extreme.
\end{proof}

\subsection{Connection to a finite group problem}
\label{sec:connection-to-finite-group}
\begin{theorem}\label{thm:1/4q test}
Fix $m \in \Z_{\geq 3}$.  Let $\pi$ be a minimal continuous piecewise linear function over $\P_q$ that is diagonally constrained.  
Then $\pi$ is extreme if and only if the system of
equations $(\mathrm E_{mq}(\pi))$ with $\tfrac{1}{mq}\Z^2$
has a unique solution.
\end{theorem}
\begin{proof}
  Since $\pi$ is  piecewise linear over $\P_q$, it is also  piecewise linear over $\P_{mq}$.  The forward direction is  the contrapositive of Theorem~\ref{theorem:systemNotUnique}\,(i), applied when we view $\pi$ as piecewise linear over $\P_{mq}$. 
  For the reverse direction, observe that 
  if the system of
  equations ($\mathrm E_{mq}(\pi)$) has a unique
  solution, then there cannot exist distinct minimal $\pi^1, \pi^2$ that are
  continuous piecewise linear over~$\P_{mq}$ such that $\pi = \tfrac12(\pi^1 + \pi^2)$.  By the contrapositive of Corollary \ref{corollary:AIG4},
  $\pi$ is affine imposing in $\Itri$.  
  Then $\pi$ is also affine imposing on $\Itri[mq]$ since it is a finer set.  
  By Theorem~\ref{theorem:systemNotUnique}\,(ii), since $\pi$ is affine imposing in $\Itri[mq]$ and the system of
  equations ($\mathrm E_{mq}(\pi)$) on $\P_{mq}$ has a unique solution, $\pi$ is extreme.  
\end{proof}

Theorem \ref{thm:main} and Theorem \ref{thm:1/4q} are direct consequences of Theorem \ref{thm:1/4q test}.

\section{Conclusions} 

In the present paper, we have focused on the diagonally constrained case. 
In a similar way, \emph{horizontally constrained} or \emph{vertically
  constrained} minimal functions can be defined, and the theory developed in
this paper can be easily adapted to these cases. 

The general case in which the restriction to diagonally constrained functions
is removed and thus all degenerations of maximal additive faces are allowed
(see \autoref{fig:constrained-functions-poset}, top) 
requires the solutions of more general functional equations and leads to the
construction of more complicated perturbation functions. This is analogous to the history of the one-dimensional ($k=1$) functions. Prior to the work in~\cite{basu-hildebrand-koeppe:equivariant}, extremality proofs were known only for full-dimensionally constrained piecewise linear functions (see Figure~\ref{fig:constrained-functions-poset}, left). A new idea was needed to handle the more general case for $k=1$.
Similar development for $k=2$ is deferred to the forthcoming paper
\cite{basu-hildebrand-koeppe:equivariant-general-2dim}.


\appendix

\section{Reflection groups and equivariant perturbations}
\label{s:reflection-groups}

We provide a general framework to motivate the definition of the functions $\psiPoint$ and $\varphiD^m$ from \autoref{s:eq-perturb}. We describe the construction at a more abstract level with the hope that it could be a useful tool to analyze infinite group problems in higher dimensions.

We follow the direction of~\cite{basu-hildebrand-koeppe:equivariant} where the relevant
arithmetics of the one-dimensional problem is captured by studying sets of
additivity relations of the form $\pi( t^i) + \pi(y) = \pi( t^i + y)$ and 
$\pi(x) + \pi( r^i-x) = \pi( r^i)$, where the points $ t^i$ and $ r^i$ are breakpoints of a one-dimensional minimal valid function~$\pi$. 
This is an important departure from the previous
  literature, which only uses additivity relations over non-degenerate
  intervals. The arithmetic nature of the problem comes into focus when one
realizes that isolated additivity relations over single points are also important for studying
extremality. These isolated additivity relations
give rise to a subgroup of the group $\Aff(\R^k)$ of invertible affine linear
transformations of~$\R^k$ as follows. 

\subsection{Reflection groups and their fundamental domains} 
For a point $\ve r \in \R^k,$ define the \emph{reflection} $\rho_{\ve r}\colon \R^k\to \R^k$, $\ve x
  \mapsto \ve r-\ve x$.  For a vector $\ve t \in \R^k,$ define the \emph{translation}
  $\tau_{\ve t} \colon \R^k\to \R$, $\ve x \mapsto \ve x + \ve t$.  We consider the reflections $\rho_{\ve r}$ and translations $\tau_{\ve t}$ as elements of the group $\Aff(\R^k)$.

Given a set $R$ of points in $\R^k $ and a set $T$ of vectors in $\R^k$, we define the \emph{reflection group}
$\Gamma= \Gamma(R,T)=\langle\, \rho_{\ve r}, \tau_{\ve t} \st \ve r\in R,\, \ve t\in T\,\rangle$. 
A \emph{group character} of~$\Gamma$ is a group homomorphism $\chi\colon \Gamma\to\C^\times$. The \emph{orbit} of a point~$\x\in\R^k$ under the group~$\Gamma \subseteq \Aff(\R^k)$ is the set
  $ \Gamma(\x) = \{\, \gamma(\x) \st \gamma\in\Gamma \,\}.$ We extend this notation to subsets of $\R^k$: for a subset $X \subseteq \R^k$, $\Gamma(X) = \bigcup_{\x \in X}\Gamma(\x)$.

In the following, we assume that $R \neq \emptyset$, i.e., at least one of the generators is a reflection.
The structure of the group $\Gamma$ is easy to describe completely.
The following lemma, which appeared in \cite{basu-hildebrand-koeppe:equivariant} for $k=1$, summarizes the structure of this group and generalizes easily from \cite{basu-hildebrand-koeppe:equivariant}.

\begin{lemma}\label{lemma:structure}
  Let $\rr_1 \in R$.  Then the group $\Gamma=\Gamma(R,T)= \langle\, \rho_{\rr},
\tau_{\t}\st \rr \in R, \t \in T\,\rangle$ is the
  semidirect product $\Gamma^+ \rtimes \langle \rho_{\ve r_1} \rangle$,  
  where the (normal) subgroup of translations is of index~$2$ in~$\Gamma$ and
  can be written as \begin{equation}\label{eq:gamma+}\Gamma^+ = \{\, \tau_{\ve t} \st \ve t \in 
  Y\,\},\end{equation} where $Y$ is the additive subgroup of~$\R^k$ 
  \begin{equation}\label{eq:lambda}Y = {\langle\, \rr - \rr_1, \t \st \rr \in R, \t \in T \,\rangle}_\Z \subseteq \R^k.\end{equation}%
  There is a unique group character~$\chi\colon\Gamma\to\{\pm1\}\subset\C^\times$ with
  $\chi(\rho) = -1$ for every reflection~$\rho \in\Gamma$ and $\chi(\tau) =
  +1$ for every translation $\tau\in\Gamma$.
\end{lemma}%

\begin{definition}
  A function $\psi\colon\R^k\to\R$ is called \emph{$\Gamma$-equivariant} if 
  it satisfies the \emph{equivariance formula} 
  \begin{equation}\label{eq:equivariance}
    \psi(\gamma(\x)) = \chi(\gamma) \psi(\x)
    \quad\text{for $\x\in \R^k$ and $\gamma \in \Gamma$}.
  \end{equation}  
\end{definition}

We will use formula~\eqref{eq:equivariance} to give an alternative derivation of the functions $\psiPoint$ and $\varphiD^m$ defined in \autoref{s:eq-perturb}. These functions provide the perturbation functions when~\autoref{corPerturb} is invoked in subsections~\ref{sec:non-extreme-by-perturbation} and \ref{sec:non-extreme-by-diag-perturbation}.

\begin{obs}\label{obs:RcapT} 
Let $\Gamma = \Gamma(R,T)$ be a reflection group with $R\cap T \neq \emptyset$ and let $\psi$ be any $\Gamma$-equivariant function. Then, $\rho_\0 \in \Gamma$ and $\psi(\0) = 0$.
\end{obs} 
\begin{proof} Let $\rr \in R\cap T$; then $\rho_\0 = \rho_\rr \circ \tau_\rr$. Also, we have  $\psi(\0) = \psi(\rho_\0 (\0)) = \chi(\rho_\0) \psi(\0) = -\psi(\0)$; hence, $\psi(\0) = 0$. \end{proof}

It follows from Observation~\ref{obs:RcapT} that when $R \cap T \neq \emptyset$ and $\psi$ is $\Gamma$-equivariant, we have $\psi \equiv 0$ on all of $\Gamma(\0)$. If we restrict ourselves to continuous $\Gamma$-equivariant functions and $Y$ defined in~\eqref{eq:lambda} is dense in $\R^k$, then $\psi \equiv 0$ is the unique $\Gamma$-equivariant function.  On the other hand, when $Y$ has inherent discreteness properties, which we make precise in the following discussion, we can construct many non-trivial continuous $\Gamma$-equivariant functions. To do so, we only need to construct a function on a subset of $\R^k$.

\begin{definition}
 A \emph{fundamental domain} of a reflection group~$\Gamma$ is a subset of~$\R^k$ that is a system of representatives of the orbits. 
\end{definition}

Given a reflection group $\Gamma$ for $k=1$, if the group $Y$ from \eqref{eq:lambda} in Lemma~\ref{lemma:structure} is discrete, a fundamental domain of $\Gamma$ can be chosen as a certain closed interval.  In higher dimensions, when $Y$ is discrete, the fundamental domain is no longer a closed set. Even so, it is easy to describe the closure of a fundamental domain.  This is made concrete in the following discussion and Lemma~\ref{lem:fundamental-domain}.

A well known fact is that for any discrete subgroup $\Lambda$ of $\R^k$ there exists a finite set of vectors $\t^1, \dots, \t^\ell \in \R^k$ such that $\Lambda = \langle \t^1, \dots, \t^\ell \rangle_\Z$. These vectors are called the {\em basis} of $\Lambda$. We say that $\Lambda$ is a {\em lattice} of the linear subspace $\langle \t^1, \dots, \t^\ell \rangle_\R$. The set $V_\Lambda = \{\, \sum_{i=1}^\ell \lambda_i \t^i \st 0 \leq \lambda_i \leq 1\,\}$ is called the {\em closed fundamental parallelepiped} of $\Lambda$ with respect to the basis $\t^1, \dots, \t^\ell$. Define $\t = \sum_{i=1}^\ell \t^i$
%
%
and set $M :=\max\{\,\ve \t\cdot \x \st \x \in V\,\} = \t\cdot \t$. Define $V^+_\Lambda = \{\,\x \in V \st \ve \t\cdot \x \leq \tfrac{M}{2}\,\}$ and $V^-_\Lambda = \{\,\x \in V \st \ve \t\cdot \x \geq \tfrac{M}{2}\,\}$. (These definitions are of course with respect to the particular basis $\{\t^1, \ldots, \t^\ell\}$; the basis will usually be fixed in a particular context). 

A {\em mixed-lattice} is a subgroup $Y \subseteq \R^k$ such that $Y = \Lambda + L$ where $\Lambda$ is a lattice of a linear subspace $L'$ of $\R^k$, $L$ is a linear subspace of $\R^k$, and $L'$ and $L$ are complementary subspaces, i.e., $\R^k = L' + L$ and $L \cap L' = \{\0\}$. 


\begin{lemma}
\label{lem:fundamental-domain}
Let $\Gamma=\Gamma(R,T)$ be a reflection group with $\emptyset \subsetneq R\subseteq T$ such that the corresponding $Y$ from \eqref{eq:lambda} is a mixed-lattice, i.e., $Y = \Lambda + L$ and let $\t^1, \ldots, \t^\ell$ be a basis of $\Lambda$. Let $L' = \langle\, \t^1, \ldots, \t^\ell \, \rangle_\R$. Let $V^+_{\Lambda}$ be defined with respect to this basis. Then there exists a fundamental domain $\tilde V$ for $\Gamma$ such that $\intr_{L'}(V^+_{\Lambda}) \subseteq \tilde V \subseteq V^+_{\Lambda}$. 
\end{lemma}
\begin{proof}
We first show that $V^+_{\Lambda}$ contains a representative for every point $\x$ in $\R^k$. Let $\x = \sum_{i=1}^\ell \lambda_i \t^i + \p$ for some $0 \leq \lambda_i \leq 1$, $\p \in Y$. We will show that $\gamma(\x) \in V^+_{\Lambda}$ for some $\gamma \in \Gamma$. Let $\x' = \sum_{i=1}^\ell \lambda_i \t^i$ and let $\t = \sum_{i=1}^k \t^i$. If $\x' \in V^+_{\Lambda}$, then we are done by taking $\gamma = \tau_{-\p}$. Otherwise, $\ve \t\cdot \x' > \tfrac{M}{2}$. Consider $\tau_\t \circ \rho_\0 ( \x') = \t -\x' = \sum_{i=1}^\ell (1-\lambda_i) \t^i$. By \autoref{obs:RcapT}, $\rho_\0 \in \Gamma$  and so $\gamma = \tau_\t \circ \rho_\0 \circ \tau_{-\p} \in \Gamma$. Further, $\ve \t\cdot (\t - \x') = M - \t\cdot \x'  < \tfrac{M}{2}$, and hence $\gamma(\x) = \t-\x' \in V^+_{\Lambda}$. Hence, $V^+_{\Lambda}$ contains a representative for every point in $\R^k$.

Next we show that every point $\x \in \intr_{L'}(V^+_{\Lambda})$ is a unique representative in $V^+_{\Lambda}$ because for any non-trivial $\tau_\t \in \Gamma^+$, $\tau_\t (\x) \notin V_{\Lambda}$, and for any $\rr \in R\subseteq T$, $\rho_\rr (\x) = \tau_\rr\circ\rho_\0(\x)$ lies in $\Gamma^+(\intr_{L'}(V^-_{\Lambda}))$, which does not intersect $V^+_{\Lambda}$ (recall that $\Gamma^+$ is the subgroup defined in \eqref{eq:gamma+} for $\Gamma$).  
\end{proof}

\begin{figure}[t!]
  \centering
  \begin{tabular}{ccc}
     \DrawTwodimGrid &\qquad \qquad& \DrawTwodimGridDiagFundDomain
    \end{tabular}
  \caption{Reflection groups $\Gamma_{q, \point}$ and $\Gamma_{q, \diag}$ 
    and the closures of their fundamental domains (\emph{blue}) for $m=3$. 
    \emph{Left}, the case $\Gamma_{q, \point}$.  Translating the closure of
    the fundamental domain, $V_{q, \point}^+ = \FundaTriangleLower$ (\emph{blue triangle}), by $\tau_\t$ 
    for $\t \in \Lambda_{q, \point} = \tfrac{1}{q}\Z^2$ gives the triangles
    labeled with $+$.  Reflections by $\rho_\rr$ for $\rr\in\frac1q\Z^2$ take
    these triangles to the triangles labeled with $-$. 
    \emph{Right}, the case $\Gamma_{q, \diag}$.  
    Translating the closure of the fundamental domain, $V^+_{q,\diag}
    $ (\emph{blue line segment}),
    by $\tau_\t$ for $\t \in L_\diag$ gives the diagonal strip, labeled~$+$, containing
    the fundamental domain.  Further translations by $\tau_\t$ for $\t \in
    \Lambda_{q, \diag}$ give the remaining diagonal strips labeled~$+$.
    The reflections $\rho_\rr$ in $\Gamma_{q, \diag}$ then take these strips to
    the strips labeled~$-$.}
  \label{fig:2d-perturbation-funda-domains}
\end{figure}

The following lemma explains how to construct $\Gamma$-equivariant functions using the fundamental domain. 

\begin{lemma}[Construction of $\Gamma$-equivariant functions]
  \label{lemma:phi}
  Let $\Gamma=\Gamma(R,T)$ be a reflection group with $\emptyset \subsetneq R\subseteq T$ such that the corresponding $Y$ from \eqref{eq:lambda} is a mixed-lattice, i.e., $Y = \Lambda + L$ and let $\t^1, \ldots, \t^\ell$ be a basis of $\Lambda$. Let $L' = \langle\, \t^1, \ldots, \t^\ell \, \rangle_\R$.  Let $V^+_{\Lambda}$ be defined with respect to this basis. Let $\psi\colon V^+_{\Lambda}\to \R$ be any function such that
  $\psi\big|_{\partial_{L'}(V^+_{\Lambda})} = 0$,
  where $\partial_{L'}(V^+_{\Lambda})$ denotes the boundary of~$V^+_{\Lambda}$ with respect to the linear subspace $L'$.
  Then the \emph{equivariance formula}~\eqref{eq:equivariance}  
  gives a well-defined extension of $\psi$ to all of~$\R^k$.
\end{lemma}

Figures \ref{fig:2d-perturbation-functions} and ~\ref{fig:2d-perturbation-funda-domains} illustrate this construction.

\begin{proof}
By Lemma~\ref{lem:fundamental-domain}, $V^+_{\Lambda}$ contains a fundamental domain.  Since $\intr_{L'}(V^+_{\Lambda})$ has unique representatives for the orbits of $\Gamma$ and $\psi = 0$ on the boundary $\partial_{L'}(V^+_{\Lambda})$, the extension is well-defined. 
%
%
\end{proof}

\subsection{Deriving the perturbation functions $\psiPoint, \varphiD^m$ using equivariance formulas} 

In \cite{basu-hildebrand-koeppe:equivariant}, the authors use $\Gamma = \langle\,\rho_g, \tau_g \st g \in \frac{1}{q}\Z\,\rangle$,  where $Y = \Lambda =  \tfrac{1}{q}\Z$.  Using the lattice basis $\{ t^1 = \tfrac{1}{q}\}$, we obtain the fundamental parallelepiped $V_\Lambda = [0, \tfrac{1}{q}]$ and hence $V^+_\Lambda = [0, \tfrac{1}{2q}]$.  In this one-dimensional case, $V^+_\Lambda$ is actually a fundamental domain for $\Gamma$.

We proceed similarly with two different reflection groups in dimension two.  We first consider the reflection group 
 $\Gamma_{q, \point} = \langle \,\rho_\g, \tau_\g \st \g \in \frac{1}{q}\Z^2\,\rangle$
generated by reflections and translations corresponding to all possible
vertices of $\P_q$; see \autoref{fig:2d-perturbation-funda-domains} (left).  
The corresponding lattice $Y_{q, \point} = \Lambda_{q, \point} = \tfrac{1}{q}\Z^2$.  Using the  lattice basis $\bigl\{\t^1 = \ColVec[q]{1}{0},\, \t^2 = \ColVec[q]{0}{1}\bigr\}$, we obtain the fundamental parallelepiped $V_{q, \point} = [0,\tfrac{1}{q}]^2$ from which we obtain $V_{q, \point}^+ = \FundaTriangleLower = \frac1q \conv(\{ \ColVec{0}{0}
, \ColVec{1}{0}
, \ColVec{0}{1}
 \})$.   
 We make this particular choice of fundamental domain in part because 
 $V_{q, \point}^+ \in \Itri$ and $\Gamma_{q,\point}(V_{q, \point}^+) \subseteq \Itri$. (Note that we have simplified the notation, e.g., $V_{\Lambda_{q, \point}}$ is now denoted by $V_{q, \point}$.)
 
For any  $m \in \Z_{\geq 3}$, we may now interpret the function $\psiPoint \colon \R^2 \to \R$ defined in~\autoref{s:eq-perturb} in the following way: at  all vertices
of $\P_{mq}$ that lie on the boundary of~$\FundaTriangleLower$, let $\psiPoint$ take the value $0$, and at all vertices of $\P_{mq}$ that lie on the interior of  of~$\FundaTriangleLower$, 
 we assign
$\psiPoint$ to have the value $1$.  
Interpolate these values to define $\psiPoint$ on $\FundaTriangleLower$. 
By Lemma~\ref{lemma:phi}, the extension of $\psiPoint$ to $\R^2$ via the
equivariance formula~\eqref{eq:equivariance} is well-defined. This is an
alternative description for the function $\psiPoint$ defined in
\autoref{s:eq-perturb}; refer back to \autoref{fig:2d-perturbation-functions}
(left) for an illustration. 
%

One possible choice of a fundamental domain for $\Gamma_{q, \point}$ is 
$$\tilde V_{q, \point} = \intr(\FundaTriangleLower) \cup \bigl[\ColVec{0}{0},
\ColVec[2q]{0}{1}\bigr] \cup \bigl[\ColVec{0}{0}, \ColVec[2q]{1}{0}\bigr] \cup
\bigl[\ColVec[2q]{1}{1}, \ColVec[2q]{1}{0}\bigr),$$ where $[\x, \y]$ and
$[\x,\y)$ denote the closed and half open line segments, respectively, between
$\x$ and $\y$. For our construction, only its closure, $V^+_{q,\point} =
\FundaTriangleLower$, matters. 

\begin{lemma}\label{lemma:our-equivariant-psi}
  The function $\psiPoint\colon\R^2\to\R$ has the following properties:
  \begin{enumerate}[\rm(i)]
  \item  $\psiPoint\big|_I = 0$ on all edges and vertices $I \in \Ipointedge$.
  \item $\psiPoint(\x)= - \psiPoint(\rho_{\g}(\x)) = - \psiPoint(\g - \x)$ for all $\x \in
    \R^2$ and $\g \in \frac{1}{q}\Z^2$.
  \item $\psiPoint(\x) = \psiPoint(\tau_\g(\x)) = \psiPoint(\g + \x)$ for all $\x \in
    \R^2$ and $\g \in \frac{1}{q}\Z^2$.
     \item Let $i=1,2$ or $3$, and let $F\in \Delta \P_q$ be such that $p_i(F) \in \Ipoint$. Then, $F \subseteq E(\psiPoint)$.
  \item $\psiPoint$ is continuous piecewise linear over $\P_{mq}$.
  
  \end{enumerate}
\end{lemma}
\begin{proof}
  Properties (i), (ii), (iii) follow directly from the equivariance formula~\eqref{eq:equivariance}.   The function is continuous because it is continuous on the interior of each $I \in \Itri$ by construction and because $\psiPoint\big|_I = 0$ on all edges $I \in \Iedge$.   Property (iv) follows from properties (i), (ii), and (iii) and the fact that $\verts(\P_q) = \tfrac{1}{q}\Z^2$.  Finally, the function is continuous piecewise linear by construction as well.
\end{proof}

We next analyze $\varphiD^m$ from \autoref{s:eq-perturb}.  Let
$\Gamma_{q, \diag} = \langle \,\rho_\y, \tau_\y \st \y \in \R^2, \ve1\cdot\y \equiv 0
\pmod{\frac{1}{q}}\,\rangle \supseteq \Gamma_{q, \point}$  be the group generated by reflections and
translations corresponding to all points on diagonal edges of $\P_q$; see \autoref{fig:2d-perturbation-funda-domains} (right). 
In this case, $Y_{q, \diag} = \Lambda_{q, \diag} + L_\diag$ where $\Lambda_{q, \diag} = \tfrac{1}{q}\Z \times \{\0\}$ and $L_\diag$ is as defined in Definition~\ref{def:Ldiag}.  We choose the lattice basis $\bigl\{\t^1 = \ColVec[q]{1}{0}\bigr\}$, which has the fundamental parallelepiped $V_{q,\diag}= \bigl[\ColVec{0}{0}, \ColVec[q]{1}{0}\bigr]$ and hence 
$V^+_{q,\diag} = \bigl[\ColVec{0}{0}, \ColVec[2q]{1}{0}\bigr]$. (Note that we have simplified the notation, e.g., $V_{\Lambda_{q, \diag}}$ is now denoted by $V_{q, \diag}$.)

We consider an alternative description for the function $\varphiD^m$, $m\geq
3$.  This is done by setting $\varphiD^m\bigl(\ColVec{0}{0}\bigr) = 0$,
$\varphiD^m\bigl(\ColVec[2q]{1}{0}\bigr) = 0$,  and for integer $i$ with $1
\leq i < \tfrac{m}{2}$  we set $\varphiD^m\bigl(\ColVec[mq]{i}{0}\bigr) = 1$.
Then the function is interpolated over the vertices of $\P_{mq}$ that lie in
$V^+_{q,\diag}$. We extend the function to all of $\R^2$ by applying the
equivariance formula~\eqref{eq:equivariance} (the extension is well-defined by
Lemma~\ref{lemma:phi}).  This results in the continuous piecewise linear
function $\varphiD^m$ defined in~\autoref{s:eq-perturb}; refer back to
\autoref{fig:2d-perturbation-functions} (right) for an illustration.

\begin{lemma}\label{lemma:our-equivariant-varphiD}
  The function $\varphiD^m\colon\R^2\to\R$ has the following properties:
  \begin{enumerate}[\rm(i)]
  \item  $\varphiD^m\big|_I = 0$ on all edges and vertices $I \in \Ipointdiag$.
  \item $\varphiD^m(\x)= - \varphiD^m(\rho_{\y}(\x)) = - \varphiD^m(\y - \x)$ for all $\x \in \R^2$ and $\y \in \R^2$ such that $\ve 1 \cdot \y \equiv 0 \pmod{\tfrac1q}$.
  \item $\varphiD^m(\x) = \varphiD^m(\tau_\y(\x)) =\varphiD^m(\y + \x)$  for all $\x \in
    \R^2$ and $\y \in \R^2$ such that $\ve 1 \cdot \y \equiv 0 \pmod{\tfrac1q}$.
    \item Let $i=1,2,$ or $3$ and let $F \in \Delta \P_q$ be such that $p_i(F) \in \Ipointdiag$. Then, $F \subseteq E(\varphiD^m)$.
  \item $\varphiD^m$ is continuous piecewise linear over $\P_{mq}$.
  \end{enumerate}
\end{lemma}
\begin{proof}
  Properties (i), (ii), (iii) follow directly from the equivariance formula~\eqref{eq:equivariance}.   Property (iv) follows from properties (i), (ii), and (iii) and the fact that all faces of $\Ipointdiag$ are contained in the set $\{\,\y \in \R^2 \st \ve 1 \cdot \y \equiv 0 \pmod{\tfrac{1}{q}}\,\}$.  The function is continuous because the restriction to $V^+_{q,\diag}$ is continuous and the function vanishes on the relative boundary of $V^+_{q,\diag}$.   Finally, the function is piecewise linear by construction as well.
\end{proof}

\section{Genuinely $k$-dimensional functions}
\label{sec:gen-k-functions}
\subsection{Preliminaries} In this section, we prove useful properties of a special class of functions called {\em genuinely k-dimensional} functions. In the process, we motivate our assumption in Theorem~\ref{thm:main} and Theorem~\ref{thm:1/4q} that $\f \in \verts(\P_q)$.  

\begin{definition}\label{def:genk}
A function $\theta\colon \R^k \rightarrow \R$ is {\em genuinely $k$-dimensional}
if there does not exist a function $\varphi \colon \R^{k-1} \rightarrow \R$ and a linear
map $T \colon \R^k \rightarrow \R^{k-1}$ such that $\theta = \varphi\circ T$. 
\end{definition}

Genuinely $k$-dimensional functions 
were studied in \cite{bhkm}. 
We will show that $\f$ must be a vertex of the complex $\P$ whenever  $\pi$ is a minimal piecewise linear function over $\P$ that is genuinely $k$-dimensional.  We will then show that it suffices to consider only genuinely $k$-dimensional functions. This is because if the function is not genuinely $k$-dimensional we can study the function in a lower dimension by instead studying its restriction to a linear subspace of $\R^k$.

We will need the following lemma, which is implied by Lemma 13 in~\cite{bccz} and is a consequence of Dirichlet's Approximation Theorem for the reals. 

\begin{lemma}[\cite{bhkm}]
\label{lem:half-line}
Let $\y\in \R^k$ be any point and $\rr \in \R^k\setminus \{\0\}$ be any direction. Then for every $\epsilon > 0$ and $\bar\lambda \geq 0$, there exists $\w \in \Z^k$ such that $\y + \w$ is at distance less than $\epsilon$ from the half line $\{\,\y + \lambda \rr \st \lambda \geq \bar\lambda\,\}$.
\end{lemma}

The proof of the next lemma is adapted from the proof of Claim 2 in~\cite{bccz}. For any linear subspace $M$ of~$\R^k$, $\proj_M(\cdot)$ will denote orthogonal projection onto $M$. Also $M^\perp$ will denote the orthogonal complement of~$M$.

\begin{lemma}\label{lemma:lattice-projection}
Let $L$ be any linear subspace of $\R^k$. Then $\proj_L(\Z^k)$ has the following form: there exists a linear subspace $L' \subseteq L$ (we allow the possibility $L' = \{\0\}$) such that $\proj_L(\Z^k) = \Lambda + D$, where $\Lambda$ is a lattice that spans $L'^\perp \cap L$ and $D$ is a dense subset of $L'$.
\end{lemma}

\begin{proof} Let $\Lambda' = \proj_L(\Z^k)$. 
Let $V_\epsilon$ be the linear subspace of $L$ spanned by the points in $\{\, \y\in\Lambda'\st \|\y\|<\epsilon\,\}$. Notice that, given $\epsilon'>\epsilon''>0$, then $V_{\epsilon'}\supseteq V_{\epsilon''}\supseteq\{\0\}$. Since $\dim(V_\epsilon)$ changes discretely as $\epsilon \to 0$, there exists $\epsilon_0>0$ such that $V_\epsilon=V_{\epsilon_0}$ for every $0<\epsilon<\epsilon_0$. Let $L'=V_{\epsilon_0}$. Observe that $\Lambda'\cap L'$ is dense in $L'$ and $\Lambda = \proj_{L'^\perp\cap L}(\Lambda')$ is discrete (i.e., $B(\0, \epsilon_0) \cap \Lambda = \{\0\}$). Since $\Lambda$ is the projection of a subgroup of $\R^k$, it is also a subgroup and therefore it is a discrete subgroup, i.e., a lattice. We thus have the result using $D = \Lambda' \cap L'$.
\end{proof}

The following lemma can be found within the proof of Lemma 2.10 in \cite{bhkm} for the case where $L$ is a one-dimensional linear space. 
\begin{lemma}
\label{lem:linearSpace}
Suppose $\theta \colon \R^k \to \R$ is a subadditive function such that $\theta = 0$ on a linear space $L$.  For any $\x,\y \in \R^k$ such that $\x - \y \in L$, we have $\theta(\x) = \theta(\y)$.
\end{lemma}
\begin{proof}
Since $\x-\y \in L$, $\theta(\x-\y) = 0$.  By subadditivity, $\theta(\y) + \theta(\x-\y) \geq \theta(\x)$, which implies $\theta(\y) \geq \theta(\x)$. Similarly,  $\theta(\x) \geq \theta(\y)$, and hence we have equality. \end{proof}

The following lemma is modified version of Lemma 2.10 from \cite{bhkm} to give detail about when we can choose  a linear map $T$  that  can be represented as a rational matrix.  We assume Lipschitz continuity because this continuity is implicit in continuous piecewise linear functions.

\begin{lemma}
\label{lem:genuinely-k-dim}
Let $\theta \colon \R^k \rightarrow \R$ be nonnegative, Lipschitz continuous,
subadditive and periodic modulo the lattice $\Z^k$. 
Suppose there exist $\rr \in \R^k\setminus\{\0\}$ and $\bar\lambda >0$ such that
$\theta(\lambda \rr) = 0$ for all $0\leq \lambda \leq \bar\lambda$. Then $\theta$ is not genuinely $k$-dimensional, i.e., there exists a linear map $T\colon \R^k\to \R^{k-1}$ and a function $\varphi\colon \R^{k-1} \to \R$ such that $\pi = \varphi\circ T$. Furthermore, if $\rr \in \Q^k$, then $T$ can be represented by a rational matrix.
\end{lemma}

\begin{proof}
Let the Lipschitz constant for $\theta$ be $K$, that is, $|\theta(\x)-\theta(\y)| \leq K\|\x-\y\|$ for all $\x,\y\in \R^k$.

We will begin by showing that $\theta(\lambda \rr) = 0$ for all $\lambda \in \R$.  Let $\lambda' \in \R$.

Suppose that $\lambda' > \bar \lambda$ and let $M \in \Z_+$ such that $0 \leq \lambda'/M \leq \bar \lambda$.  From the hypothesis, we have that $\theta(\frac{\lambda'}{M} \rr)= 0$.  By nonnegativity and subadditivity of $\theta$ we see
$
0 \leq \theta(\lambda' \rr) \leq M \theta(\frac{\lambda'}{M} \rr)= 0,
$
and therefore, $\theta(\lambda' \rr) = 0$. This shows that $\theta(\lambda \rr) = 0$ for all $\lambda \geq 0$.

Next suppose $\lambda' < 0$.  By \autoref{lem:half-line}, for all $\epsilon >0$ there exists a $\w \in \Z^k$ such that $\lambda' \rr + \w$ is at distance less than $\epsilon$ from the half line $\{\,\lambda' \rr + \lambda \rr\st \lambda \geq -\lambda'\,\} = \{\lambda \rr \st \lambda \geq 0\}$.  That is, there exists a $\tilde\lambda \geq 0$ such that $\|\lambda' \rr + \w - \tilde \lambda \rr\| \leq \epsilon$. Since $\theta(\tilde \lambda \rr) = 0$, by periodicity and then Lipschitz continuity, we see that
$
0 \leq \theta(\lambda' \rr) = \theta(\lambda'\rr + \w) = \theta( \lambda' \rr + \w) - \theta(\tilde \lambda \rr) \leq K \epsilon.
$
This holds for every $\epsilon >0$ and therefore $\theta(\lambda'\rr) = 0$.  Thus, we have shown that $\theta(\lambda \rr) = 0$ for all $\lambda \in \R$.

Let $L = \{ \, \lambda \rr \st \lambda \in \R\,\}$.  By Lemma~\ref{lem:linearSpace}, for any $\x,\y$ such that $\x-\y \in L$, we have $\theta(\x) = \theta(\y)$.

We conclude that $\theta = \varphi \circ \proj_{L^\perp}$ for some function $\varphi \colon \R^{k-1} \rightarrow \R$ and therefore $\theta$ is not  genuinely $k$-dimensional.   Finally, if $\rr  \in \Q^k$, then $\proj_{L^\perp}$ can be represented by a rational matrix.
\end{proof}

\subsection{Dimension reduction for functions that are not genuinely $k$-dimensional}

We now show that it suffices to consider only genuinely $k$-dimensional functions for testing extremality of continuous piecewise linear functions. 


\begin{remark}
  Given a piecewise linear continuous valid function~$\zeta\colon \R\to\R$ for
  the one-dimensional infinite group problem $R_\f(\R,\Z)$,
  Dey--Richard \cite[Construction 6.1]{dey3} consider the
  function~$\kappa\colon 
  \R^2\to\R$, $\kappa(\ve x) =\zeta(\ve1\cdot\ve x)$, where $\ve1=\ColVec{1}{1}$, and show that $\kappa$ is
  minimal and extreme if and only if $\zeta$ is minimal and extreme,
  respectively. 
  If $\zeta$ has rational breakpoints in $\frac1q\Z$ with $q \in \Z_+$, then $\kappa$ belongs to
  our class of diagonally constrained continuous piecewise linear functions
  over $\P_q$.  However, these functions are not genuinely $2$-dimensional,
  and as Dey--Richard point out, we can study the one-dimensional function
  $\zeta$ instead of the $2$-dimensional function $\kappa$. We call the function~$\kappa$ a {\em diagonal embedding} of $\zeta$.

\end{remark}

The following two theorems can be found in \cite{dey3} for the special case of diagonal embeddings.   We also refer the interested reader to \cite{dey2} where the authors exhibit a sequential merge procedure, creating extreme functions in higher dimensions from extreme functions in lower dimensions and vice versa.

\begin{lemma}
\label{lem:mapMinimaliff}
Let $T\colon \R^k \to \R^\ell$ be a linear map. 
Suppose  $\pi \colon \R^k \to \R$ and $\varphi \colon T\R^k \rightarrow \R$ satisfy $\pi = \varphi \circ T$.  Then $\pi$ is minimal for $R_\f(\R^k, \Z^k)$ if and only if $\varphi$ is minimal for $R_{T\f}(T\R^k, T\Z^k)$.  
\end{lemma}
\begin{proof} ($\Longleftarrow$) Suppose $\varphi$ is minimal for $R_{T\f}(T\R^k, T\Z^k)$. We demonstrate that $\pi$ satisfies the criterion from \autoref{thm:minimal} to be minimal.

\begin{enumerate}
\item For any $\z \in \Z^k$, $0 = \varphi(T\z) = \pi(\z)$.
\item For any $\x,\y \in \R^{k}$ we have 
$$\pi(\x) + \pi(\y) - \pi(\x + \y) = \varphi(T\x) + \varphi(T\y) - \varphi(T(\x + \y))
= \varphi(T\x) + \varphi(T\y) - \varphi(T\x + T\y) \geq 0.$$
\item For any $\x\in \R^k$, we have 
$$
\pi(\x) + \pi(\f-\x) = \varphi(T\x) + \varphi(T(\f - \x)) =  \varphi(T\x) + \varphi(T\f - T\x) = 1.
 $$ 
\end{enumerate}

Therefore $\pi$ is minimal by \autoref{thm:minimal}.\smallbreak

\noindent ($\Longrightarrow$) Suppose $\pi$ is minimal for $R_\f(\R^k, \Z^k)$. We demonstrate that $\varphi$ satisfies the criterion from \autoref{thm:minimal} to be minimal.

\begin{enumerate}
\item For any $\z \in \Z^k$, $0 = \pi(\z) = \varphi(T\z)$.
\item For any $\x,\y \in T\R^k$, let $\hat \x \in T^{-1}\x$, $\hat \y \in T^{-1}\y$. Then 
$$
0 \leq \pi(\hat \x) + \pi(\hat \y) - \pi(\hat \x + \hat \y) = 
 \varphi(\x) + \varphi(\y) - \varphi( \x + \y).
$$
\item Similarly, for any $\x \in T\R^k$, let $\hat \x \in T^{-1} \x$.
Then $$
1 = \pi(\hat \x) + \pi(\f -\hat \x) = \varphi(\x) + \varphi(T\f - \x).
$$
\end{enumerate}

Therefore $\varphi$ is minimal by \autoref{thm:minimal}.\end{proof}

\begin{lemma}
\label{lemma:non-gen-k-test}
Let $\pi \colon \R^k \to \R$ be a minimal valid function.  Let $T \colon \R^k \rightarrow \R^{\ell}$ be a linear map and let $\varphi \colon T\R^k \rightarrow \R$
such that $\pi = \varphi\circ T$.  Then $\pi$ is extreme for $R_\f(\R^k,\Z^k)$ if and only if 
$\varphi$ is extreme for $R_{T\f}(T\R^k, T\Z^k)$.
\end{lemma}

\begin{proof}
%

($\Longrightarrow$) We prove the contrapositive. Suppose $\varphi$ is not extreme for $R_{T\f}(T\R^k, T\Z^k)$.  Then, by Lemma~\ref{lem:tightness}, there exist distinct minimal valid functions $\varphi^1, \varphi^2$ for $R_{T\f}(T\R^k, T\Z^k)$ such that $\varphi = \tfrac12( \varphi^1 + \varphi^2)$. But then $\pi^1 = \varphi^1 \circ T$ and $\pi^2 = \varphi^2 \circ T$ are distinct functions, and $\pi = \tfrac12( \pi^1 + \pi^2)$.  By \autoref{lem:mapMinimaliff}, $\pi^1, \pi^2$  are minimal for $R_\f(\R^k, \Z^k)$.  Therefore $\pi$ is not extreme.

($\Longleftarrow$) We again prove the contrapositive. Suppose that $\pi$ is not extreme for $R_\f(\R^k, \Z^k)$.  Then there exist distinct minimal valid functions $\pi^1, \pi^2$ for $R_\f(\R^k, \Z^k)$ such that $\pi = \tfrac12 (\pi^1 + \pi^2)$.  Since $\pi, \pi^1, \pi^2$ are minimal by Lemma~\ref{lem:tightness}, $\pi(\0)= \pi^1(\0) = \pi^2(\0) = 0$.  Since $E(\pi) \subseteq E(\pi^1), E(\pi^2)$ by \autoref{lemma:tight-implies-tight}, and $0 = \pi(\x) + \pi(-\x) - \pi(\0) = \Delta\pi(\x,-\x)$ for all $\x \in T^{-1}(\0)$, it follows that $\pi^i(\x) = -\pi^i(-\x)$ for $i=1,2$.   Since $\pi^i$ are valid functions, $\pi^i \geq 0$, therefore we must have $\pi^i(\x) =0 $ for all $\x \in T^{-1}(\0)$.   By Lemma~\ref{lem:linearSpace}, $\pi^i(\x) = \pi^i(\y)$ whenever $\x - \y \in T^{-1}(\0)$.  Therefore, we must have $\varphi^1, \varphi^2$ such that 
$\pi^1 = \varphi^1 \circ T$ and $\pi^2 = \varphi^2\circ T$.  Since $\pi^1, \pi^2$ are distinct, the functions $\varphi^1, \varphi^2$ are distinct as well.  Also since $\pi = \tfrac12( \pi^1 + \pi^2)$, we have $\varphi = \tfrac12( \varphi^1 + \varphi^2)$.  By \autoref{lem:mapMinimaliff}, $\varphi^1, \varphi^2$ are minimal for $R_{T\f}(T\R^k, T \Z^k)$.  Therefore $\varphi$ is not extreme.
\end{proof}

Given any family of polyhedra $\mathcal{F}$ (not necessarily a polyhedral complex), we say a polyhedral complex $\P$ is a {\em refinement} of $\mathcal{F}$ if every polyhedron of $\mathcal{F}$ is a union of polyhedra from $\P$.

\begin{prop}[Dimension reduction]\label{prop:dim-reduction}
\citedinsurveyas{Proposition B.9}
Let $\P$ be a pure and complete polyhedral complex in $\R^k$ that is periodic modulo $\Z^k$. Let $\pi\colon\R^k \to \R$ be a piecewise linear function over $\P$, such that $\pi$ is nonnegative, subadditive, periodic modulo $\Z^k$ and $\pi(\0) = 0$. If $\pi$ is not genuinely $k$-dimensional, then there exists a natural number $0 \leq \ell < k$, a pure and complete polyhedral complex $\mathcal{X}$ in $\R^\ell$ that is periodic modulo $\Z^\ell$, a nonnegative and subadditive function $\phi\colon\R^\ell \to \R$ that is piecewise linear over $\mathcal{X}$, and a point $\f' \in \R^\ell \setminus \Z^\ell$ with the following properties:
\begin{enumerate}
\item $\pi$ is minimal for $R_\f(\R^k, \Z^k)$ if and only if $\phi$ is minimal for $R_{\f'}(\R^\ell, \Z^\ell)$. 
\item $\pi$ is extreme for $R_\f(\R^k, \Z^k)$ if and only if $\phi$ is extreme for $R_{\f'}(\R^\ell, \Z^\ell)$.
\end{enumerate}
\end{prop}

\begin{proof} Since $\pi$ is not genuinely $k$-dimensional, it follows by iteratively
  applying the definition of genuinely $k$-dimensional functions that there
  exist a number $0 \leq \ell < k$, a function $\varphi\colon \R^\ell \to \R$, and a linear map $T\colon \R^k \to \R^\ell$ such that $\varphi\colon \R^\ell \to \R$ is genuinely $\ell$-dimensional and $\pi = \varphi \circ T$. Since $\pi$ is nonnegative, $\varphi$ must also be nonnegative. Since $\pi$ is subadditive and $T$ is additive, $\varphi$ must be subadditive.
\bigskip

{\em Claim 1. $T\Z^k$ is a lattice that spans $\R^\ell$.}\medskip

Since every linear map is a projection composed with an isomorphism, Lemma~\ref{lemma:lattice-projection} implies that there exists a linear subspace $L \subseteq \R^\ell$ such that $T\Z^k = \Lambda + D$, where $\Lambda$ is a lattice spanning $L^\perp$ and $D$ is dense in $L$. If $L = \{\0\}$ then we are done. So we assume $\dim(L) \geq 1$. Since $\pi$ is continuous (it is piecewise linear over a locally finite polyhedral complex),  and $T$ is linear map, it follows that $\varphi$ is continuous. Also, since $\pi$ vanishes over $\Z^k$, $\varphi$ vanishes over $T\Z^k$. But this implies that $\varphi$ vanishes over $D$, and thus over $L$. By Lemma~\ref{lem:linearSpace}, $\varphi$ is constant on the affine subspaces parallel to $L$. This contradicts the assumption that $\varphi$ is genuinely $\ell$-dimensional. This concludes the proof of Claim 1.
\bigskip

Let $\mathcal{U} = \bigcup_{I \in \P} \{I \cap [0,1]^n\}$. Since $\pi$ is piecewise linear over $\P$, $\pi$ is also piecewise linear over a refinement of $\P$, in particular, over the polyhedral complex $\bigcup_{I \in \mathcal{U}, \w \in \Z^k}\{ I + \w\}$, that is periodic modulo $\Z^k$. Since $T\Z^k$ is a lattice and for every $I \in \mathcal{U}$, $TI$ is a polytope (it is the projection of the polytope $I$), we can find a refinement of the family of polytopes $\bigcup_{I \in \mathcal{U}, \w \in T\Z^k} \{TI + \w\}$; we denote this refinement by $\P'$, which is a pure and complete polyhedral complex of $\R^\ell$. We observe that $\varphi$ is piecewise linear over $\P'$ and $\P'$ is a polyhedral complex that is periodic modulo $T\Z^k$. 

Now simply find an invertible linear transformation $A\colon T\Z^k \to \Z^\ell$ and let $\phi := \varphi\circ A^{-1}$ be the piecewise linear function defined over the pure and complete polyhedral complex $\mathcal{X}:= A\P'$ and let $\f' := AT\f$. Then $\f' \not\in \Z^\ell$, since $1 = \pi(\f) = \phi(\f')$ and $\phi(\Z^\ell) = \pi(\Z^k) = 0$. The two properties now follow from Lemmas~\ref{lem:mapMinimaliff} and \ref{lemma:non-gen-k-test}.
\end{proof}

\begin{remark}[Dimension reduction]
\label{remark:dimension-reduction}
\citedinsurveyas{Remark B.10}
Using Proposition~\ref{prop:dim-reduction}, the extremality/minimality question for $\pi$ that is not genuinely $k$-dimensional can be reduced to the same question for a lower-dimensional genuinely $\ell$-dimensional function (so $\ell < k$.) When $\P$ is a rational polyhedral complex, this reduction can be done algorithmically. The question of making this effective for the irrational case is beyond the scope of this paper.

\end{remark}

\subsection{The assumption of $\f \in \verts(\P)$}

We will show that $\f$ is a vertex for any minimal continuous piecewise linear function that is genuinely $k$-dimensional.
\begin{theorem}
\label{thm:faVertex-gen-k-dim}
\citedinsurveyas{Theorem B.11}
Let $\P$ be a pure and complete polyhedral complex in $\R^k$ that is periodic modulo $\Z^k$. Let $\theta \colon \R^k \rightarrow \R$ be minimal, piecewise linear function over $\P$ that is genuinely $k$-dimensional.  Then $\f \in \verts(\P)$.
\end{theorem}

\begin{proof}
For the sake of contradiction, suppose $\f \notin \verts(\P)$.  Therefore, there exists some $I \in \P$ with $\f \in \relint(I)$ and the dimension of $I$ is at least one.   Since $\pi$ is minimal, 
$0 \leq \pi \leq 1$.  Since $\pi(\f) = 1$, $\pi \leq 1$, $\pi$ is affine on $I$ and $\f \in \relint(I)$, we have that $\pi(\x) = 1$ for all $\x \in I$.
Now consider $\pi$ on $\f - I$ and note that $\0 \in \f - I$.  By symmetry,
$\pi(\x) = 0 $ for all $\x \in \f-I$.  Since the dimension of $I$ is at least
one, there exists $\rr \in (\f - I)\setminus \{\0\}$.  But then $\pi(\lambda\rr) = 0$ for all $\lambda \in [0,1]$.  Since $\pi$ is continuous piecewise linear over $\P$, by \autoref{lem:lipschitz}, it satisfies the hypotheses of  \autoref{lem:genuinely-k-dim}.
 Therefore, $\pi$ is not genuinely $k$-dimensional, which is a contradiction.  Therefore, we must have $\f \in \verts(\P)$.
\end{proof}

\begin{remark}
Using Proposition~\ref{prop:dim-reduction} and Theorem~\ref{thm:faVertex-gen-k-dim}, we can achieve dimension reduction when $\f \not\in \verts(\P)$. Thus, although the results presented in this paper assume that $\f\in \verts(\P)$, this assumption is actually not very restrictive.
\end{remark}

\subsection{Boundedness of cells for genuinely $k$-dimensional functions}\label{sec:bounded-cells}
In this subsection, we show that for genuinely $k$-dimensional minimal valid functions that are piecewise linear over a pure and complete polyhedral complex $\P$ in $\R^k$ that is periodic modulo $\Z^k$, the cells of $\P$ are full-dimensional bounded polytopes (so they cannot be unbounded polyhedra).
 
\begin{lemma}\label{lemma:closed-proj}
Let $\rr \in \R^k$ be any vector and let $L = \rr^\perp$ be the orthogonal complement of the subspace spanned by $\rr$. Let $U$ be a compact convex set with nonempty interior in $\R^k$. Then $\proj_L(U + \Z^k)$ is a closed set.
\end{lemma}

\begin{proof} Since orthogonal projections onto linear subspaces are linear operators, $\proj_L(U + \Z^k) = \proj_L(U) + \proj_L(\Z^k)$. Observe that $\proj_L(U)$ is also a compact convex set with nonempty interior with respect to~$L$. By \autoref{lemma:lattice-projection}, there exists a linear subspace $L' \subseteq L$ such that $\proj_L(\Z^k) = \Lambda + D$, where $\Lambda$ is a lattice that spans $L'^\perp \cap L$ and $D$ is a dense subset of $L'$. Since $\proj_L(U)$ is convex with nonempty interior, $\proj_L(U) + D = \proj_L(U) + L'$. Let $U'$ be the orthogonal projection of $\proj_L(U)$ onto $L'^\perp \cap L$; so $U'$ is compact convex set. Thus, we have 
$$\begin{array}{rcl}\proj_L(U + \Z^k) & = & \proj_L(U) + \proj_L(\Z^k) \\
& = & \proj_L(U) + \Lambda + D \\
& = & \proj_L(U) + L' + \Lambda \\
& = & U' + L' + \Lambda.
\end{array}$$

Since $U'$ is a compact set and $\Lambda$ is a closed set, $U' + \Lambda$ is closed  (see, e.g.,~\cite{inf-dim-an}  Lemma 5.3 (4)). Moreover, $U' + \Lambda \subseteq L'^\perp \cap L$. Therefore, $U' + \Lambda + L'$ is closed.\end{proof}

Let $H := [0,1]^k$ denote the unit hypercube.

\begin{lemma}\label{lemma:finite-rep} Let $\P$ be a locally finite polyhedral complex that is periodic modulo $\Z^k$. Then for any full-dimensional polyhedron $I \in \P$, the set $I + \Z^k$ is a finite union of the form $\bigcup_{j\in J}(I_j + \Z^k)$ where $J$ is a finite index set and each $I_j$ is a full-dimensional polytope contained in $H$. 
\end{lemma}

\begin{proof}We can take $I_j$ to be all full-dimensional polytopes contained in $(I + \Z^k) \cap H$. There are only finitely many of these polytopes by the locally finite property of $\P$ (see Definition~\ref{def:polyhedralComplex}\,(iv)).\end{proof}

  \begin{lemma}
  \label{lemma:genKvertices}
  Let $\theta \colon \R^k \rightarrow \R$ be a piecewise linear minimal valid function over a polyhedral complex $\P$ that is pure, complete and periodic modulo the lattice $\Z^k$.  If $\theta$ is genuinely $k$-dimensional, then the cells of $\P$ and $\Delta \P$ are full-dimensional polytopes. 
  \end{lemma}

\begin{proof}


Suppose to the contrary that a cell $I^*$ has a recession direction $\rr$. Let $L$ be the linear subspace orthogonal to $\rr$, i.e., $L = \langle \rr\rangle^\perp$. 
Let $U = \bigcup\{\,I \in \P\st \rr \textrm{ is a recession direction for }I\,\}$. Define $S = \proj_L(U)$.
\bigskip

{\em Claim 1. $S = L$.} \medskip

First, notice that $H \cap \P$ contains finitely many full-dimensional polytopes by the local finiteness of $\P$. Combining this observation with \autoref{lemma:finite-rep}, we can express $U = \bigcup_{j \in J}(I_j + \Z^k)$ where $J$ is a finite index set and each $I_j$ is a full-dimensional polytope. Therefore, $S = \proj_L(U) = \bigcup_{j \in J}\proj_L(I_j + \Z^k)$, which is a finite union of closed sets by \autoref{lemma:closed-proj}. Therefore, $S$ is closed. The set $S$ is nonempty because $I^*$ has recession direction $\rr$. If $S \neq L$, then there exists a boundary point $\x$ of $S$ (considered as a subset of $L$). Thus, there exist $Q_0\in \P$ and $\y \in Q_0$ such that $\x = \proj_L(\y)$ and $Q_0$ has $\rr$ as a recession direction. Moreover, we can choose $\y$ so that $\y$ is in the relative interior of a face $F_0 \subseteq Q_0$ where $F_0$ also has $\rr$ as a recession direction. Let $Q_1, \ldots, Q_p\in \P$ be the cells that also have $F_0$ as their face (using the local finiteness of $\P$). We set $p=0$ if $F_0 = Q_0$. Since $\P$ is complete and $\y \in \relint(F_0)$, we can choose $\delta > 0$ such that $B(\y, \delta) \subseteq Q_0 \cup Q_1 \cup \dots \cup Q_p$. Since $F_0$ is a face of each of these polyhedra, $\rr$ is a recession direction for each $Q_0, Q_1, \ldots, Q_p$. Thus, $B(\y,\delta) \subseteq U$ and thus, $\proj_L(B(\y,\delta)) \subseteq S$. But $\proj_L(\y) = \x$ and $\x$ is a boundary point of $S$. This is a contradiction. Therefore, $S=L$. This concludes the proof of Claim 1.
\bigskip

Fix $\x \in L = S$. Let $Q\in \P$ be the cell such that $\x \in \proj_L(Q)$ and $\rr$ is a recession direction of $Q$. Thus, there exists a constant $\lambda(\x) \in \R$ such that $\x + \mu\rr \in Q$ for all $\mu \geq \lambda(\x)$. Since $\theta$ is bounded and affine over $Q$, $\theta$ must be constant on the half-line $\x + \mu\rr$, $\mu \geq \lambda(\x)$. Thus, there exists a constant $C(\x)$ such that $\theta(\x + \mu \rr) = C(\x)$ for all $\mu \geq \lambda(\x)$. 
We now show that $\theta(\x + \mu \rr) = C(\x)$ for all $\mu \in \R$. Let $\mu' < \lambda(\x)$ and let $\y = \x + \mu'\rr$. By \autoref{lem:half-line}, for all $\epsilon >0$ there exists $\w \in \Z^k$ such that $\y + \w$ is at distance less than $\epsilon$ from the half line $\{\,\y +\mu\rr \st \mu \geq \lambda(\x) - \mu'\,\} = \{\,\x + \lambda \rr \st \lambda \geq \lambda(\x)\,\}$.  That is, there exists $\tilde\lambda \geq \lambda(\x)$ such that $\|\y + \w - (\x + \tilde \lambda \rr)\| \leq \epsilon$. Since $\theta(\x + \tilde \lambda \rr) = \theta(\x + \lambda(\x)\rr)$, by periodicity and then Lipschitz continuity, we see that
$
|\theta(\x + \lambda(\x)\rr) - \theta(\y)| = |\theta(\x + \tilde \lambda \rr) - \theta(\y+\w)|\leq K \epsilon.
$
This holds for every $\epsilon >0$ and therefore $\theta(\y)=\theta(\x + \lambda(\x)\rr)$.  Thus, we have shown that for any $\x \in L$, $\theta$ is constant on the line $\x + \mu \rr$, $\mu \in \R$. But this contradicts the fact that $\theta$ is genuinely $k$-dimensional.

Finally, if all cells of $\P$ are polytopes, then this property also holds for $\Delta\P$.\end{proof}

\section{Additional proofs}
\label{s:additional-proofs}

\subsection{Proofs of lemmas on polyhedral complexes}
\label{s:additional-proofs:polyhedral-complexes}

\begin{proof}[Proof of \autoref{prop:projection}]
  First of all, we have
\begin{align*}
p_1(F(I,J,K)) 
&= \{\, \x \in I \st \exists \y \in J, \z \in K \text{ such that } \x + \y = \z \,\}\\
&= \{\, \x \in \R^k \st \exists \y \in J, \z \in K \text{ such that } \x + \y = \z\,\} \cap I\\
&= \{\, \z - \y \st  \y \in J, \z \in K  \,\} \cap I\\
&= (K + (-J)) \cap I.
\end{align*}
A similar calculation shows $p_2(F(I,J,K)) = (K+ (-I)) \cap J$.  Finally,
\begin{align*}
p_3(F(I,J,K)) 
&= \{\, \z \in K \st \exists \x \in I, \y \in J \text{ such that } \x + \y = \z \,\}\\
&= \{\, \z \in \R^k \st \exists  \x \in I, \y \in J \text{ such that } \x + \y = \z \,\} \cap K\\
&= \{\, \x +\y \st  \x \in I, \y \in J \,\} \cap K\\
&= (I+J) \cap K. \tag*{\qedhere}
\end{align*}
\end{proof}

\begin{proof}[Proof of \autoref{lemma:delta-p-is-complex}]
We show the 4 conditions of Definition \ref{def:polyhedralComplex}. 
\begin{enumerate}[\rm(i)]
\item Since $\emptyset \in \P$, we have $F(\emptyset, \emptyset, \emptyset) = \emptyset \in \Delta \P$.
\item Let $I,J,K\in \P$.  Let $\hat F$ be a face of $F(I,J,K)$.  Write $I,J,K$ as inequality systems as $A_I \x \leq \b_I$, $A_J \x \leq \b_J$, $A_K\x \leq \b_K$.  Then 
$$
 F(I,J,K) = \{\,(\x,\y) \st A_I \x \leq \b_I,\ A_J \y \leq \b_J,\ A_K (\x + \y) \leq \b_K\,\}.
$$
The face $\hat F$ is obtained by setting certain inequalities to equalities.  This corresponds to restricting to faces of $I,J,K$.  Therefore, there exist $I', J',K'  \in \P$ such that $F(I',J',K') = \hat F$.  Therefore $\hat F \in \Delta \P$.
\item Let $I,J,K,I', J',K'\in \P$.  Then $F(I,J,K) \cap F(I',J',K') = F(I \cap I', J\cap J', K \cap K')$.  Since $\P$ is closed under intersection, $I \cap I', J\cap J', K \cap K' \in \P$.  Therefore $F(I \cap I', J\cap J', K \cap K') \in \Delta \P$.
\item Since $\P$ is locally finite, it follows that $\Delta \P$ is locally finite. 
\end{enumerate}
Hence, $\Delta \P$ is a polyhedral complex.  Finally, consider any $(\x,\y) \in \R^k \times \R^k$.  Let $I, J,K \in \P$ such that $\x \in I$, $\y \in J$, $\x+\y \in K$.  These faces $I,J,K$ exist since $\P$ is complete in $\R^k$.  Therefore, $(\x,\y) \in F(I,J,K) \in \Delta \P$.  Thus, $\Delta \P$ is complete.   Since it is a locally finite complete polyhedral complex, it is also pure. This follows from the following argument. Suppose to the contrary, there is a face $F$ in $\Delta\P$ that is maximal but not full-dimensional. Let $(\x,\y) \in F$ be a point in the relative interior of $F$ and note that $(\x,\y)$ cannot be contained in any other face of $\Delta \P$ by the maximality of $F$. By the locally finite property of $\Delta \P$, there exists an open ball $B$ around $(\x,\y)$ such that $B$ intersects $\Delta \P$ only in $F$. Since $B\cap F$ is a strict subset of $B$, this contradicts that $\Delta \P$ is complete.
\end{proof}


\clearpage
\providecommand\ISBN{ISBN }
\bibliographystyle{../amsabbrvurl}
\bibliography{../bib/MLFCB_bib}

\end{document}